\documentclass[review]{elsarticle}

\usepackage{amsmath}
\usepackage{amsfonts}
\usepackage{amssymb}
\usepackage{amsthm}
\usepackage{bm}
\usepackage{indentfirst}
\usepackage{graphicx}
\usepackage{subfigure}
\usepackage{array}
\usepackage{fullpage}

\DeclareMathAlphabet{\mathsfsl}{OT1}{cmss}{m}{sl}

\newcommand{\PreserveBackslash}[1]{\let\temp=\\#1\let\\=\temp}
\newcolumntype{C}[1]{>{\PreserveBackslash\centering}p{#1}}
\newcolumntype{R}[1]{>{\PreserveBackslash\raggedleft}p{#1}}
\newcolumntype{L}[1]{>{\PreserveBackslash\raggedright}p{#1}}

\numberwithin{equation}{section}
\newtheorem{thm}{Theorem}[section]

\newtheorem{cor}{Corollary}
\theoremstyle{definition}

\usepackage{anyfontsize}
\usepackage{enumerate}
\usepackage{esint}
\usepackage{etoolbox}
\usepackage{isomath}
\usepackage{mathrsfs}
\usepackage{mathtools}
\usepackage{multicol}
\usepackage{pgfplots}
\usepackage{scalerel}
\usepackage{stackengine}
\usepackage{tabulary}
\usepackage{tikz}

\makeatletter
\newcommand*\bdot{\mathpalette\bdot@{.65}}
\newcommand*\bdot@[2]{\mathbin{\vcenter{\hbox{\scalebox{#2}{$\m@th#1\bullet$}}}}}
\makeatother

\makeatletter
\newcommand*\bddot{\mathpalette\bddot@{.65}}
\newcommand*\bddot@[2]{\mathbin{\vcenter{\hbox{\scalebox{#2}
    {$\m@th#1\smash{{}_{\bullet}^{\bullet}}$}}}}}
\makeatother

\newcommand{\circled}[2][]{%
  \tikz[baseline=(char.base)]{%
    \node[shape = circle, draw, inner sep = .5pt]
    (char) {\phantom{\ifblank{#1}{#2}{#1}}};%
    \node at (char.center) {\makebox[0pt][c]{#2}};}}
\robustify{\circled}

\makeatletter
\newcommand{\opnorm}{\@ifstar\@opnorms\@opnorm}
\newcommand{\@opnorms}[1]{%
  \left|\mkern-1.5mu\left|\mkern-1.5mu\left|
   #1
  \right|\mkern-1.5mu\right|\mkern-1.5mu\right|
}
\newcommand{\@opnorm}[2][]{%
  \mathopen{#1|\mkern-1.5mu#1|\mkern-1.5mu#1|}
  #2
  \mathclose{#1|\mkern-1.5mu#1|\mkern-1.5mu#1|}
}
\makeatother

\stackMath
\newcommand\reallywidecheck[1]{%
\savestack{\tmpbox}{\stretchto{%
  \scaleto{%
    \scalerel*[\widthof{\ensuremath{#1}}]{\kern-.6pt\bigwedge\kern-.6pt}%
    {\rule[-\textheight/2]{1ex}{\textheight}}%WIDTH-LIMITED BIG WEDGE
  }{\textheight}%
}{0.5ex}}%
\stackon[1pt]{#1}{\scalebox{-1}{\tmpbox}}%
}
\usepackage[utf8]{inputenc}
\usepackage{bm}
\usepackage{mathrsfs}
\usepackage{amsmath}
\usepackage{amsfonts}
\usepackage{amsthm}
\usepackage{color}
\usepackage{setspace}
\usepackage{float}
\usepackage{import}
\usepackage{url}
\usepackage{multicol}
\usepackage{subfigure}
\biboptions{sort&compress}
\usepackage[top=1in,bottom=1in,left=1in,right=1in]{geometry}

\newtheorem{remark}{Remark}

\newcommand{\real}{\mathbb{R}}

\newcommand{\mcL}{\mathcal{L}}

\newcommand{\mcO}{\mathcal{O}}

\usepackage{xcolor}

\def\omg{{\Omega}}
\def\omgi{\mathcal{I}{\Omega}}
\def\omgb{\mathcal{B}\Omega}
\def\omgbb{\mathcal{B}\mathcal{B}\Omega}

\def \bb{\mathbf{b}}

\def \fb{\mathbf{f}}

\def \ub{\mathbf{u}}

\def \wb{\mathbf{w}}
\def \vb{\mathbf{v}}
\def \xb{\mathbf{x}}

\def \zb{\mathbf{z}}

\def \nb{\mathbf{n}}
\def \yb{\mathbf{y}}

\newcommand{\vertii}[1]{{\left\vert\left\vert #1
    \right\vert\right\vert}}    

\begin{document}

\begin{frontmatter}

\title{An asymptotically compatible treatment of traction loading in linearly elastic peridynamic fracture}

\address[yy]{Department of Mathematics, Lehigh University, Bethlehem, PA}
\address[nt]{Center for Computing Research, Sandia National Laboratories, Albuquerque, NM}

\author[yy]{Yue Yu\corref{cor1}}\ead{yuy214@lehigh.edu}
\author[yy]{Huaiqian You}\ead{huy316@lehigh.edu}
\author[nt]{Nathaniel Trask}\ead{natrask@sandia.gov}

%\date{January 2020}

\begin{abstract}
Meshfree discretizations of state-based peridynamic models are attractive due to their ability to naturally describe fracture of general materials. However, two factors conspire to prevent meshfree discretizations of state-based peridynamics from converging to corresponding local solutions as resolution is increased: quadrature error prevents an accurate prediction of bulk mechanics, and the lack of an explicit boundary representation presents challenges when applying traction loads. In this paper, we develop a reformulation of the linear peridynamic solid (LPS) model to address these shortcomings, using improved meshfree quadrature, a reformulation of the nonlocal dilitation, and a consistent handling of the nonlocal traction condition to construct a model with rigorous accuracy guarantees. In particular, these improvements are designed to enforce discrete consistency in the presence of evolving fractures, whose {\it a priori} unknown location render consistent treatment difficult.  In the absence of fracture, when a corresponding classical continuum mechanics model exists, our improvements provide asymptotically compatible convergence to corresponding local solutions, eliminating surface effects and issues with traction loading which have historically plagued peridynamic discretizations. When fracture occurs, our formulation automatically provides a sharp representation of the fracture surface by breaking bonds, avoiding the loss of mass. We provide rigorous error analysis and demonstrate convergence for a number of benchmarks, including manufactured solutions, free-surface, nonhomogeneous traction loading, and composite material problems. Finally, we validate simulations of brittle fracture against a recent experiment of dynamic crack branching in soda-lime glass, providing evidence that the scheme yields accurate predictions for practical engineering problems.

\end{abstract}


\end{frontmatter}

%\begin{document}
%\maketitle

%\tableofcontents

\section{Introduction}

Peridynamics provides a description of continuum mechanics in terms of integral operators rather than classical differential operators \cite{silling_2000,seleson2009peridynamics,parks2008implementing,zimmermann2005continuum,emmrich2007analysis,du2011mathematical,bobaru2016handbook}. These nonlocal models are defined in terms of a lengthscale $\delta$, referred to as a horizon, which denotes the extant of nonlocal interaction. The nonlocal viewpoint allows a natural description of processes requiring reduced regularity in the relevant solution, such as fracture mechanics \cite{bazant2002nonlocal,du2013nonlocal}. An important feature of such models is that { when classical continuum models still apply,} they revert back to classical continuum models as $\delta \rightarrow 0$. Discretizations which preserve this limit under refinement $h \rightarrow 0$ are termed asymptotically compatible (AC) \cite{tian2014asymptotically}, and there has been significant work in recent years toward establishing such discretizations - for an incomplete list see  \cite{tian2014asymptotically,d2020numerical,leng2019asymptotically,pasetto2018reproducing,hillman2020generalized,seleson2016convergence,du2016local,trask2019asymptotically,You_2019,you2020asymptotically,tao2017nonlocal}. Broadly, strategies either involve adopting traditional finite element shape functions and carefully performing geometric calculations to integrate over relevant horizon/element subdomains, or adopt a strong-form meshfree discretization where particles are associated with abstract measure. The former is more amenable to mathematical analysis due to a better variational setting, while the latter is simple to implement and generally faster \cite{silling2005meshfree,bessa2014meshfree}. In this paper we pursue the meshfree viewpoint.

For fracture mechanics problems one often refines both $\delta$ and $h$ at the same rate under so-called M-convergence, $\delta = M h$, for $M>0$ \cite{bobaru2009convergence}. In this setting, one obtains banded stiffness matrices allowing scalable implementations. Typically in the literature a scheme is termed AC if it recovers the solution in both the finite $\delta$ and M-convergence limit - in this work we abuse the definition slightly and only require the M-convergence case for asymptotic compatibility as the relevant limit for {problems with a corresponding local limit}. This AC property is only one necessary ingredient in achieving a convergent simulation, and our recent work focused upon establishing convergence in this setting for boundary value problems \cite{You_2019,you2020asymptotically}. To achieve similar convergence for problems involving fracture, one must also consider the interplay between consistency of quadrature for discrete operators and the imposition of traction loads as fracture surfaces open up and evolve \cite{lipton2014dynamic}. For peridynamic fracture problems where the free surface evolves implicitly via the breaking of bonds \cite{parks2008pdlammps,trask2019asymptotically}, one lacks an explicit boundary representation over the course of a simulation. In addition to providing challenges regarding accurate imposition of traction loads, the breaking of bonds also renders higher-order numerical quadrature inaccurate, as consistent AC quadrature weights are typically derived in the absence of damage.

Our goal is to provide a comprehensive treatment of fracture, nonlocal quadrature, and traction loading which is able to {perform more accurate state-based peridynamic fracture simulations free of spurious surface effects. In particular, when no fracture occurs and therefore the classical continuum theory applies, the formulation should preserve the AC limit under M-convergence. When fracture occurs, the formulation should be able to capture the material damage and the evolving fracture surfaces via bond breaking.} This practically means that one is able to incorporate all of the necessary ingredients to perform non-trivial simulations of fracture mechanics while maintaining a scalable implementation and guaranteeing convergence. Such a capability is elusive in the peridynamic literature; while peridynamics has been shown to provide a powerful modeling platform for a broad range of applications \cite{diehl2019review,javili2019peridynamics}, the development of efficient discretizations with rigorous underpinnings has lagged behind until the last few years.
%{\color{red}Differs from the traditional definition of asymptotic compatibility which only exists in the absence of damage, the convergence here in fact requires a correct discretization of the non-local model which captures aspects of existing local models. In particular, when the classical continuum theory applies, the numerical nonlocal solution should converge to the correct local limit. When fracture occurs, here is no known field theory that is approachedby peridynamics with bond breaking.  (It is noted that there is an emerging theory of localfracture modeling that is approached by peridynamic models with bond softening, see forreference [1], [2], [3], [4].)}

The challenge in incorporating traction loading into a peridynamic framework stems from the fact that, in contrast to local mechanics, peridynamic boundary conditions must be defined on a finite volume region outside the surface \cite{cortazar2008approximate,du2013nonlocal,tao2017nonlocal}. Theoretical and numerical challenges arise in how to mathematically impose nonhomogeneous Neumann boundary conditions properly in the nonlocal model. {In peridynamic models, careless imposition of traction loads leads to a smaller effective material stiffness close to the boundary, since the integral on those material points is over a smaller region. Therefore, an unphysical strain energy concentration is induced, leading in turn to an artificial softening of the material near the boundary.} Such undesirable phenomena are referred to in the literature as a ``surface'' or ``skin'' effects \cite{ha2011characteristics,bobaru2011adaptive}. We propose a novel treatment of nonlocal traction-type boundary conditions which avoid the surface effect by designing a loading aimed to recover the corresponding local traction boundary condition as $\delta \rightarrow 0$. The approach requires no explicit representation of the boundary, imposing the traction volumetrically using the same information that would be available during a traditional meshfree bond-based peridynamics simulation. Although the Neumann-constrained nonlocal problem and its AC limit were investigated in nonlocal diffusion models \cite{cortazar2008approximate,du2015integral,tao2017nonlocal,d2020physically,You_2019,you2020asymptotically}, to the authors' best knowledge, the development of AC peridynamic formulations with traction-type boundary conditions remains restricted to weak formulations, simple traction loadings and/or simple geometries. Several modeling and numerical approaches have been proposed to correct the surface effect \cite{le2018surface,bobaru2016handbook,madenci2014coupling,oterkus2010peridynamic,macek2007peridynamics,du2017peridynamic,madenci2014peridynamic,oterkus2015peridynamics} but mostly restricted to free surfaces. For nonzero loadings, the tractions are often applied as prescribed body forces through a layer of finite thickness at the material boundary \cite{le2018surface,javili2019peridynamics,lipton2019classic}, as a surface integral through a weak form \cite{madenci2018weak}, or by modifying the nonlocal operator through eigenvalues analysis \cite{aksoylu2020nonlocal}. 
%In contrast to variational discretizations of nonlocal operators, meshfree discretizations are often used due to their ease of implementation, efficiency, and natural treatment of fracture problems \cite{silling_2005_2}.
Therefore, developing an AC meshfree discretization method for peridynamics which is capable to handle nonhomogeneous traction loadings on complex boundaries is critical for the general practice of peridynamics in realistic engineering applications.

We consider the linear peridynamic solid (LPS) model \cite{emmrich2007well} as a prototypical state-based model appropriate for brittle fracture. The LPS model may be interpreted as a nonlocal generalization of the mixed form of linear elasticity, evolving both displacements and a dilitation. We will show that consistent treatment evolving traction loading will require a modification to the definition of dilitation to guarantee consistency in the presence of fractures; conceptually this corresponds to the fact that dilitation is a kinematic variable without associated boundary conditions, and should be estimated consistently independently of whether a fracture is occurring in the vicinity of a given point. Based on the modified nonlocal dilitation, we further propose a new nonlocal generalization of classical traction loads in the LPS model. Particularly, we convert the local traction loads to a correction term in the momentum balance equation, which provides an estimate for the nonlocal interactions of each material point with points outside the domain. Based on this traction-type boundary condition, a meshfree formulation is developed for the LPS model based on the optimization-based quadrature rule \cite{trask2019asymptotically}, which preserves the AC limit under M-convergence and naturally represents the evolving free surfaces in dynamic fracture problems. {We note that asymptotic compatibility is not well-defined for dynamic fracture, as there is no known corresponding local theory for peridynamics with bond breaking}\footnote{There is an emerging theory of local fracture modeling that is approached by peridynamic models with bond softening, see \cite{lipton2014dynamic,lipton2016cohesive,jha2020kinetic,lipton2020plane}.}. { However, our modified LPS formulation preserves the AC limit for the linear elastic model with traction loading on the evolving fracture surfaces. This fact, together with the consistent discretization introduced here, provide the opportunity for efficient and accurate peridynamic fracture simulations.}

We remark that the paper is organized to first establish the rigorous mathematical underpinnings of the approach, while the second half focuses on a more engineering-oriented exploration of its application. Readers with more applied interests may skip many of the proofs in the work without issue. The work is organized as follows. We recall first the linear peridynamic solid (LPS) model definition in Section \ref{sec:LPS}. In Section \ref{sec:Neumann}, we introduce a novel approach to apply classical traction loads on the LPS model. After establishing the continuous limits of the scheme, we next pursue a consistent discretization. In Section \ref{sec:disc} we introduce meshfree quadrature which preserves asymptotic compatibility in the $\delta\rightarrow0$ limit, and establish the discrete scheme for boundary value problems in the absence of fracture. We proceed to investigate a number of two-dimensional statics problems with analytic solutions for the local limit in Section \ref{sec:num}. These test cases include: linear patch tests (Subsection \ref{sec:patch}); manufactured local limits to illustrate asymptotic convergence rates (Subsection \ref{sec:mconv}); homogeneous materials with free-surfaces or non-zero traction loading on curvilinear surfaces (Subsection \ref{ssec:isoSol}); composite materials with internal interfaces (Subsection \ref{ssec:anisoSol}). In Section \ref{sec:exp}, we further extend the proposed formulation to handle dynamic brittle fracture, and provide preliminary validation results by comparing our numerical results with available numerical simulations and experimental measurements on three benchmark problems. Section \ref{sec:conclusion} summarizes our findings and discusses future research.

\section{A Linear State-Based Peridynamic Model}\label{sec:LPS}

We consider the state-based linear peridynamic solid (LPS) model in a body occupying the domain $\Omega\subset\mathbb{R}^d$, { $d = 2$ or $3$}. Let $\theta$ be the nonlocal dilitation, generalizing the local divergence of displacement, and $K(r)$ denote a positive radial function with compactly supported on the $\delta$-ball $B_\delta(\mathbf{x})$. The momentum balance and nonlocal dilitation are then given by the following, 
\begin{equation}\label{eq:nonlocElasticity}
    \mcL_\delta \mathbf{u}:=-\frac{C_\alpha}{m(\delta)}  \int_{B_\delta (\mathbf{x})} \left(\lambda- \mu\right) K(\left|\mathbf{y}-\mathbf{x}\right|) \left(\mathbf{y}-\mathbf{x} \right)\left(\theta(\mathbf{x}) + \theta(\mathbf{y}) \right) d\mathbf{y}
    \end{equation}
    \begin{equation*}
  -  \frac{C_\beta}{m(\delta)}\int_{B_\delta (\mathbf{x})} \mu K(\left|\mathbf{y}-\mathbf{x}\right|)\frac{\left(\mathbf{y}-\mathbf{x}\right)\otimes\left(\mathbf{y}-\mathbf{x}\right)}{\left|\mathbf{y}-\mathbf{x}\right|^2}  \left(\mathbf{u}(\mathbf{y}) - \mathbf{u}(\mathbf{x}) \right) d\mathbf{y} = \mathbf{f}(\xb),
\end{equation*}
\begin{equation}\label{eqn:oritheta}
\theta(\mathbf{x}):=\dfrac{d}{m(\delta)}\int_{B_\delta (\mathbf{x})} K(\left|\mathbf{y}-\mathbf{x}\right|) (\mathbf{y}-\mathbf{x})\cdot \left(\mathbf{u}(\mathbf{y}) - \mathbf{u}(\mathbf{x}) \right)d\mathbf{y},
\end{equation}
where $\mathbf{u}\in\mathbb{R}^d$ denotes the displacement, $\fb\in\mathbb{R}^d$ denotes the body load, the weighted volume 
$$m(\delta):=\int_{B_\delta (\mathbf{x})}K(\left|\mathbf{y}-\mathbf{x}\right|)\left|\mathbf{y}-\mathbf{x}\right|^2d\mathbf{y},$$
and $\mu$, $\lambda$ denote the shear and Lame modulus, respectively. With appropriate choice of scaling parameters $C_\alpha>0$, $C_\beta>0$ and the weighting function $K(r)$, it can be shown that the system converges to the Navier equations \cite{mengesha2012nonlocal,mengesha2014bond,mengesha2014nonlocal}:
\begin{equation}\label{eqn:local}
\mcL_0 \mathbf{u}:=-\nabla\cdot(\lambda tr(\mathbf{E})\mathbf{I}+2\mu \mathbf{E})=-(\lambda-\mu)\nabla [\text{tr}(\mathbf{E})]-\mu \nabla\cdot(2\mathbf{E}+\text{tr}(\mathbf{E})\mathbf{I})=\mathbf{f},
\end{equation}
where the strain tensor $\mathbf{E}:=\dfrac{1}{2}(\nabla \mathbf{u}+(\nabla \mathbf{u})^T)$ and we note that $\text{tr}(\mathbf{E})=\nabla\cdot\mathbf{u}$. To recover parameters for 3D linear elasticity, $C_\alpha=3$, $C_\beta=30$. For 2D problems, $C_\alpha=2$, $C_\beta=16$. In this paper we consider 2D problems ($d=2$) and the following popular scaled kernel:
\begin{equation}\label{eqn:K}
{K(r)}=\left\{\begin{array}{cl}
\dfrac{1}{r},\; &\text{ for }r\leq \delta; \\
0,\;& \text{ for }r> \delta,\\
\end{array}\right.\qquad\dfrac{K(r)}{m(\delta)}=\left\{\begin{array}{cl}
\dfrac{3}{2\pi\delta^3r},\; &\text{ for }r\leq \delta; \\
0,\;& \text{ for }r> \delta. \\
\end{array}
\right.
\end{equation}  
although the idea may be generalized to more general kernels and 3D cases. 
As shown in \cite{mengesha2014nonlocal}, we can define the nonlocal strain energy density as
\begin{align*}
W_\delta(\mathbf{u})=&\frac{C_\alpha d\left(\lambda - \mu\right)}{(m(\delta))^2}\int_\Omega \left[\int_{B_\delta (\mathbf{x})} K(\left|\mathbf{y}-\mathbf{x}\right|) (\mathbf{y}-\mathbf{x})\cdot \left(\mathbf{u}(\mathbf{y}) - \mathbf{u}(\mathbf{x}) \right)d\mathbf{y}\right]^2 d\mathbf{x}\\
&+  \frac{C_\beta\mu}{2m(\delta)}\int_\Omega\left[\int_{B_\delta (\mathbf{x})}  \dfrac{K(\left|\mathbf{y}-\mathbf{x}\right|)}{\left|\mathbf{y}-\mathbf{x}\right|^2}\left[\left(\mathbf{u}(\mathbf{y}) - \mathbf{u}(\mathbf{x}) \right)\cdot\left(\mathbf{y}-\mathbf{x}\right)\right]^2 d\mathbf{y}\right]d\mathbf{x},
\end{align*}
and the energy space $S_\delta(\Omega)$ as
$$S_\delta(\Omega):=\left\{\ub\in L^2(\Omega):|\ub|_{S_\delta(\Omega)}:=\frac{1}{m(\delta)}\int_\Omega\int_\Omega\dfrac{K(\left|\mathbf{y}-\mathbf{x}\right|)}{\left|\mathbf{y}-\mathbf{x}\right|^2}\left[\left(\mathbf{u}(\mathbf{y}) - \mathbf{u}(\mathbf{x}) \right)\cdot\left(\mathbf{y}-\mathbf{x}\right)\right]^2 d\mathbf{y}d\mathbf{x}<\infty \right\}.$$
Note that $|\ub|_{S_\delta(\Omega)}=0$ if and only if $\ub$ represents an infinitesimally rigid displacement, i.e.:
$$\ub(\xb)\in\Pi:=\{\mathbb{Q}\xb+\bb,\mathbb{Q}\in\real^{d\times d},\mathbb{Q}^T=-\mathbb{Q},\bb\in\real^d\}.$$

\section{Neumann and Mixed-type Constraint Problems}\label{sec:n}

%For notation simplicity, in the following we denote regions close to the domain boundary $\partial\Omega$ as
%\begin{align*}
%\omgi&:=\{\mathbf{x}\in\Omega|\text{dist}(\mathbf{x},\partial\Omega)<\delta\},\,\omgb:=\{\mathbf{x}\notin\Omega|\text{dist}(\mathbf{x},\partial\Omega)<\delta\},\,\omgbb:=\{\mathbf{x}\notin\Omega|\text{dist}(\mathbf{x},\partial\Omega)<2\delta\}.
%\end{align*}

\begin{figure}[!htb]\centering
 %\subfigure{\includegraphics[width=0.53\textwidth]{domain1.eps}}\quad
 \subfigure{\includegraphics[width=0.25\textwidth]{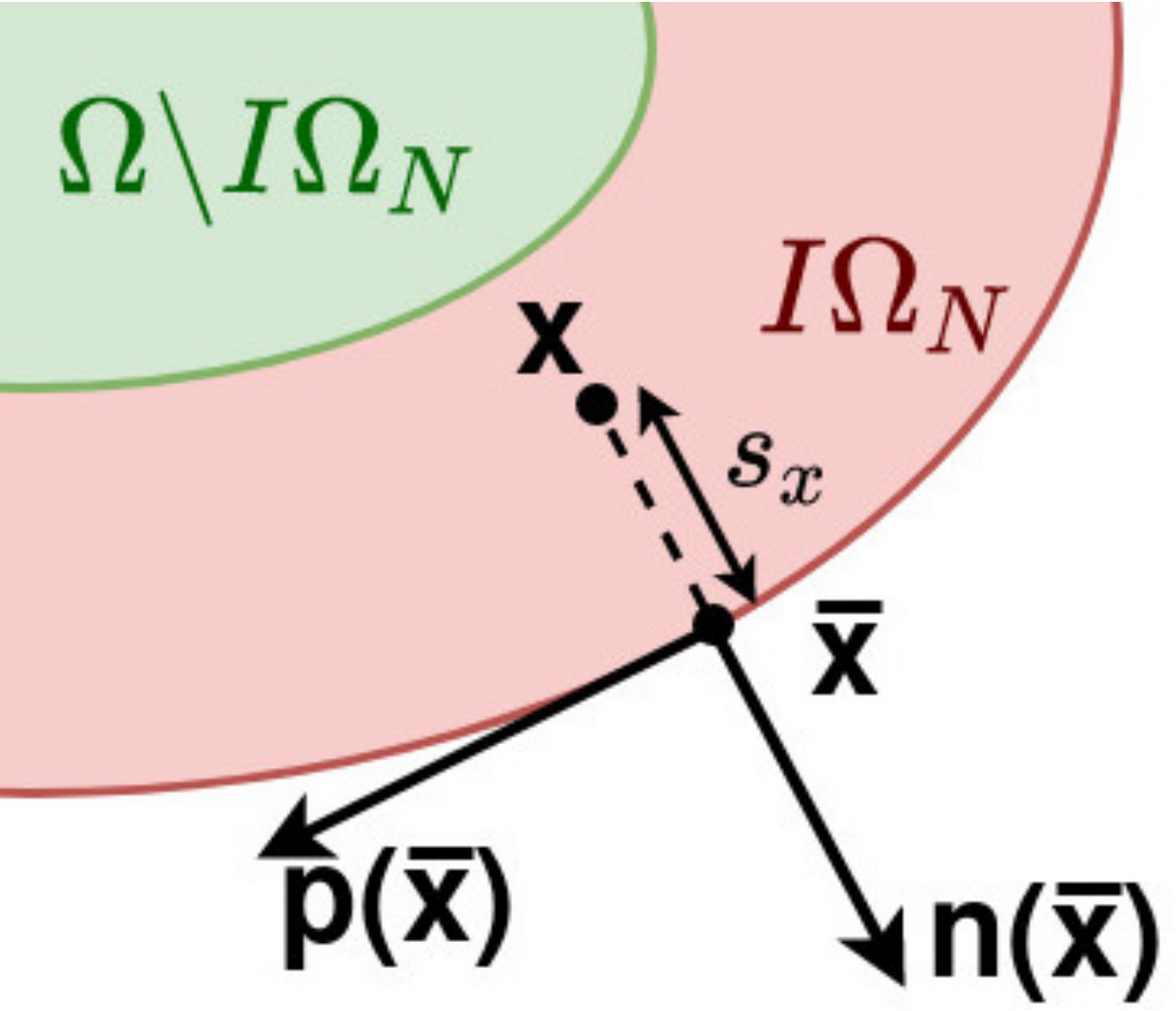}}
\caption{Notations for the projection of point $\mathbf{x}\in \omgi$, the corresponding unit tangential vector $\mathbf{p}(\overline{\mathbf{x}})$ and the unit normal vector $\mathbf{n}(\overline{\mathbf{x}})$.
 }
 \label{fig1}
\end{figure}

We now consider a state-based peridynamic problem with general mixed boundary conditions: $\partial\Omega=\partial\Omega_D\bigcup \partial\Omega_N$ %
and $(\partial\Omega_D)^o\bigcap (\partial\Omega_N)^o=\emptyset$. %
Here $\partial \Omega_D$ and $\partial \Omega_N$ are both 1D curves. We denote the regions near the boundary $\partial\Omega$ as
\begin{align*}
\omgi&:=\{\mathbf{x}\in\Omega|\text{dist}(\mathbf{x},\partial\Omega)<\delta\},\,\omgb:=\{\mathbf{x}\notin\Omega|\text{dist}(\mathbf{x},\partial\Omega)<\delta\},\,\omgbb:=\{\mathbf{x}\notin\Omega|\text{dist}(\mathbf{x},\partial\Omega)<2\delta\}.
\end{align*}
Note that to apply the nonlocal Dirichlet-type boundary condition, $\mathbf{u}(\mathbf{x})=\mathbf{u}_D(\mathbf{x})$ is required in a layer with non-zero volume outside $\Omega$, while the proposed traction load is applied as a Neumann boundary condition on the sharp interface $\partial\Omega_N$. To define a Dirichlet-type constraint, we denote
\begin{align*}
\omgi_D&:=\{\mathbf{x}\in\Omega|\text{dist}(\mathbf{x},\partial\Omega_D)<\delta\},\,\omgb_D:=\{\mathbf{x}\notin\Omega|\text{dist}(\mathbf{x},\partial\Omega_D)<\delta\},\,\omgbb_D:=\{\mathbf{x}\notin\Omega|\text{dist}(\mathbf{x},\partial\Omega_D)<2\delta\},
\end{align*}
and assume that the value of $\ub$ is given on $\omgbb_D$. Similarly, to apply the Neumann constraint, we denote
\begin{align*}
\omgi_N&:=\{\mathbf{x}\in\Omega|\text{dist}(\mathbf{x},\partial\Omega_N)<\delta\},\,\omgb_N:=\{\mathbf{x}\notin\Omega|\text{dist}(\mathbf{x},\partial\Omega_N)<\delta\},\,\omgbb_N:=\{\mathbf{x}\notin\Omega|\text{dist}(\mathbf{x},\partial\Omega_N)<2\delta\}.
\end{align*}
Unless stated otherwise, in this paper we further assume sufficient regularity in the boundary that we may take $\delta$ sufficiently small so that for any $\mathbf{x}\in \omgi_N$ (see Figure \ref{fig1} for illustration), there exists a unique orthogonal projection of $\mathbf{x}$ onto $\partial\Omega_N$. We denote this projection as $\overline{\mathbf{x}}$. Therefore, one has $\overline{\mathbf{x}}-\mathbf{x}=s_x\mathbf{n}(\overline{\mathbf{x}})$ for $\mathbf{x}\in \omgi_N$, where $0<s_x<\delta$. Here $\mathbf{n}$ denotes the normal direction pointing out of the domain for each $\mathbf{x}\in\omgi_N$, and $\mathbf{p}$ denotes the tangential direction. Moreover, we employ the following notations for the directional components of the Hessian matrix of a scalar function $v$:
\begin{align*}
[v({{\mathbf{x}}})]_{pp}:=\mathbf{p}^T(\overline{\mathbf{x}}) \nabla^2v({{\mathbf{x}}}) \mathbf{p}(\overline{\mathbf{x}}),\quad
[v({{\mathbf{x}}})]_{nn}:=\mathbf{n}^T(\overline{\mathbf{x}}) \nabla^2v({{\mathbf{x}}}) \mathbf{n}(\overline{\mathbf{x}}),\quad
[v({{\mathbf{x}}})]_{pn}:=\mathbf{p}^T(\overline{\mathbf{x}}) \nabla^2v({{\mathbf{x}}}) \mathbf{n}(\overline{\mathbf{x}}).
\end{align*}
%{\color{red}
%and the higher order derivative components are similarly defined.
%}

\subsection{Formulation for Non-Homogeneous Traction Loading}\label{sec:Neumann}

In this section, we consider an LPS model subject to local traction loads on the sharp interface $\partial\Omega_N$, by developing nonlocal Neumann constraint formulation with proper correction terms for $\xb\in\omgi_N$. 

Firstly, we propose a corrected formulation for the nonlocal dilitation $\theta$ in \eqref{eqn:oritheta}. When $\mathbf{u} \in C^1(\Omega)$ and $B_\delta(\mathbf{x})\backslash \Omega = \emptyset$, the definition of $\theta(\mathbf{x})$ limits to a local divergence operator $\nabla\cdot\ub(\xb)$ as $\delta \rightarrow 0$ by taking the Taylor series expansion of $\ub$ as $\mathbf{u}(\mathbf{y}) = \mathbf{u}(\mathbf{x}) + \nabla \mathbf{u}(\mathbf{x}) \cdot \left( \mathbf{y} - \mathbf{x} \right) + O(\delta^2)$ and employing a symmetry argument. However, for $\xb\in\omgi_N$ the domain of integration is non-spherical due to proximity to the boundary, and the loss of symmetry results in an inconsistent $\theta$. Thus, surface-effects manifest in the definition of dilitation before any modeling assumptions are made regarding the material response. To address the surface-effect we modify the definition of nonlocal dilitation in \eqref{eqn:oritheta} to enforce consistency for linear displacement fields, independent of whether the horizon intersects the boundary of the domain. In the spirit of correspondence models and corrected { smoothed particle hydrodynamics (SPH)} schemes \cite{oger2007improved}, we introduce a correction tensor $\mathbf{M}(\mathbf{x})$ to \eqref{eqn:oritheta}:
\begin{equation}\label{eq:continuousNonlocDilitation3}
  \theta^{corr}(\mathbf{x}) =  \frac{d}{m(\delta)} \int_{B_\delta (\mathbf{x})\cap\Omega} K(\left|\mathbf{y}-\mathbf{x}\right|) \left(\mathbf{y}-\mathbf{x}\right) \cdot \mathbf{M}(\mathbf{x})\cdot \left(\mathbf{u}(\mathbf{y}) - \mathbf{u}(\mathbf{x}) \right) d\mathbf{y},
\end{equation}
\begin{equation}\label{eq:continuousdilCorr}
  \mathbf{M}(\mathbf{x}) =  \left[ \frac{d}{m(\delta)} \int_{B_\delta (\mathbf{x})\cap \Omega} K(\left|\mathbf{y}-\mathbf{x}\right|) \left(\mathbf{y}-\mathbf{x}\right) \otimes \left(\mathbf{y}-\mathbf{x}\right)  d\mathbf{y} \right]^{-1}.
\end{equation}
Note that for $\xb\in\Omega\backslash\omgi_N$, $\mathbf{M}(\mathbf{x})$ reverts to the identity matrix and \eqref{eq:continuousNonlocDilitation3} reverts to \eqref{eqn:oritheta}. With a slight abuse of notation, we denote $\theta^{corr}(\xb)$ as $\theta$ in the remainder. In the next section, we will further show that for sufficiently smooth domain $\Omega$ and $\ub\in C^1(\Omega)$, the modified diliation is well-posed and consistent with the local dilitation.

We next introduce a Neumann constraint to impose a traction load $\mathbf{T}$ on $\partial \Omega_N$ by modifying the state-peridynamic peridynamic model \eqref{eq:nonlocElasticity} in $\omgi_N$. Denoting $T_p$ and $T_n$ as the tangential and normal components of $\mathbf{T}$, respectively, we propose the following formulation:
 \begin{align}
    \nonumber\mcL_{N\delta}\ub(\xb):=&-\frac{C_\alpha}{m(\delta)}  \int_{B_\delta (\mathbf{x}) \cap \Omega} \left(\lambda - \mu\right) K(\left|\mathbf{y}-\mathbf{x}\right|) 
     \left(\mathbf{y}-\mathbf{x} \right)\left(\theta(\mathbf{x}) + \theta(\mathbf{y}) \right) d\mathbf{y}\\
   \nonumber&-\frac{C_\beta}{m(\delta)}\int_{B_\delta (\mathbf{x}) \cap \Omega} \mu
     K(\left|\mathbf{y}-\mathbf{x}\right|)\frac{\left(\mathbf{y}-\mathbf{x}\right)\otimes\left(\mathbf{y}-\mathbf{x}\right)}{\left|\mathbf{y}-\mathbf{x}\right|^2} 
      \left(\mathbf{u}(\mathbf{y}) - \mathbf{u}(\mathbf{x}) \right) d\mathbf{y}\\
    \nonumber&-\frac{2C_\alpha\theta(\mathbf{x})}{m(\delta)}  \int_{B_\delta (\mathbf{x}) \backslash \Omega} 
     \left(\lambda - \mu\right) K(\left|\mathbf{y}-\mathbf{x}\right|) 
     \left(\mathbf{y}-\mathbf{x} \right) d\mathbf{y}\\
\nonumber&-\frac{C_\beta\theta(\mathbf{x})}{2m(\delta)}\int_{B_\delta (\mathbf{x}) \backslash\Omega}(\lambda+2\mu) K(\left|\mathbf{y}-\mathbf{x}\right|)
     \frac{[\left(\mathbf{y}-\mathbf{x} \right)\cdot \mathbf{n}][\left(\mathbf{y}-\mathbf{x} \right)\cdot \mathbf{p}]^2}{\left|\mathbf{y}-\mathbf{x}\right|^2}
      \mathbf{n} d\mathbf{y} \\
 \nonumber&+\frac{C_\beta\theta(\mathbf{x})}{2m(\delta)} \int_{B_\delta (\mathbf{x}) \backslash \Omega} \lambda K(\left|\mathbf{y}-\mathbf{x}\right|)   \frac{[\left(\mathbf{y}-\mathbf{x} \right)\cdot \mathbf{n}]^3}{\left|\mathbf{y}-\mathbf{x}\right|^2}\mathbf{n} d\mathbf{y}\\
\nonumber=& \fb(\xb)+\frac{C_\beta}{m(\delta)} \int_{B_\delta (\mathbf{x}) \backslash \Omega} K(\left|\mathbf{y}-\mathbf{x}\right|) 
     \frac{[\left(\mathbf{y}-\mathbf{x} \right)\cdot \mathbf{n}]}{\left|\mathbf{y}-\mathbf{x}\right|^2} 
      [\left(\mathbf{y}-\mathbf{x} \right)\cdot \mathbf{p}]^2[T_p(\bar{\mathbf{x}})\mathbf{p}] d\mathbf{y}\\
\nonumber &+\frac{C_\beta}{2m(\delta)} \int_{B_\delta (\mathbf{x}) \backslash \Omega} K(\left|\mathbf{y}-\mathbf{x}\right|) 
     \frac{[\left(\mathbf{y}-\mathbf{x} \right)\cdot \mathbf{n}]}{\left|\mathbf{y}-\mathbf{x}\right|^2} 
      \left([\left(\mathbf{y}-\mathbf{x} \right)\cdot \mathbf{n}]^2-[\left(\mathbf{y}-\mathbf{x} \right)\cdot \mathbf{p}]^2\right)[T_n(\bar{\mathbf{x}})\mathbf{n}] d\mathbf{y}\\
      :=&\fb_{N\delta}(\xb),\label{eq:newform1}
 \end{align}
 where $\bar{\mathbf{x}}$ is the projection of $\mathbf{x}$ on the boundary. In the next section, we will show that this formulation provides an approximation for the corresponding linear elastic model with local traction loadings in the case of linear displacement fields.

To summarize, we obtain a formulation for a static state-based peridynamic problem with general mixed boundary conditions:
 \begin{equation}\label{eqn:probn}
\left\{\begin{array}{ll}
    \mcL_\delta\ub(\xb) = \mathbf{f}(\xb),&\quad \text{ in }\Omega\backslash\omgi_N\\
\mcL_{N\delta}\ub(\xb) = \mathbf{f}_{N\delta}(\xb),&\quad \text{ in }\omgi_N\\
%%%%%%%%%%%%%%%%%%%%%%%%%%%%%%%     
\theta(\mathbf{x})=\dfrac{d}{m(\delta)}\int_{B_\delta (\mathbf{x})} K(\left|\mathbf{y}-\mathbf{x}\right|) (\mathbf{y}-\mathbf{x})^T \left(\mathbf{u}(\mathbf{y}) - \mathbf{u}(\mathbf{x}) \right)d\mathbf{y},&\quad \text{ in }\omg\cup\omgb_D\backslash\omgi_N\\
\theta(\mathbf{x})=\dfrac{d}{m(\delta)}\int_{B_\delta (\mathbf{x})\cap \Omega} K(\left|\mathbf{y}-\mathbf{x}\right|) (\mathbf{y}-\mathbf{x})^T \mathbf{M}(\mathbf{x})\left(\mathbf{u}(\mathbf{y}) - \mathbf{u}(\mathbf{x}) \right)d\mathbf{y},&\quad \text{ in }\omgi_N\\
\mathbf{u}(\mathbf{x})=\mathbf{u}_D(\mathbf{x}), &\quad \text{ in }\omgbb_D
\end{array}\right.
\end{equation}
 where the correction tensor is defined as
 $$\mathbf{M}:=\left[\dfrac{d}{m(\delta)}\int_{B_\delta (\mathbf{x})\cap \Omega}K(\left|\mathbf{y}-\mathbf{x}\right|)(\mathbf{y}-\mathbf{x})\otimes(\mathbf{y}-\mathbf{x})d\mathbf{y}\right]^{-1},$$
 and a body load $\fb_{N\delta}$ is defined on $\xb\in\omgi_N$ as
 \begin{align*}
\mathbf{f}_{N\delta}(\xb):=&\fb(\xb)+\frac{C_\beta}{m(\delta)} \int_{B_\delta (\mathbf{x})\backslash\Omega} K(\left|\mathbf{y}-\mathbf{x}\right|) 
     \frac{[\left(\mathbf{y}-\mathbf{x} \right)\cdot \mathbf{n}]}{\left|\mathbf{y}-\mathbf{x}\right|^2} 
      [\left(\mathbf{y}-\mathbf{x} \right)\cdot \mathbf{p}]^2[T_p(\bar{\mathbf{x}})\mathbf{p}] d\mathbf{y}\\
 &+\frac{C_\beta}{2m(\delta)} \int_{B_\delta (\mathbf{x})\backslash\Omega} K(\left|\mathbf{y}-\mathbf{x}\right|) 
     \frac{[\left(\mathbf{y}-\mathbf{x} \right)\cdot \mathbf{n}]}{\left|\mathbf{y}-\mathbf{x}\right|^2} 
      \left([\left(\mathbf{y}-\mathbf{x} \right)\cdot \mathbf{n}]^2-[\left(\mathbf{y}-\mathbf{x} \right)\cdot \mathbf{p}]^2\right)[T_n(\bar{\mathbf{x}})\mathbf{n}] d\mathbf{y}.    
\end{align*}

\subsection{Well-posedness and Consistency Analysis}

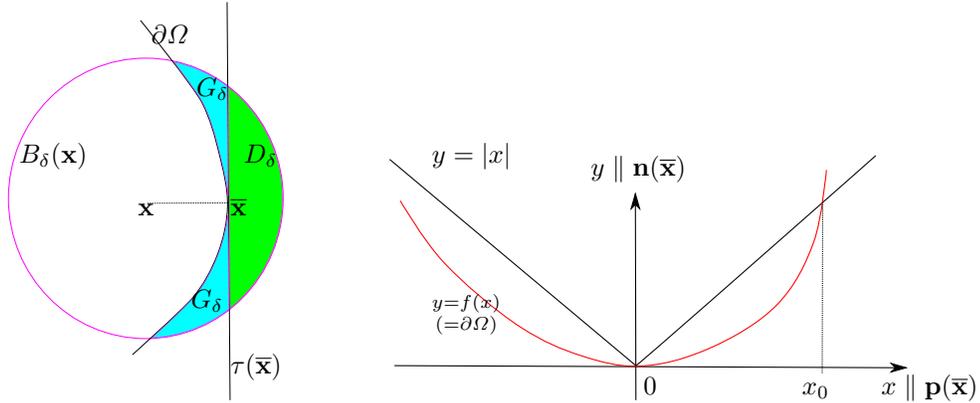
\begin{figure}[!htb]
\centering
\subfigure{\includegraphics[scale=.4]{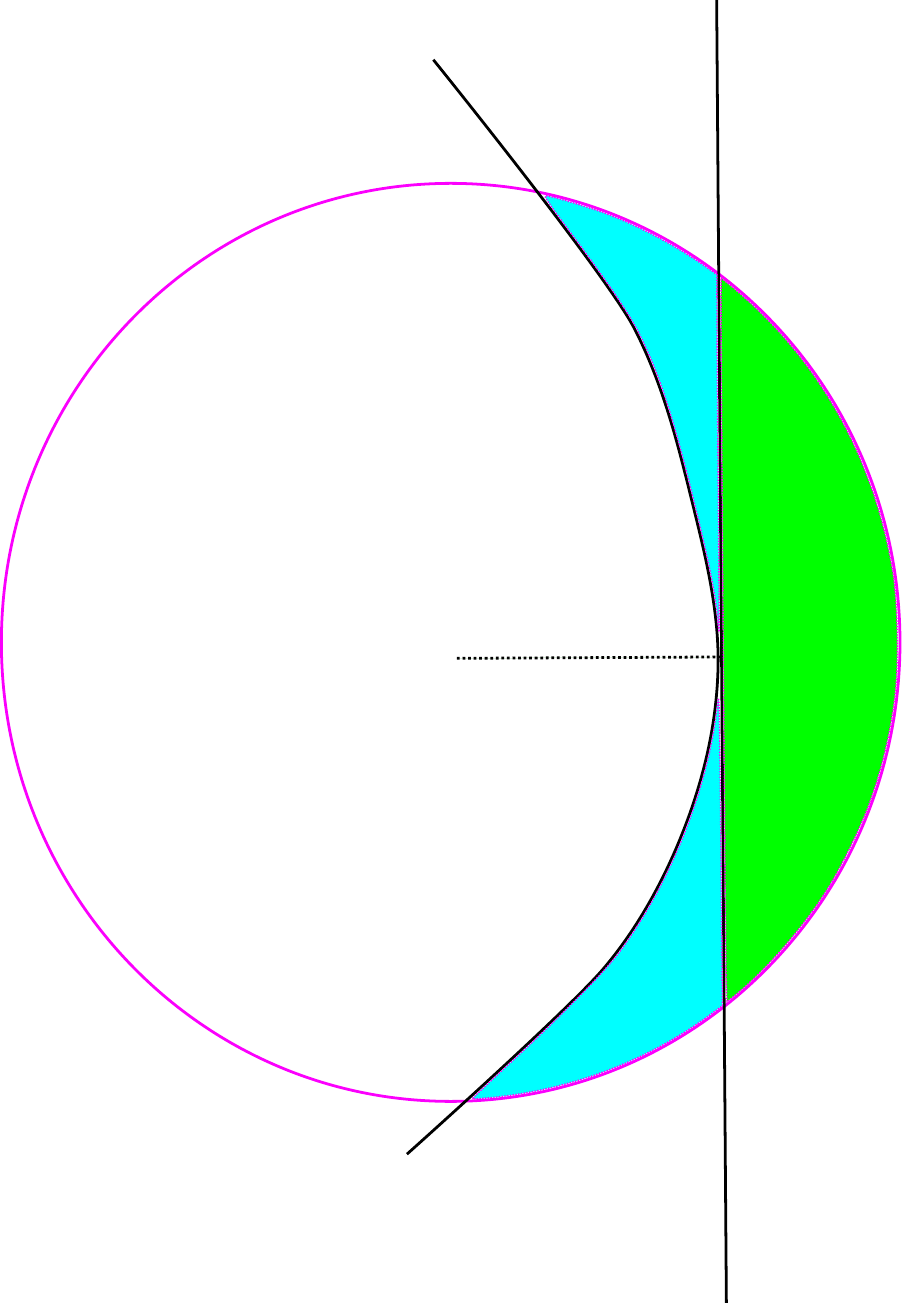}
 \put(-55,70){$\mathbf{x}$}\put(-20,70){$\overline{\mathbf{x}}$}
 \put(-20,10){$\tau(\overline{\mathbf{x}})$}
 \put(-15,90){$D_\delta$}
 \put(-33,115){$G_\delta$}\put(-35,35){$G_\delta$}
 \put(-50,135){$\partial\Omega$}
 \put(-100,90){$B_\delta(\mathbf{x})$}}\qquad\qquad
\subfigure{\includegraphics[scale=0.5]{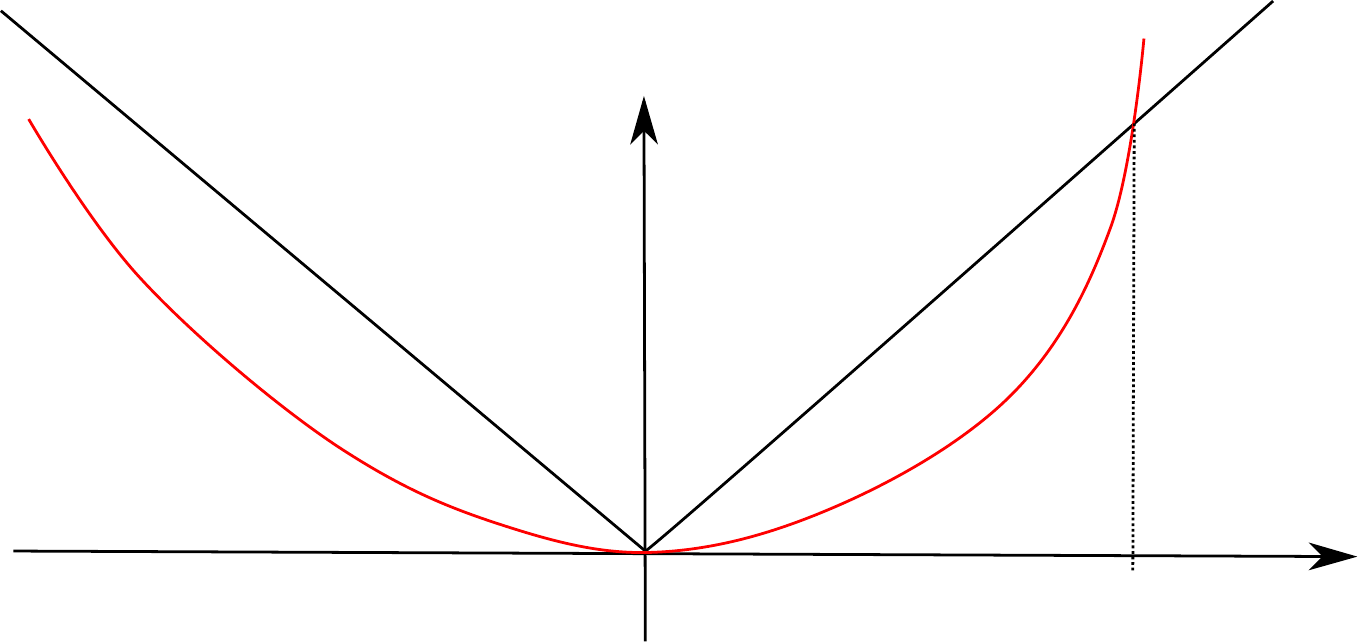}
 \put(-100,2){$0$}
 \put(-120,85){$y\parallel\mathbf{n}(\overline{\mathbf{x}}) $}
 \put(-10,2){$x\parallel\mathbf{p}(\overline{\mathbf{x}}) $}
 \put(-40,2){$x_0$}
 \put(-180,30){$\substack{y=f(x)\\ (=\partial\Omega)}$}
 \put(-180,90){$y=|x|$}} 
 \caption{Notation for geometric estimates: Left: illustration of regions $D_\delta$ and $G_\delta$. Green represents $D_\delta$, the region in $B_\delta(\mathbf{x})$ which lies opposite the boundary tangent at $\overline{\mathbf{x}}$. Cyan represents $G_\delta$, the region in $B_\delta(\mathbf{x})$ which lies between $\partial\Omega$ and the boundary tangent. Right: Local Cartesian coordinate system in neighborhood of $\overline{\mathbf{x}}$. Here, the region $G_\delta$ lies below the red curve $y=f(x)$ when $\mathbf{x}=\overline{\mathbf{x}}$.}
 \label{graph}
\end{figure}

In this section, we will show that the modified diliation is well-posed and consistent with the local dilitation. Specifically, we prove that for sufficiently smooth domain, the correction tensor $\mathbf{M}$ is invertible, and that for $\mathbf{u}\in C^1(\Omega)$, $\theta\rightarrow\nabla\cdot\mathbf{u}$ as $\delta\rightarrow 0$. Moreover, we will demonstrate that for linear displacement $\ub$ and under certain geometric assumptions, the modified formulation \eqref{eqn:probn} is consistent with the classical linear elastic problem with traction loadings. For simplicity of notation, we indicate a generic constant independent of $\delta$ as $C$, and write $K(|\yb-\xb|)$ as $K$.

We first analyze existence and bounds of $\mathbf{M}$: 
\begin{thm}\label{lem:M}
Given that $\Omega\in \mathbb{R}^d$ $(d=2)$ is a {$C^{3}$} domain, then there exists a $\overline{\delta}>0$  such that for $0<\delta\leq \overline{\delta}$ the correction tensor is a well-defined symmetric matrix, and
\begin{align*}
    \mathbf{M}&=\left[\begin{array}{cc}
    \dfrac{d}{m(\delta)}\int_{B_\delta (\mathbf{x})\cap \Omega}K((\mathbf{y}-\mathbf{x})\cdot{\mathbf{p}})^2 d\mathbf{y}&O(\delta^2)\\
    O(\delta^2)&\dfrac{d}{m(\delta)}\int_{B_\delta (\mathbf{x})\cap \Omega}K((\mathbf{y}-\mathbf{x})\cdot{\mathbf{n}})^2d\mathbf{y}\\
    \end{array}\right]^{-1}\\
    &=\left[\begin{array}{cc}
    \left(\dfrac{d}{m(\delta)}\int_{B_\delta (\mathbf{x})\cap \Omega}K((\mathbf{y}-\mathbf{x})\cdot{\mathbf{p}})^2 d\mathbf{y}\right)^{-1}+O(\delta^4)&O(\delta^2)\\
    O(\delta^2)&\left(\dfrac{d}{m(\delta)}\int_{B_\delta (\mathbf{x})\cap \Omega}((\mathbf{y}-\mathbf{x})\cdot{\mathbf{n}})^2d\mathbf{y})\right)^{-1}+O(\delta^4)\\
    \end{array}\right].\\
\end{align*}
\end{thm}
\begin{proof}
% $\kappa({\mathbf{x}})$ denotes the curvature of $\partial\Omega$ at ${\mathbf{x}}$. 
To show that the correction tensor $\mathbf{M}$ is well-defined, it suffices to show that $|det(\mathbf{M}^{-1})|>0$. We adopt notation in Figure \ref{graph}, with a Cartesian coordinate system oriented so that $\overline{\mathbf{x}}$ coincides with the origin, and the vectors $\mathbf{p}(\overline{\mathbf{x}})$ and $\mathbf{n}(\overline{\mathbf{x}})$ are oriented along the positive $x$-axis and negative $y$-axis, respectively. We note that
\begin{align*}
    \mathbf{M}^{-1}&=\dfrac{d}{m(\delta)}\left[\begin{array}{cc}
    \int_{B_\delta (\mathbf{x})\cap \Omega}K((\mathbf{y}-\mathbf{x})\cdot{\mathbf{p}})^2 d\mathbf{y}&\int_{B_\delta (\mathbf{x})\cap \Omega}K((\mathbf{y}-\mathbf{x})\cdot{\mathbf{p}})((\mathbf{y}-\mathbf{x})\cdot{\mathbf{n}})d\mathbf{y}\\
    \int_{B_\delta (\mathbf{x})\cap \Omega}K((\mathbf{y}-\mathbf{x})\cdot{\mathbf{p}})((\mathbf{y}-\mathbf{x})\cdot{\mathbf{n}})d\mathbf{y}&\int_{B_\delta (\mathbf{x})\cap \Omega}K((\mathbf{y}-\mathbf{x})\cdot{\mathbf{n}})^2d\mathbf{y}\\
    \end{array}\right]\\
    &=\left[\begin{array}{cc}
    1&0\\
    0&1\\
    \end{array}\right]
    -\dfrac{d}{m(\delta)}\left[\begin{array}{cc}
    \int_{D_\delta}K((\mathbf{y}-\mathbf{x})\cdot{\mathbf{p}})^2 d\mathbf{y}&\int_{D_\delta}K((\mathbf{y}-\mathbf{x})\cdot{\mathbf{p}})((\mathbf{y}-\mathbf{x})\cdot{\mathbf{n}})d\mathbf{y}\\
    \int_{D_\delta}K((\mathbf{y}-\mathbf{x})\cdot{\mathbf{p}})((\mathbf{y}-\mathbf{x})\cdot{\mathbf{n}})d\mathbf{y}&\int_{D_\delta}K((\mathbf{y}-\mathbf{x})\cdot{\mathbf{n}})^2d\mathbf{y}\\
    \end{array}\right]\\
    &-\dfrac{d}{m(\delta)}\left[\begin{array}{cc}
    \int_{G_\delta}K((\mathbf{y}-\mathbf{x})\cdot{\mathbf{p}})^2 d\mathbf{y}&\int_{G_\delta}K((\mathbf{y}-\mathbf{x})\cdot{\mathbf{p}})((\mathbf{y}-\mathbf{x})\cdot{\mathbf{n}})d\mathbf{y}\\
    \int_{G_\delta}K((\mathbf{y}-\mathbf{x})\cdot{\mathbf{p}})((\mathbf{y}-\mathbf{x})\cdot{\mathbf{n}})d\mathbf{y}&\int_{G_\delta}K((\mathbf{y}-\mathbf{x})\cdot{\mathbf{n}})^2d\mathbf{y}\\
    \end{array}\right].
\end{align*}
%With $\tau(\overline{\mathbf{x}})$ representing the tangent line to $\partial\Omega$ at  $\overline{\mathbf{x}}$, here $D_\delta$ is the region of $B_\delta(\mathbf{x})$ on the side of $\tau(\overline{\mathbf{x}})$ {\em not} containing $\Omega$ (as shown in the green region in the left plot of Figure \ref{graph}), and $G_\delta:=B_\delta(\mathbf{x}) \backslash (D_\delta\cup \Omega)$ (as shown in the cyan region in the left plot of Figure \ref{graph}). %
We estimate first the $D_\delta$ part. Rewriting %
$\mathbf{y}\in D_\delta%B(\mathbf{x},\delta)\backslash\Omega
$ as $\mathbf{x}+(r\cos(\theta),r\sin(\theta))$ with $s_x<r<\delta$ and %
$-\pi/2\leq-\arccos(s_x/r)\leq\theta\leq \arccos(s_x/r)\leq \pi/2$, we obtain
$$\int_{D_\delta}K((\mathbf{y}-\mathbf{x})\cdot{\mathbf{p}})^2 d\mathbf{y}=\int_{s_x}^\delta\int_{-\arccos(s_x/r)}^{\arccos(s_x/r)} K(r)r^3\sin^2\theta d\theta dr\leq\int_{0}^\delta\int_{-\pi/2}^{\pi/2} K(r)r^3\sin^2\theta d\theta dr=m(\delta)/4,$$
$$\int_{D_\delta}K((\mathbf{y}-\mathbf{x})\cdot{\mathbf{n}})^2 d\mathbf{y}=\int_{s_x}^\delta\int_{-\arccos(s_x/r)}^{\arccos(s_x/r)} K(r)r^3\cos^2\theta d\theta dr\leq\int_{0}^\delta\int_{-\pi/2}^{\pi/2} K(r)r^3\cos^2\theta d\theta dr=m(\delta)/4,$$
$$\int_{D_\delta}K((\mathbf{y}-\mathbf{x})\cdot{\mathbf{p}})((\mathbf{y}-\mathbf{x})\cdot{\mathbf{n}})d\mathbf{y}=0,$$
since the first two terms decrease monotonically with increasing $s_x$. We then have
\begin{align*}
    \dfrac{d}{m(\delta)}\left[\begin{array}{cc}
    \int_{D_\delta}K((\mathbf{y}-\mathbf{x})\cdot{\mathbf{p}})^2 d\mathbf{y}&\int_{D_\delta}K((\mathbf{y}-\mathbf{x})\cdot{\mathbf{p}})((\mathbf{y}-\mathbf{x})\cdot{\mathbf{n}})d\mathbf{y}\\
    \int_{D_\delta}K((\mathbf{y}-\mathbf{x})\cdot{\mathbf{p}})((\mathbf{y}-\mathbf{x})\cdot{\mathbf{n}})d\mathbf{y}&\int_{D_\delta}K((\mathbf{y}-\mathbf{x})\cdot{\mathbf{n}})^2d\mathbf{y}\\
    \end{array}\right]=\left[\begin{array}{cc}
    d_p&0\\
    0&d_n\\
    \end{array}\right]
\end{align*}
where $0\leq d_p,d_n\leq 1/2$. We now proceed to show that for a domain with $C^3$ regularity, the magnitude of all elements in the matrix
\begin{align*}
    \dfrac{d}{m(\delta)}\left[\begin{array}{cc}
    \int_{G_\delta}K((\mathbf{y}-\mathbf{x})\cdot{\mathbf{p}})^2 d\mathbf{y}&\int_{G_\delta}K((\mathbf{y}-\mathbf{x})\cdot{\mathbf{p}})((\mathbf{y}-\mathbf{x})\cdot{\mathbf{n}})d\mathbf{y}\\
    \int_{G_\delta}K((\mathbf{y}-\mathbf{x})\cdot{\mathbf{p}})((\mathbf{y}-\mathbf{x})\cdot{\mathbf{n}})d\mathbf{y}&\int_{G_\delta}K((\mathbf{y}-\mathbf{x})\cdot{\mathbf{n}})^2d\mathbf{y}\\
    \end{array}\right]
\end{align*}
are bounded by $O(\delta^2)$. Note that with the Cartesian coordinate system in Figure \ref{graph}, $\overline{\mathbf{x}}=(0,0)$ and $\tau(\overline{\mathbf{x}})=\{y=0\}$. Let: $y=f(x)$ be the curve describing $\partial\Omega$; $\kappa({\mathbf{x}})$ denote the curvature of $\partial\Omega_N$ at $\overline{\mathbf{x}}$; and $(c_1(l),c_2(l))$ be the paramaterization of the boundary $\partial\Omega$ by the arclength $l$. {Note that the range of $l$ depends upon the particular geometry of $\partial\omg$.} Then we have $\mathbf{x}_l=(c_1(l),c_2(l))^T$, and
\begin{displaymath}
\mathbf{x}_l=\overline{\mathbf{x}}+\left(\begin{array}{c}
l\\
0\\
\end{array}\right)+\left(\begin{array}{c}
0\\
\frac{\kappa(\overline{\mathbf{x}})l^2}{2}\\
\end{array}\right)+\left(\begin{array}{c}
c^{'''}_1(0)\frac{l^3}{6}\\
c^{'''}_2(0)\frac{l^3}{6}\\
\end{array}\right)+O(l^4).
\end{displaymath}
The area $|G_\delta|\leq |\kappa(\overline{\mathbf{x}})|\dfrac{\delta^3}{3}+O(\delta^4)$. Therefore, when $\delta$ is sufficiently small, for the kernel $K$ in \eqref{eqn:K} we have
$$\dfrac{d}{m(\delta)}\left|\int_{G_\delta}K((\mathbf{y}-\mathbf{x})\cdot{\mathbf{p}})^2 d\mathbf{y}\right|\leq \dfrac{3}{\pi\delta^3}|{G_\delta}|\delta \leq \dfrac{1}{\pi}|\kappa(\overline{\mathbf{x}})|\delta+O(\delta^2)\leq O(\delta).$$
A similar bound follows for $\dfrac{d}{m(\delta)}\left|\int_{G_\delta}K((\mathbf{y}-\mathbf{x})\cdot{\mathbf{n}})^2 d\mathbf{y}\right|$. 
%\footnote{The above bounds hold for more general kernels $K=C/r^\alpha$ with $\alpha\leq d+2$.}. 
For $\dfrac{d}{m(\delta)}\left|\int_{G_\delta}K((\mathbf{y}-\mathbf{x})\cdot{\mathbf{n}})((\mathbf{y}-\mathbf{x})\cdot{\mathbf{p}}) d\mathbf{y}\right|$, following from the symmetry of $K$,
\begin{displaymath}
\dfrac{d}{m(\delta)}\left|\int_{G_\delta}K((\mathbf{y}-\mathbf{x})\cdot{\mathbf{p}})((\mathbf{y}-\mathbf{x})\cdot{\mathbf{n}}) d\mathbf{y}\right|=\dfrac{d}{m(\delta)}\left|\int_{E_\delta}K((\mathbf{y}-\mathbf{x})\cdot{\mathbf{p}})((\mathbf{y}-\mathbf{x})\cdot{\mathbf{n}}) d\mathbf{y}\right|,
\end{displaymath}
where $E_\delta$ denotes the region in $G_\delta$ which is asymmetric with respect to the $y$ axis in the right plot of Figure \ref{graph}. As shown in \cite{You_2019}, the area of $E_\delta$ has $|E_\delta|\leq O(\delta^4)$. Therefore
$$\dfrac{d}{m(\delta)}\left|\int_{E_\delta}K((\mathbf{y}-\mathbf{x})\cdot{\mathbf{p}})((\mathbf{y}-\mathbf{x})\cdot{\mathbf{n}}) d\mathbf{y}\right|\leq O(\delta^2).$$
For sufficiently small $\delta$ we have 
\begin{align*}
|det(\mathbf{M}^{-1})|\geq& (1-d_p)(1-d_n)-C(2-d_p-d_n)\delta-C\delta^2\geq \dfrac{1}{4}-C\delta>0.
\end{align*}
\end{proof}

\begin{remark}\label{rmk:1}
From the proof of Lemma \ref{lem:M}, we note that when $|E_\delta|=0$, i.e., when $B_\delta(\xb)\backslash\Omega$ is symmetric with respect to $\mathbf{n}(\overline{\xb})$, then
\begin{equation}\label{eqn:M1}
    \mathbf{M}=\left[\begin{array}{cc}
    \left(\dfrac{d}{m(\delta)}\int_{B_\delta (\mathbf{x})\cap \Omega}K((\mathbf{y}-\mathbf{x})\cdot{\mathbf{p}})^2 d\mathbf{y}\right)^{-1}&0\\
    0&\left(\dfrac{d}{m(\delta)}\int_{B_\delta (\mathbf{x})\cap \Omega}((\mathbf{y}-\mathbf{x})\cdot{\mathbf{n}})^2d\mathbf{y})\right)^{-1}\\
    \end{array}\right].
\end{equation}
\end{remark}
% \begin{remark}\label{rmk:2}
% Note that the well-posedness property of $\mathbf{M}$ also holds when $\Omega$ is $C^2$, but it is no longer guaranteed that the off diagonal terms are second-order in $\delta$, and thus we may not be able to prove consistency.

% \YY{[YY: for $C^2$ we don't have the $O(\delta^2)$ estimate for the off-diagonal terms in $M$, only the existence of $M$.]}
% \end{remark}

We now show that the nonlocal dilitation $\theta^{corr}$ is consistent with the local dilitation:

\begin{thm}\label{thm:thetaconsis}
Assume that $\mathbf{u}\in C^1$ and $\Omega$ is a $C^3$ domain, then there exists $\overline{\delta}>0$ such that for any $0<\delta\leq \overline{\delta}$,   
$$|\theta^{corr}(\mathbf{x})-\nabla\cdot\mathbf{u}(\mathbf{x})|=O(\delta)$$
for $\mathbf{x}\in\omgi_N$. If $\ub$ further satisfies $\ub\in C^2$, then
$$|\theta^{corr}(\mathbf{x})-\nabla\cdot\mathbf{u}(\mathbf{x})|=O(\delta^2)+O(\delta)|\ub|_{2,\infty},$$
where $|\ub|_{2,\infty}$ is the Sobolev seminorm representing the maximum of the Hessian matrix elements for each component of $\ub$.
\end{thm}
\begin{proof}
We again adopt the coordinate system from Figure \ref{graph}. Denoting $u_p,u_n$ as the displacement components along the directions of $\mathbf{p}(\overline{\mathbf{x}})$ and $\mathbf{n}(\overline{\mathbf{x}})$, respectively, for $\ub\in C^2$ we have
\begin{align}
\nonumber&\theta^{corr}-\nabla\cdot\mathbf{u}\\
\nonumber=&\dfrac{d}{m(\delta)}\left(M_{11}\int_{B_\delta (\mathbf{x})\cap\Omega}K((\mathbf{y}-\mathbf{x})\cdot{\mathbf{p}})(u_p(\mathbf{y})-u_p(\mathbf{x}))d\mathbf{y}+M_{22}\int_{B_\delta (\mathbf{x})\cap\Omega}K((\mathbf{y}-\mathbf{x})\cdot{\mathbf{n}})(u_n(\mathbf{y})-u_n(\mathbf{x}))d\mathbf{y}\right.\\
\nonumber&\left.+M_{12}\int_{B_\delta (\mathbf{x})\cap\Omega}K((\mathbf{y}-\mathbf{x})\cdot{\mathbf{p}})(u_n(\mathbf{y})-u_n(\mathbf{x}))+K((\mathbf{y}-\mathbf{x})\cdot{\mathbf{n}})(u_p(\mathbf{y})-u_p(\mathbf{x}))d\mathbf{y}\right)-\dfrac{\partial u_p}{\partial \mathbf{p}}(\mathbf{x})-\dfrac{\partial u_n}{\partial \mathbf{n}}(\mathbf{x})\\
\nonumber=&\dfrac{d}{m(\delta)}\left(M_{11}\dfrac{\partial u_p}{\partial \mathbf{p}}(\mathbf{x})\int_{B_\delta (\mathbf{x})\cap\Omega}K|(\mathbf{y}-\mathbf{x})\cdot{\mathbf{p}}|^2d\mathbf{y}+M_{22}\dfrac{\partial u_n}{\partial \mathbf{n}}(\mathbf{x})\int_{B_\delta (\mathbf{x})\cap\Omega}K|(\mathbf{y}-\mathbf{x})\cdot{\mathbf{n}}|^2d\mathbf{y}\right.\\
\nonumber&+M_{11}[u_p(\xb)]_{np}\int_{B_\delta (\mathbf{x})\cap\Omega}K[(\mathbf{y}-\mathbf{x})\cdot{\mathbf{p}}]^2[(\mathbf{y}-\mathbf{x})\cdot{\mathbf{n}}]d\mathbf{y}\\
\nonumber&+M_{22}\int_{B_\delta (\mathbf{x})\cap\Omega}K[(\mathbf{y}-\mathbf{x})\cdot{\mathbf{n}}]([u_n(\xb)]_{nn}[(\mathbf{y}-\mathbf{x})\cdot{\mathbf{n}}]^2+[u_n(\xb)]_{pp}[(\mathbf{y}-\mathbf{x})\cdot{\mathbf{p}}]^2)d\mathbf{y}\\
\nonumber&\left.+M_{12}\left(\dfrac{\partial u_p}{\partial \mathbf{n}}+\dfrac{\partial u_n}{\partial \mathbf{p}}\right)\int_{E_\delta}K((\mathbf{y}-\mathbf{x})\cdot{\mathbf{p}})((\mathbf{y}-\mathbf{x})\cdot{\mathbf{n}})d\mathbf{y}+O(\delta^5)\right)-\dfrac{\partial u_p}{\partial \mathbf{p}}(\mathbf{x})-\dfrac{\partial u_n}{\partial \mathbf{n}}(\mathbf{x})\\
%=&\left(\dfrac{\partial u}{\partial \mathbf{p}}(\mathbf{x})+\dfrac{\partial v}{\partial \mathbf{n}}(\mathbf{x})\right)\dfrac{1}{G}\int_{B_\delta (\mathbf{x})\cap\Omega}K|(\mathbf{y}-\mathbf{x})\cdot{\mathbf{p}}|^2d\mathbf{y}\int_{B_\delta (\mathbf{x})\cap\Omega}K|(\mathbf{y}-\mathbf{x})\cdot{\mathbf{n}}|^2d\mathbf{y}\\
%&+O(\delta)-\dfrac{\partial u}{\partial \mathbf{p}}(\mathbf{x})-\dfrac{\partial v}{\partial \mathbf{n}}(\mathbf{x})\\
%=&O(\delta)\\
&=O(\delta^2)+A_1(\xb)[u_p(\mathbf{x})]_{np}+A_2(\xb)[u_n(\mathbf{x})]_{nn}+A_3(\xb)[u_n(\mathbf{x})]_{pp}=O(\delta),\label{eqn:thetataylor}
\end{align}
where
\begin{align*}
&A_1(\xb):=\dfrac{\int_{B_\delta (\mathbf{x})\cap\Omega} K(\left|\mathbf{y}-\mathbf{x}\right|) [(\mathbf{y}-\mathbf{x})\cdot\mathbf{p}]^2[(\mathbf{y}-\mathbf{x})\cdot\mathbf{n}]d\mathbf{y}}{\int_{B_\delta (\mathbf{x})\cap \Omega}K[(\mathbf{y}-\mathbf{x})\cdot{\mathbf{p}}]^2 d\mathbf{y}},\\
&A_2(\xb):=\dfrac{\int_{B_\delta (\mathbf{x})\cap\Omega} K(\left|\mathbf{y}-\mathbf{x}\right|) [(\mathbf{y}-\mathbf{x})\cdot\mathbf{n}]^3d\mathbf{y}}{2\int_{B_\delta (\mathbf{x})\cap \Omega}K[(\mathbf{y}-\mathbf{x})\cdot{\mathbf{n}}]^2 d\mathbf{y}},\\
&A_3(\xb):=\dfrac{\int_{B_\delta (\mathbf{x})\cap\Omega} K(\left|\mathbf{y}-\mathbf{x}\right|) [(\mathbf{y}-\mathbf{x})\cdot\mathbf{p}]^2[(\mathbf{y}-\mathbf{x})\cdot\mathbf{n}]d\mathbf{y}}{2\int_{B_\delta (\mathbf{x})\cap \Omega}K[(\mathbf{y}-\mathbf{x})\cdot{\mathbf{n}}]^2 d\mathbf{y}}.
\end{align*}
For $\ub\in C^1$, the conclusion can be shown with Taylor expansion following a similar procedure as above.
\end{proof}

%\begin{remark}
%From the proof of Theorem \ref{thm:thetaconsis}, we note that when $|E_\delta|=0$, with the symmetry of $B_\delta(\xb)\cap\Omega$ and Remark \ref{rmk:1} we have

%\end{remark}

%\remark{Note that historically, solution of the LPS model in peridynamics requires the imposition of a volumetric Dirichlet condition over a collar region of width $2\delta$. With the modified definition of dilitation, this is no longer true and one may apply boundary conditions on a width $\delta$ region, identical to e.g. bond-based peridynamics.}

Having proven well-posedness and accuracy of the nonlocal dilitation, we next show that the formulation in \eqref{eqn:probn} approximately passes the linear patch test in the local limit.
\begin{thm}\label{thm:linearpatch}
Given that $\Omega\in \mathbb{R}^d$ $(d=2)$ is a {$C^{3}$} domain, and a linear displacement field $\ub$ which is a solution of the classical linear elastic problem in the absence of forcing term $\fb$:
\begin{displaymath}
 \left\{\begin{array}{ll}
 -\nabla\cdot(\lambda tr(\mathbf{E})\mathbf{I}+2\mu \mathbf{E})=0,\quad \text{where }\mathbf{E}=\dfrac{1}{2}(\nabla \mathbf{u}+(\nabla \mathbf{u})^T),&\quad \text{in }\Omega,\\
 \lambda tr(\mathbf{E})\mathbf{n}+2\mu \mathbf{E}\mathbf{n}= \mathbf{T},&\quad \text{on }\partial\Omega_N,\\
 \mathbf{u}= \mathbf{u}_D,&\quad \text{on }\omgbb_D.
 \end{array}\right.
 \end{displaymath}
 When $|E_\delta|=0$, i.e., $B_\delta(\xb)\backslash\Omega$ is symmetric with respect to $\mathbf{n}(\overline{\xb})$ for all $\xb\in\omgi_N$, $\ub$ is also the solution of the state-based peridynamic problem \eqref{eqn:probn} in the absence of forcing term $\fb$. When $B_\delta(\xb)\backslash\Omega$ is not symmetric, $\ub$ passes the linear patch test approximately in $\omgi_N$, i.e., $\mcL_\delta\ub(\xb) = \mathbf{f}(\xb)$ for $\xb\in \Omega\backslash\omgi_N$, and $\mcL_{N\delta}\ub(\xb) = \mathbf{f}_{N\delta}(\xb)+O(\delta)\mathbf{1}$ for $\xb\in\omgi_N$.
\end{thm}
\begin{proof}
Taking a linear displacement field $\ub=\mathbb{D}\xb+\bb$ where $\mathbb{D}\in\real^{d\times d}$, $\bb\in\real^d$, for $\xb\in\Omega\backslash\omgi_N$ the proof can be found in, e.g., \cite{silling2008convergence}. We therefore focus on $\xb\in\omgi_N$, and again employ the notation in Figure \ref{graph}. Moreover, we denote the elements of $\mathbb{D}$ as $D_{ij}$, $i,j\in\{1,2\}$.

We first consider the case when $|E_\delta|=0$ for all $\xb\in\omgi_N$. Substituting \eqref{eqn:M1} into the definition of $\theta$ in  \eqref{eqn:probn}, we obtain $\theta(\xb)=D_{11}+D_{22}$ for all $\xb\in\Omega\cup\omgb_D$. Note that 
$$T_p=\mu\left(\dfrac{\partial u_p}{\partial \mathbf{n}}+\dfrac{\partial u_n}{\partial \mathbf{p}}\right)=\mu(D_{12}+D_{21}),\;T_n=\lambda\nabla\cdot\ub+2\mu \dfrac{\partial u_n}{\partial \mathbf{n}}=(\lambda+2\mu)D_{22}+\lambda D_{11},$$
the proof of $\mcL_{N\delta}\ub(\xb) = \mathbf{f}_{N\delta}(\xb)$ is obtained via a straightforward substitution of $\ub$ and $\theta$ into \eqref{eq:newform1}.

We now consider the general case. Combining $\ub=\mathbb{D}\xb+\bb$ and \eqref{eqn:thetataylor} yields $\theta(\mathbf{x})=D_{11}+D_{22}+O(\delta^2)$. 
%In the above derivation, we have used the following properties:
%$$\dfrac{d}{m(\delta)}\int_{B_\delta (\mathbf{x})\cap \Omega} K((\yb-\xb)\cdot\mathbf{n})((\yb-\xb)\cdot\mathbf{p})d\mathbf{y}=\dfrac{d}{m(\delta)}\int_{E_\delta} K((\yb-\xb)\cdot\mathbf{n})((\yb-\xb)\cdot\mathbf{p})d\mathbf{y}\leq C\delta^2,$$
%$$\dfrac{d}{m(\delta)}\int_{B_\delta (\mathbf{x})\cap \Omega} KM_{12}(\xb) ((\yb-\xb)\cdot\mathbf{p})^2d\mathbf{y}\leq \dfrac{C}{\delta}\int_{B_\delta (\mathbf{x})\cap \Omega} K((\yb-\xb)\cdot\mathbf{p})^2d\mathbf{y}\leq C\delta^2.$$
Substituting the definitions of $\theta$ and $\ub$ into \eqref{eq:newform1}, 
%and note that
%$$\frac{1}{m(\delta)}\int_{B_\delta(\xb)\cap\Omega}K[(\yb-\xb)\cdot\mathbf{n}]d\xb \leq C|B_\delta(\xb)\cap\Omega|/\delta=O(\delta^{-1}),$$
we then have:
 \begin{align*}
    \nonumber&\mcL_{N\delta}\ub(\xb)-\fb_{N\delta}(\xb)\\
    =&O(\delta)\mathbf{1}-\frac{2C_\alpha(D_{11}+D_{22})}{m(\delta)}  \int_{B_\delta (\mathbf{x})} \left(\lambda - \mu\right) K 
     \left(\mathbf{y}-\mathbf{x} \right) d\mathbf{y}-\frac{C_\beta}{m(\delta)}\int_{B_\delta (\mathbf{x}) \cap \Omega} \mu     K\frac{\left(\mathbf{y}-\mathbf{x}\right)\otimes\left(\mathbf{y}-\mathbf{x}\right)}{\left|\mathbf{y}-\mathbf{x}\right|^2}      \cdot \mathbb{D}(\yb-\xb) d\mathbf{y}\\
\nonumber&-\frac{C_\beta(D_{11}+D_{22})}{2m(\delta)}\int_{B_\delta (\mathbf{x}) \backslash\Omega}(\lambda+2\mu) K
     \frac{[\left(\mathbf{y}-\mathbf{x} \right)\cdot \mathbf{n}][\left(\mathbf{y}-\mathbf{x} \right)\cdot \mathbf{p}]^2}{\left|\mathbf{y}-\mathbf{x}\right|^2}
      \mathbf{n} d\mathbf{y} \\
 \nonumber&+\frac{C_\beta(D_{11}+D_{22})}{2m(\delta)} \int_{B_\delta (\mathbf{x}) \backslash \Omega} \lambda K \frac{[\left(\mathbf{y}-\mathbf{x} \right)\cdot \mathbf{n}]^3}{\left|\mathbf{y}-\mathbf{x}\right|^2}\mathbf{n} d\mathbf{y}-\frac{C_\beta\mu(D_{12}+D_{21})}{m(\delta)} \int_{B_\delta (\mathbf{x}) \backslash \Omega} K
     \frac{[\left(\mathbf{y}-\mathbf{x} \right)\cdot \mathbf{n}]}{\left|\mathbf{y}-\mathbf{x}\right|^2} 
      [\left(\mathbf{y}-\mathbf{x} \right)\cdot \mathbf{p}]^2\mathbf{p} d\mathbf{y}\\
 &-\frac{C_\beta[(\lambda+2\mu)D_{22}+\lambda D_{11}]}{2m(\delta)} \int_{B_\delta (\mathbf{x}) \backslash \Omega} K
     \frac{[\left(\mathbf{y}-\mathbf{x} \right)\cdot \mathbf{n}]}{\left|\mathbf{y}-\mathbf{x}\right|^2} 
      \left([\left(\mathbf{y}-\mathbf{x} \right)\cdot \mathbf{n}]^2-[\left(\mathbf{y}-\mathbf{x} \right)\cdot \mathbf{p}]^2\right)\mathbf{n} d\mathbf{y}\\
      =&\frac{C_\beta}{m(\delta)}\int_{E_\delta} \mu     K\frac{D_{11}[\left(\mathbf{y}-\mathbf{x}\right)\cdot\mathbf{p}]^3+D_{22}[\left(\mathbf{y}-\mathbf{x}\right)\cdot\mathbf{p}][\left(\mathbf{y}-\mathbf{x}\right)\cdot\mathbf{n}]^2}{\left|\mathbf{y}-\mathbf{x}\right|^2} \mathbf{p}d\mathbf{y}\\
      &+\frac{C_\beta}{m(\delta)}\int_{E_\delta} \mu     K\frac{(D_{12}+D_{21})[\left(\mathbf{y}-\mathbf{x}\right)\cdot\mathbf{p}][\left(\mathbf{y}-\mathbf{x}\right)\cdot\mathbf{n}]^2}{\left|\mathbf{y}-\mathbf{x}\right|^2} \mathbf{n}d\mathbf{y}+O(\delta)=O(\delta)\mathbf{1}.
 \end{align*}
\end{proof}

\begin{cor}
Given that $\Omega\in \mathbb{R}^d$ $(d=2)$ is a {$C^{3}$} domain and $|E_\delta|=0$, then the set of rigid deformations $\Pi$ is in the solution set of \eqref{eqn:probn} with $\fb=0$ and $\mathbf{T}=0$.
\end{cor}

We now investigate the consistency of the proposed mixed-type volume constraint formulation for general $\ub$, by considering the truncation estimate of the local solution. We denote $\mathbf{u}_\delta$ as the solution of the nonlocal problem \eqref{eqn:probn} and $\mathbf{u}_0$ as the solution of the corresponding linear elasticity problem:
 \begin{equation}\label{eq:u0}
 \left\{\begin{array}{ll}
 \mathcal{L}_0\ub=-\nabla\cdot(\lambda tr(\mathbf{E})\mathbf{I}+2\mu \mathbf{E})=\mathbf{f},\quad \text{where }\mathbf{E}=\dfrac{1}{2}(\nabla \mathbf{u}+(\nabla \mathbf{u})^T),&\quad \text{in }\Omega,\\
 \lambda tr(\mathbf{E})\mathbf{n}+2\mu \mathbf{E}\mathbf{n}= \mathbf{T},&\quad \text{on }\partial\Omega_N,\\
 \mathbf{u}= \mathbf{u}_D,&\quad \text{on }\omgbb_D.
 \end{array}\right.
 \end{equation}

 Denoting the truncation estimate $\mathbf{e}_\delta(\mathbf{x}):=\mcL_\delta\mathbf{u}_\delta(\mathbf{x})-\mcL_\delta \mathbf{u}_0(\mathbf{x})$ for $\xb\in\Omega\backslash\omgi_N$ and $\mathbf{e}_\delta(\mathbf{x}):=\mcL_{N\delta}\mathbf{u}_\delta(\mathbf{x})-\mcL_{N\delta} \mathbf{u}_0(\mathbf{x})$ for $\xb\in\omgi_N$, we may obtain the following bound for $\mathbf{e}_\delta$:

  \begin{thm}\label{thm:truncation}
  Assume that the local solution $\mathbf{u}_0\in C^1$, then $|\mathbf{e}_\delta|=O(\delta^2)$ for $\mathbf{x}\in \Omega\backslash\omgi_N$ and $|\mathbf{e}_\delta|=O(1)$ for $\mathbf{x}\in \omgi_N$.
 \end{thm}
 
 \begin{proof}
 For $\mathbf{x}\in \Omega\backslash\omgi_N$, from  $\mathbf{e}_\delta=\fb(\xb)-\mcL_{\delta} \mathbf{u}_0(\mathbf{x})=\mcL_{0} \mathbf{u}_0(\mathbf{x})-\mcL_{\delta} \mathbf{u}_0(\mathbf{x})$, the bound of $\mathbf{e}_\delta$ may be obtained via Taylor expansion of $\ub_0$ following a similar derivation as in \cite{You_2019}.  Denoting $u_{0p}$, $u_{0n}$ as the components of $\ub_0$ along the tangential and normal directions, respectively, for $\mathbf{x}\in\omgi_N$, with the calculation in \eqref{eqn:thetataylor} we have $\theta(\mathbf{u}_0)-\nabla\cdot\mathbf{u}_0=O(\delta)$. Note that the tangential and normal components of the traction load satisfies
$$T_p=\mu\left(\dfrac{\partial u_{0p}}{\partial \mathbf{n}}+\dfrac{\partial u_{0n}}{\partial \mathbf{p}}\right),\quad T_n=\lambda\nabla\cdot\ub_0+2\mu \dfrac{\partial u_{0n}}{\partial \mathbf{n}},$$
%and
%$$\frac{1}{m(\delta)}\int_{B_\delta(\xb)\cap\Omega}K[(\yb-\xb)\cdot\mathbf{n}]d\xb \leq C|B_\delta(\xb)\cap\Omega|/\delta=O(\delta),$$
and with the Taylor expansion of $\ub_0$, for $\mathbf{x}\in\omgi_N$ we have
\begin{align*}
\nonumber \mathbf{e}_\delta =&(L_0\mathbf{u}_0-L_{N\delta} \mathbf{u}_0)+(\mathbf{f}_{N\delta}-\mathbf{f})=-\dfrac{1}{2}\nabla\cdot(\lambda tr(\nabla \mathbf{u}_0+(\nabla \mathbf{u}_0)^T)\mathbf{I}+2\mu(\nabla \mathbf{u}_0+(\nabla \mathbf{u}_0)^T))\\
&+\frac{C_\alpha\left(\lambda - \mu\right)}{m(\delta)}  \int_{B_\delta (\mathbf{x}) \cap \Omega}  K
     \left(\mathbf{y}-\mathbf{x} \right)\left(\nabla\cdot\mathbf{u}_0(\mathbf{x}) + \nabla\cdot\mathbf{u}_0(\mathbf{y}) \right) d\mathbf{y}\\
   \nonumber&+\frac{C_\beta}{m(\delta)}\int_{B_\delta (\mathbf{x}) \cap \Omega} \mu
     K\frac{\left(\mathbf{y}-\mathbf{x}\right)\otimes\left(\mathbf{y}-\mathbf{x}\right)}{\left|\mathbf{y}-\mathbf{x}\right|^2}  \left(\mathbf{u}_0(\mathbf{y}) - \mathbf{u}_0(\mathbf{x}) \right) d\mathbf{y}
    \nonumber+\frac{2C_\alpha\nabla\cdot\ub_0(\mathbf{x})}{m(\delta)}  \int_{B_\delta (\mathbf{x}) \backslash \Omega} 
     \left(\lambda - \mu\right) K 
     \left(\mathbf{y}-\mathbf{x} \right) d\mathbf{y}\\
\nonumber&+\frac{C_\beta\nabla\cdot\ub_0(\mathbf{x})}{2m(\delta)}\int_{B_\delta (\mathbf{x}) \backslash\Omega}(\lambda+2\mu) K
     \frac{[\left(\mathbf{y}-\mathbf{x} \right)\cdot \mathbf{n}][\left(\mathbf{y}-\mathbf{x} \right)\cdot \mathbf{p}]^2}{\left|\mathbf{y}-\mathbf{x}\right|^2}
      \mathbf{n} d\mathbf{y} \\
 \nonumber&-\frac{C_\beta\nabla\cdot\ub_0(\mathbf{x})}{2m(\delta)} \int_{B_\delta (\mathbf{x}) \backslash \Omega} \lambda K   \frac{[\left(\mathbf{y}-\mathbf{x} \right)\cdot \mathbf{n}]^3}{\left|\mathbf{y}-\mathbf{x}\right|^2}\mathbf{n} d\mathbf{y}
\nonumber+\frac{C_\beta}{m(\delta)} \int_{B_\delta (\mathbf{x}) \backslash \Omega} K
     \frac{[\left(\mathbf{y}-\mathbf{x} \right)\cdot \mathbf{n}]}{\left|\mathbf{y}-\mathbf{x}\right|^2} 
      [\left(\mathbf{y}-\mathbf{x} \right)\cdot \mathbf{p}]^2[{T}_p(\bar{\mathbf{x}})\mathbf{p}] d\mathbf{y}\\
 &+\frac{C_\beta}{2m(\delta)} \int_{B_\delta (\mathbf{x}) \backslash \Omega} K     \frac{[\left(\mathbf{y}-\mathbf{x} \right)\cdot \mathbf{n}]}{\left|\mathbf{y}-\mathbf{x}\right|^2} 
      \left([\left(\mathbf{y}-\mathbf{x} \right)\cdot \mathbf{n}]^2-[\left(\mathbf{y}-\mathbf{x} \right)\cdot \mathbf{p}]^2\right)[{T}_n(\bar{\mathbf{x}})\mathbf{n}] d\mathbf{y}+\mcO(1)\mathbf{1}=\mcO(1)\mathbf{1}.
\end{align*}
 \end{proof}

\begin{remark}
From Thm. \ref{thm:geo} in the next section, we will see that given the possible numerical error from boundary approximations in the meshfree formulation, the $O(1)$ truncation estimate of $\mathbf{e}_\delta$ for $\xb\in\omgi_N$ is of optimal in M-convergence tests.
\end{remark}

\begin{remark}
To theoretically show the $L^2$ convergence of $\ub_\delta$ to the local limit $\ub_0$, a nonlocal Poincare-Korn's inequality would be required which will be addressed in the future work. In this work we demonstrate the asymptotic convergence rate with numerical examples in Section \ref{sec:num}, where a first order convergence $O(\delta)$ is observed for $\vertii{\ub_\delta-\ub_0}_{L^2(\Omega)}$, which indicates that the $O(1)$ truncation estimate in $\omgi_N$ is sufficient to obtain asymptotic convergence when $\delta=O(h)$. A similar phenomenon was also observed on the Neumann constraint nonlocal diffusion problem in \cite{You_2019}.

%To see the relation between the truncation estimate in $\omgi_N$ and the asymptotic convergence of nonlocal solutions in $L^2$, we take a nonlocal diffusion equation in \cite{You_2019} for simple demonstration. Denoting $u_\delta$, $u_0$ as the solutions from nonlocal diffusion and the corresponding local diffusion problems, respectively, we assume that $u_\delta=u_0$ in $\real^2\backslash\Omega$, and the truncation estimate $e_\delta$ satisfies $e_\delta=O(1)$ for $\xb\in\omgi_N$ and $e_\delta=0$ for $\xb\in\Omega\backslash\omgi_N$. When $\partial\Omega_N$ is sufficiently smooth, the nonlocal trace theorem in \cite{Wright2019} yields:
%\begin{align*}
%\vertii{u_\delta-u_0}^2_{L^2(\omgi)}&\leq C\delta \vertii{u_\delta-u_0}^2_{S^{Diff}_\delta(\Omega)}\leq C \delta A_\delta^{Diff}(u_\delta-u_0,u_\delta-u_0)\\
%&\leq C \delta \vertii{\mathcal{L}^{Diff}_{N\delta}u_\delta-\mathcal{L}^{Diff}_{N\delta}u_0}_{L^2(\omgi)}\vertii{u_\delta-u_0}_{L^2(\omgi)}\leq C\delta\vertii{u_\delta-u_0}^2_{L^2(\omgi)}. 
%\end{align*}
%Here $S^{Diff}_\delta(\Omega)$ is the nonlocal energy space, $A^{Diff}_\delta$ and $\mathcal{L}^{Diff}_{N\delta}$ are the nonlocal operators for the nonlocal diffusion problem in \cite{You_2019} which correspond to the operators defined in \eqref{eqn:bilinear} and \eqref{eq:newform1}, respectively. Therefore, we have $\vertii{u_\delta-u_0}_{L^2(\omgi)}\leq C\delta$. 
\end{remark}
 
 \section{Optimization-Based Meshfree Quadrature Rules}\label{sec:disc}
 
 \begin{figure}[!htb]\centering
 %\subfigure{\includegraphics[width=0.53\textwidth]{domain1.eps}}\quad
 \subfigure{\includegraphics[width=0.5\textwidth]{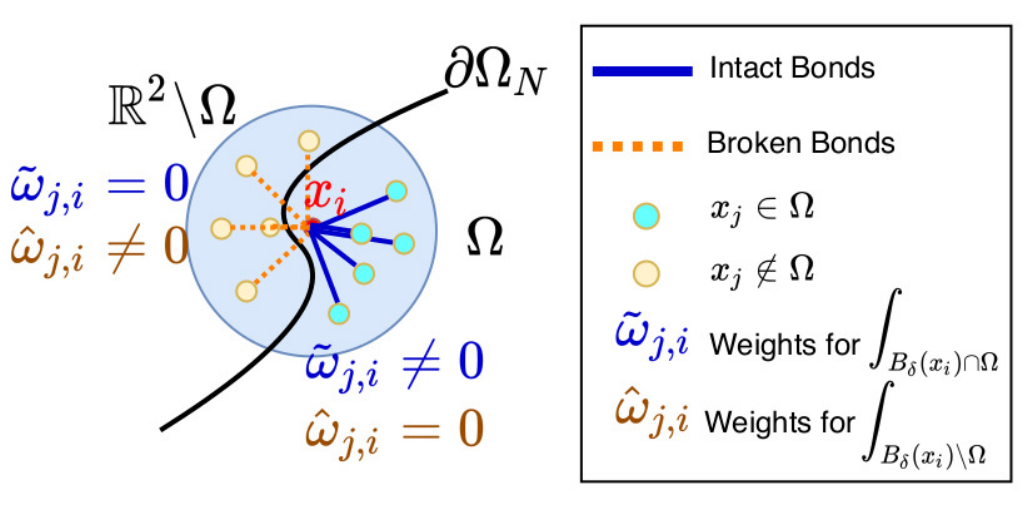}}
\qquad\qquad\subfigure{\includegraphics[width=0.25\textwidth]{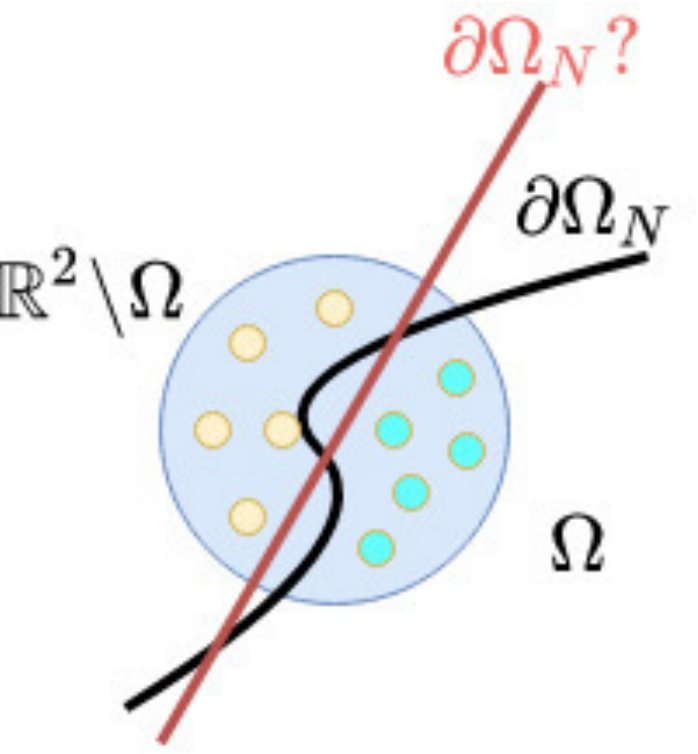}}
\caption{Illustration of neighbor points and bonds for $\xb_i\in\omgi_N$, where the yellow points represent $\xb_j\notin\Omega$ and blue points represent $\xb_j\in\Omega$. Left: An illustration of broken/intact bonds in the meshfree formulation. Right: An approximate boundary provided by breaking bonds in the meshfree formulation. The black and red curves indicate two possible boundaries represented by the same formulation in \eqref{eq:discreteNonlocElasticity2} and \eqref{eq:discreteNonlocDilitation2}. With a purely volumetric particle representation of the boundary, it is not possible to avoid first order truncation estimates from errors in the representation of the geometry.
 }
 \label{fig:Neumann}
\end{figure}

In this section, we introduce a strong-form particle discretizations of the state-based peridynamics introduced above. Discretizing the whole interaction region $\Omega\cup\omgbb$ by a collection of points $X_{h} = \{\mathbf{x}_i\}_{\{i=1,2,\cdots,N_p\}} \subset \Omega\cup\omgbb$, we aim to solve for the displacement $\ub_i\approx\ub(\xb_i)$ and the nonlocal dilitation $\theta_i\approx \theta(\xb_i)$ on all $\xb_i\in X_h$. We first characterize the distribution of collocation points as follows. Recall the definitions \cite{wendland2004scattered} of \textit{fill distance}
\begin{equation*}
h_{\chi_h,\Omega} = \underset{\ub_i \in \Omega\cup\omgb}{\sup}\, \underset{\xb_i \in \chi_h}{\min}||\xb_i - \xb_j||_2,
\end{equation*}
and \textit{separation distance}
\begin{equation*}
{q_{\chi_h} = \frac12 \underset{i \neq j}{\min} ||\xb_i - \xb_j||_2}.
\end{equation*}
For simplicity we drop subscripts and simply write $h$ and $q$. We assume that $\chi_h$ is \textit{quasi-uniform}, namely that there exists $c_{qu} > 0$ such that
\begin{equation*}
q_{\chi_h} \leq h_{\chi_h,\Omega} \leq c_{qu} q_{\chi_h}.
\end{equation*} 
To maintain an easily scalable implementation, in this paper we assume $\delta$ to be chosen such that the ratio $\frac{h}{\delta}$ is bound by a constant $M$ as $\delta \rightarrow 0$, restricting ourselves to the ``M-convergence'' scenario \cite{bobaru2009convergence}. 

As the first step, for the original LPS model \eqref{eq:nonlocElasticity}, we pursue a discretization through the following one point quadrature rule at a collection of collocation points $X_h$ \cite{silling_2010}:
\begin{align}
    \nonumber(\mathcal{L}_\delta^h \ub)_i:=&-\frac{C_\alpha}{m(\delta)} \sum_{\xb_j \in B_\delta (\mathbf{x}_i)} \left(\lambda - \mu\right) K_{ij} \left(\mathbf{x}_j-\mathbf{x}_i\right)\left(\theta_i + \theta_j \right) \omega_{j,i}\\
  &-  \frac{C_\beta}{m(\delta)}\sum_{\xb_j \in B_\delta (\mathbf{x}_i)} \mu K_{ij}\frac{\left(\mathbf{x}_j-\mathbf{x}_i\right)\otimes\left(\mathbf{x}_j-\mathbf{x}_i\right)}{\left|\mathbf{x}_j-\mathbf{x}_i\right|^2} \cdot \left(\mathbf{u}_j- \mathbf{u}_i \right) \omega_{j,i} = \mathbf{f}_i,\label{eq:discreteNonlocElasticity}
\end{align}
\begin{equation}
  \theta(\mathbf{x}) =  \frac{d}{m(\delta)} \sum_{\xb_j \in B_\delta (\mathbf{x}_i)}  K_{ij} \left(\mathbf{x}_j-\mathbf{x}_i\right) \cdot  \left(\mathbf{u}_j - \mathbf{u}_i \right) \omega_{j,i},\label{eq:discreteNonlocDilitation}
\end{equation}
where we adapt the notations $f(\xb_i) = f_i$ and $f(\xb_i,\xb_j)=f_{ij}$, and we specify $\left\{\omega_{j,i}\right\}$ as a to-be-determined collection of quadrature weights admitting interpretation as a measure associated with each collocation point $\xb_i$. We define in this section an optimization-based approach to defining these weights extending previous work \cite{trask2019asymptotically}, constructed to ensure consistency guarantees. Specifically, we seek quadrature weights for integrals supported on balls of the form
\begin{equation}
I[f] := \int_{B_\delta (\mathbf{x}_i)} f(\xb,\yb) d\yb \approx I_h[f] := \sum_{\xb_j \in B_\delta (\mathbf{x}_i)} f(\xb_i,\xb_j) \omega_{j,i}
\end{equation}
where we include the subscript $i$ in $\left\{\omega_{j,i}\right\}$ to denote that we seek a different family of quadrature weights for different subdomains $B_\delta(\mathbf{x}_i)$. We obtain these weights from the following optimization problem
\begin{align}\label{eq:quadQP}
  \underset{\left\{\omega_{j,i}\right\}}{\text{argmin}} \sum_{\xb_j \in B_\delta (\mathbf{x}_i)} \omega_{j,i}^2 \quad
  \text{such that}, \quad
  I_h[p] = I[p] \quad \forall p \in \mathbf{V}_h
\end{align}
where $\mathbf{V}_h$ denotes a Banach space of functions which should be integrated exactly. We refer to previous work \cite{trask2019asymptotically} for further information, analysis, and implementation details.

Provided the quadrature points are unisolvent over the desired reproducing space, this problem may be proven to have a solution by interpreting it as a generalized moving least squares (GMLS) problem \cite{leng2019asymptotically}. For certain choices of $V_h$, such as $m^{th}$-order polynomials, unisolvency holds under the assumptions that: the domain $\Omega$ satisfies a cone condition, the pointset $X_h \cap B_\delta(\bm{x}_i)$ is quasi-uniform, and $\delta$ is sufficiently large \cite{wendland2004scattered}.

In previous work we have provided truncation error estimates relating the quadrature error convergence rate to the order of singularity in the kernel. As discussed in \cite{trask2019asymptotically}, the key to obtaining these quadrature weights is that they may be evaluated analytically, either via analytic rules \cite{folland2001integrate} or the aid of symbolic integration software. In this work, we choose a reproducing space sufficient to integrate \eqref{eq:discreteNonlocElasticity} and \eqref{eq:discreteNonlocDilitation} exactly in the case where $\ub$ and $\theta$ are quadratic polynomials.

\begin{thm}\label{thm:Diri}
Let $\mathbf{V}_h = \left\{ q = \frac{p(\yb)}{|\yb-\xb|^3} \,|\, p \in P_5(\mathbb{R}^d) \text{ such that } \int_{B_\delta(\xb)} q(\yb) d \yb < \infty \right\}$ where $p \in P_5(\mathbb{R}^d)$ is the space of quintic polynomials, and assume $B_\delta(\xb) \subset \Omega\cup\omgbb_D$ and that the optimization problem \eqref{eq:quadQP} has a solution. Then the meshfree optimization-based quadrature (OBQ) approximations to \eqref{eq:nonlocElasticity} and \eqref{eqn:oritheta} are exact for $\ub \in (P_3(\mathbb{R}^d))^d$ and $\theta \in P_2(\mathbb{R}^d)$. Further, for $\ub \in C^3$ and $\theta \in C^2$ the truncation error for all nonlocal operators in \eqref{eq:nonlocElasticity} converge to its local limit with an $O(\delta^2)$ rate in the limit $\delta \rightarrow 0$.
\end{thm}

\begin{proof}
We prove only for the nonlocal gradient; the other operators follow similarly. Rewriting the gradient as $\int_{B_\delta (\mathbf{x})} \frac{|\xb-\yb|^2}{|\xb-\yb|^3}\left(\mathbf{y}-\mathbf{x} \right)\left(\theta(\mathbf{x}) + \theta(\mathbf{y}) \right) d\mathbf{y}$ and assuming $\theta \in P_2$, we obtain a component-wise form $\int_{B_\delta (\mathbf{x})} \frac{p(\yb)}{|\xb-\yb|^3} d\yb$ where $p \in P_5$, and thus the quadrature is exact via the equality constraint of \eqref{eq:quadQP}. To prove convergence, note that we may rewrite $\int_{B_\delta (\mathbf{x})} \frac{|\xb-\yb|^2}{|\xb-\yb|^3}\left(\mathbf{y}-\mathbf{x} \right)\left(\theta(\mathbf{y}) - \theta(\mathbf{x}) \right) d\mathbf{y}$ because constants are in the null-space of the nonlocal gradient. The proof then follows from Thm. 2.1 of \cite{trask2019asymptotically}, by approximating $\theta(\mathbf{y}) - \theta(\mathbf{x})$ via a third order Taylor series, converting to polar coordinates, and bounding terms.
\end{proof}

\begin{remark}
We have selected this particular choice of reproducing space so that the same quadrature weights may be used for all nonlocal operators. The cost of constructing the quadrature scales as $dim(\mathbf{V}_h)^3$, and significant savings may result by instead generating quadrature rules for each operator specifically. For example, one may obtain the same convergence by selecting $\mathbf{V}_h = \left\{ \frac{\yb-\xb}{|\yb-\xb|}\left( p(\yb)+p(\xb) \right) \,|\, p \in P_2(\mathbb{R}^d)\right\}$ in the nonlocal gradient. As solving \eqref{eq:quadQP} amounts to inverting a small dense matrix, we may expect a substantial $(dim(P_2)/dim(P_5))^3$ speed-up in this case ($7^3\times$speedup in $2D$ and $(56/10)^3\times$speedup in $3D$). Because the requisite optimization problems are amenable to fine-grained parallelism on GPUs using libraries such as the Compadre toolkit \cite{compadrev101}, we prefer in this work to use a single quadrature rule for all operators; in our implementation the cost of constructing the quadrature weights is negligible compared to solving the resultant stiffness matrix after discretization.
\end{remark}

\begin{remark}
  In many quadrature schemes it is desirable to enforce the positivity of quadrature weights (i.e. $\omega_{j,i}>0$). While we do not pursue this property in the current work, \eqref{eq:quadQP} may be modified to enforce this property via an inequality constraint. In the context of maximum principle preserving meshfree discretizations this has been considered \cite{seibold2008minimal}. We note that for quasi-uniform particle distributions, we have seen only a small number of negative quadrature weights, and that these are generally very small compared to the other positive weights.
\end{remark}  
 
We now apply the above quadrature rule to the LPS model with traction loads applied on a sharp boundary $\partial\Omega_N$. %Specifically, for $\xb_i\in\Omega\backslash\omgi_N$, we solve for $\ub_i$ and $\theta_i$ with  \eqref{eq:discreteNonlocElasticity} and \eqref{eq:discreteNonlocDilitation}. 
For $\xb_i\in\omgi_N$, we note that $B_\delta(\xb_i)\backslash\omg\neq\emptyset$. In the meshfree formulation, the boundary $\partial\omg_N$ is represented by breaking bonds between $\xb_i$ and $\xb_j\in B_\delta(\xb_i)\backslash\omg$, as demonstrated in Figure \ref{fig:Neumann}. For $\xb_j\in\Omega$, we denote the bond between $\xb_i$ and $\xb_j$ as ``intact'' and the change of displacement on material point $\xb_j$ may have an impact on the displacement at $\xb_i$. On the other hand, when $\xb_j\notin\Omega$, we consider the bonds between $\xb_i$ and $\xb_j$ as ``broken''. To discretize \eqref{eq:continuousNonlocDilitation3} and \eqref{eq:newform1}, the quadrature weights associated with intact bonds will be employed in the calculation of integrals inside $B_\delta(\xb_i)\cap\omg$ and the weights associated with broken bonds will be employed for integrals inside $B_\delta(\xb_i)\backslash\omg$. Particularly, we express the quadrature weights associated with intact bonds $\tilde{\omega}_{j,i}$ and the quadrature weights associated with broken bonds $\hat{\omega}_{j,i}$ in terms of the scalar mask $\gamma_{j,i}$:
\begin{align}\label{eqn:gamma}
     \gamma_{j,i} &= \begin{cases}
    1, \quad \text{if }\xb_j\in B_\delta(\xb_i)\cap\Omega,\\
    0, \quad \text{otherwise}, \\
    \end{cases}\quad \tilde{\omega}_{j,i}={\omega}_{j,i}\gamma_{j,i}, \quad \hat{\omega}_{j,i}={\omega}_{j,i}(1-\gamma_{j,i}).
\end{align}
Numerical quadrature of a given function $a(\xb)$ over $B_\delta(\xb_i)\cap\omg$ and $B_\delta(\xb_i)\backslash\omg$ may thus be calculated via
$$\int_{B_\delta(\xb_i)\cap\omg}a(\yb)d\yb\approx\sum_{\xb_j\in B_\delta(\xb_i)} \tilde{\omega}_{j,i} a(\xb_j),\qquad \int_{B_\delta(\xb_i)\backslash\omg}a(\yb)d\yb\approx\sum_{\xb_j\in B_\delta(\xb_i)} \hat{\omega}_{j,i} a(\xb_j).$$
This process is consistent with how damage is typically induced in bond-based peridynamics, such as the prototype microelastic brittle model \cite{silling_2005_2}.
Applying the above formulation in \eqref{eq:continuousNonlocDilitation3} and \eqref{eq:newform1} we propose the following meshfree scheme:
{\begin{align}
   \nonumber &(\mathcal{L}_{N\delta}^h \ub)_i:=\sum_{\xb_j \in B_\delta (\mathbf{x}_i)}\frac{K_{ij}}{m(\delta)}\left[\left(-{C_\alpha} \left(\lambda - \mu\right)  \left(\mathbf{x}_j-\mathbf{x}_i\right)\left(\theta_i + \theta_j \right)
  -{C_\beta }\mu \frac{\left(\mathbf{x}_j-\mathbf{x}_i\right)\otimes\left(\mathbf{x}_j-\mathbf{x}_i\right)}{\left|\mathbf{x}_j-\mathbf{x}_i\right|^2} \cdot \left(\mathbf{u}_j- \mathbf{u}_i \right)\right) \tilde{\omega}_{j,i}\right.\\
  \nonumber&~~+\left(-{2C_\alpha\left(\lambda - \mu\right)}     \left(\mathbf{x}_j-\mathbf{x}_i \right) -\frac{C_\beta(\lambda+2\mu)\mathbf{n}_i }{2}
     \frac{[\left(\mathbf{x}_j-\mathbf{x}_i \right)\cdot \mathbf{n}_i][\left(\mathbf{x}_j-\mathbf{x}_i \right)\cdot \mathbf{p}_i]^2}{\left|\mathbf{x}_j-\mathbf{x}_i\right|^2}\left.+\frac{C_\beta\lambda\mathbf{n}_i }{2}  \frac{[\left(\mathbf{x}_j-\mathbf{x}_i \right)\cdot \mathbf{n}_i]^3}{\left|\mathbf{x}_j-\mathbf{x}_i\right|^2}\right)\theta_i\hat{\omega}_{j,i}\right]\\
\nonumber&= \fb(\xb_i)+\sum_{\xb_j \in B_\delta (\mathbf{x}_i)}\frac{K_{ij}\hat{\omega}_{j,i}}{m(\delta)}\left({C_\beta T_p(\bar{\mathbf{x}}_i)\mathbf{p}_i}     \frac{[\left(\mathbf{x}_j-\mathbf{x}_i \right)\cdot \mathbf{n}_i]}{\left|\mathbf{x}_j-\mathbf{x}_i\right|^2} 
      [\left(\mathbf{x}_j-\mathbf{x}_i \right)\cdot \mathbf{p}_i]^2\right.\\
&~~+\left.\frac{C_\beta T_n(\bar{\mathbf{x}}_i)\mathbf{n}_i}{2}\frac{[\left(\mathbf{x}_j-\mathbf{x}_i \right)\cdot \mathbf{n}_i]}{\left|\mathbf{x}_j-\mathbf{x}_i\right|^2} 
      \left([\left(\mathbf{x}_j-\mathbf{x}_i \right)\cdot \mathbf{n}_i]^2-[\left(\mathbf{x}_j-\mathbf{x}_i \right)\cdot \mathbf{p}_i]^2\right)\right):=(\fb_{N\delta})_i,\label{eq:discreteNonlocElasticity2}
\end{align}
\begin{equation}
\theta_i =\frac{d}{m(\delta)} \sum_{\xb_j \in B_\delta (\mathbf{x}_i)} K_{ij} \left(\mathbf{x}_j-\mathbf{x}_i\right) \cdot \mathbf{M}_i\cdot \left(\mathbf{u}_j - \mathbf{u}_i \right) \tilde{\omega}_{j,i},\label{eq:discreteNonlocDilitation2}
\end{equation}
where 
\begin{equation}
\mathbf{M}_i:=\left[\dfrac{d}{m(\delta)}\sum_{\xb_j \in B_\delta (\mathbf{x}_i)} K_{ij}(\mathbf{x}_j-\mathbf{x}_i)\otimes(\mathbf{x}_j-\mathbf{x}_i)\tilde{\omega}_{j,i}\right]^{-1}.
\end{equation}}
Note that although we have shown in Thm. \ref{lem:M} that $\mathbf{M}(\xb_i)$ exists when $\partial\Omega$ is sufficiently smooth, the numerical evaluation of the correction tensor further requires that $\underset{\xb_j \in B_\delta (\mathbf{x}_i)}{\sum} K_{ij}(\mathbf{x}_j-\mathbf{x}_i)\otimes(\mathbf{x}_j-\mathbf{x}_i)\tilde{\omega}_{j,i}$ be invertible. This is true as long as there are at least $d$ non-collinear bonds within the horizon. In some settings, such as violent dynamic fracture, for a given particle all bonds may break, leaving an isolated particle. In this case the matrix inverse $\mathbf{M}_i^{-1}$ may be replaced with a pseudo-inverse $\mathbf{M}_i^+$ to improve robustness of the scheme.

% In this work, we choose a reproducing space sufficient to integrate Equations \ref{eq:discreteNonlocElasticity} and \ref{eq:discreteNonlocDilitation} exactly in the case where $\bm{u}$ and $\theta$ are quadratic polynomials. Specifically, we take $\mathbf{V}_h = \left\{ q = \frac{p(\bm{y})}{|\bm{y}-\bm{x}|^3} \,|\, p \in P_5(\mathbb{R}^d) \text{ such that } \int_{B_\delta(\bm{x})} q(\bm{y}) d \bm{y} < \infty \right\}$. where $p \in P_5(\mathbb{R}^d)$ is the space of quintic polynomials, and assume $B_\delta(x) \subset \Omega$ and that Eqn. \ref{eq:quadQP} has a solution. 

Note that in the traction-type boundary condition formulation \eqref{eq:newform1}, the unit normal vector $\mathbf{n}(\bar{\xb})$ is required, and the unit tangential vector $\mathbf{p}(\bar{\xb})$ can then be calculated as the orthogonal unit vector of $\mathbf{n}(\bar{\xb})$. However, in realistic settings the analytical form of $\mathbf{n}(\bar{\xb})$ is often unavailable. To approximate the normal vector at $\bar{\xb}_i$ for each $\xb_i$, we note that
$$\nb(\bar{\xb})\approx -\dfrac{\int_{B_\delta (\mathbf{x}) \cap \Omega}(\yb-\xb) d\yb}{\vertii{\int_{B_\delta (\mathbf{x}) \cap \Omega}(\yb-\xb) d\yb}}.$$
Therefore numerically we calculate the normal direction as $\mathbf{n}_i$ as
\begin{equation}\label{eqn:approx_n}
\nb_i = -\dfrac{\underset{{\xb_j \in B_\delta (\mathbf{x}_i)}}{\sum}(\xb_j-\xb_i)\tilde{\omega}_{j,i} }{\vertii{\underset{\xb_j \in B_\delta (\mathbf{x}_i)}{\sum}(\xb_j-\xb_i)\tilde{\omega}_{j,i}}},    
\end{equation}
and the tangential vector $\mathbf{p}_i$ is calculated as the orthogonal direction to $\nb_i$.

Note that the formulation \eqref{eqn:approx_n} provides a practical approximation of the unit normal vector for each $\xb\in\omgi_N$ instead of each $\bar{\xb}\in\partial\Omega_N$, which therefore induces possible numerical errors in \eqref{eq:discreteNonlocElasticity2}. Moreover, in \eqref{eq:discreteNonlocElasticity2} and \eqref{eq:discreteNonlocDilitation2} we only solve for $\ub$ and $\theta$ in $\Omega$, which is equivalent to breaking any bond intersecting the Neumann boundary $\partial\Omega_N$. We highlight that quadrature weights $\omega_{j,i}$ are computed in the reference configuration before bonds are broken, and therefore no remeshing or calculation of quadrature weights will be required as fracture progresses. This property offers an efficient and sharp treatment of boundary geometry which may be easily implemented in popular particle mechanics codes. However, as illustrated in Figure \ref{fig:Neumann}, such a numerical approximation for the boundary shape $\partial\Omega_N$ introduces an $O(h)$ error to the boundary shape and correspondingly to the provided traction load $\mathbf{T}$ on $\partial\Omega_N$, and therefore errors in \eqref{eq:discreteNonlocElasticity2} and \eqref{eq:discreteNonlocDilitation2}.

To characterize the resulting numerical error, in the following we consider an equivalent problem: a perturbed traction load $\hat{\mathbf{T}}(\xb)$ is provided on $\partial\Omega_N$, i.e., there exists a constant $C$ which is independent of $h$ and $\delta$, such that 
$$|\hat{\mathbf{T}}(\xb)-{\mathbf{T}}(\xb)|\leq Ch,\quad \forall \xb\in\partial\Omega_N.$$
Moreover, due to the presumed perturbation of the geometry and the numerical error in \eqref{eqn:approx_n}, for $\xb\in\omgi_N$ we assume that the normal and tangential directions are also perturbed such that $\hat{\mathbf{n}}-\mathbf{n}=O(h)$ and $\hat{\mathbf{p}}-\mathbf{p}=O(h)$. 
%Moreover, we assume that $\hat{\mathbf{n}}$ and $\hat{\mathbf{p}}$ still satisfy the following approximated symmetry conditions:
%\begin{align}
%\nonumber&\frac{C_\beta}{m(\delta)} \int_{B_\delta (\mathbf{x}) \cap {\Omega}} K(\left|\mathbf{y}-\mathbf{x}\right|) \frac{[\left(\mathbf{y}-\mathbf{x} \right)\cdot \hat{\mathbf{n}}]^2[\left(\mathbf{y}-\mathbf{x} \right)\cdot \hat{\mathbf{p}}]}{\left|\mathbf{y}-\mathbf{x}\right|^2} = O(h),\\
%&\frac{C_\beta}{m(\delta)} \int_{B_\delta (\mathbf{x}) \cap {\Omega}} K(\left|\mathbf{y}-\mathbf{x}\right|) \frac{[\left(\mathbf{y}-\mathbf{x} \right)\cdot \hat{\mathbf{p}}]^3}{\left|\mathbf{y}-\mathbf{x}\right|^2} = O(h).\label{eq:symassume}
%\end{align}
With the perturbed traction loads and perturbed unit vectors specified above, we denote the (perturbed) nonlocal operator defined in \eqref{eq:newform1} as $\hat{\mathcal{L}}_{N\delta}$, the (perturbed) nonlocal dilitation defined in \eqref{eq:continuousNonlocDilitation3} as $\hat{\theta}$, and
\begin{align}
\nonumber\hat{\fb}_{N\delta}(\xb):= & \fb(\xb)+\frac{C_\beta}{m(\delta)} \int_{B_\delta (\mathbf{x}) \backslash {\Omega}} K(\left|\mathbf{y}-\mathbf{x}\right|) 
     \frac{[\left(\mathbf{y}-\mathbf{x} \right)\cdot \hat{\mathbf{n}}]}{\left|\mathbf{y}-\mathbf{x}\right|^2} 
      [\left(\mathbf{y}-\mathbf{x} \right)\cdot \hat{\mathbf{p}}]^2[\hat{T}_p({\bar{\mathbf{x}}})\hat{\mathbf{p}}] d\mathbf{y}\\
&+\frac{C_\beta}{2m(\delta)} \int_{B_\delta (\mathbf{x}) \backslash {\Omega}} K(\left|\mathbf{y}-\mathbf{x}\right|) 
     \frac{[\left(\mathbf{y}-\mathbf{x} \right)\cdot \hat{\mathbf{n}}]}{\left|\mathbf{y}-\mathbf{x}\right|^2} 
      \left([\left(\mathbf{y}-\mathbf{x} \right)\cdot \hat{\mathbf{n}}]^2-[\left(\mathbf{y}-\mathbf{x} \right)\cdot \hat{\mathbf{p}}]^2\right)[\hat{T}_n({\bar{\mathbf{x}}})\hat{\mathbf{n}}] d\mathbf{y}.\label{eq:hatf}
\end{align}
We then provide the truncation estimate corresponding to the above perturbations as follows.
\begin{thm}\label{thm:geo}
Assume that $\ub\in C^1$, %is a solution of \eqref{eq:nonlocElasticity}-\eqref{eqn:oritheta}, $\hat{\Omega}$ 
$\Omega$ is $C^3$ smooth,  $\hat{\mathbf{T}}$ is a perturbed approximation of the local traction load as defined in \eqref{eq:u0}:
$$\hat{T}_p(\bar{\xb})=\mu\left(\dfrac{\partial u_{p}}{\partial {\mathbf{n}}}({{\xb}})+\dfrac{\partial u_{n}}{\partial {\mathbf{p}}}({{\xb}})\right)+O(h),\;\hat{T}_n(\bar{\xb})= \lambda\nabla\cdot\ub(\xb)+2\mu \dfrac{\partial u_n}{\partial {\mathbf{n}}}(\xb)+O(h),\;\forall\xb\in\partial\Omega_N,$$
where $u_{p}$, $u_{n}$ are the components of $\ub$ along the non-perturbed tangential and normal directions respectively. The truncation estimates from perturbed $\mathbf{T}$, $\mathbf{n}$ and $\mathbf{p}$ for nonlocal operator in \eqref{eq:newform1} is bounded by $O(1)$ under the M-convergence condition. Specifically,
$$|\hat{\mathcal{L}}_{N\delta}\ub-\hat{\fb}_{N\delta}-{\mathcal{L}}_{N\delta}\ub+{\fb}_{N\delta}|=O(1).$$
\end{thm}
\begin{proof}
%The bound for \eqref{eq:continuousNonlocDilitation3} is a corollary of Thm. \ref{thm:thetaconsis}. Denoting $K(|\yb-\xb|)$ as $K$ for simplicity, %Denote 
%$$\hat{\sigma}_p(\ub):=\mu\left(\dfrac{\partial u_{p}}{\partial \hat{\mathbf{n}}}({{\xb}})+\dfrac{\partial u_{n}}{\partial \hat{\mathbf{p}}}({{\xb}})\right),\,\hat{\sigma}_n(\ub):= \lambda\nabla\cdot\ub(\xb)+2\mu \dfrac{\partial u_n}{\partial \hat{\mathbf{n}}}(\xb) ,$$
%$${\sigma}_p(\ub):=\mu\left(\dfrac{\partial u_{p}}{\partial {\mathbf{n}}}({{\xb}})+\dfrac{\partial u_{n}}{\partial {\mathbf{p}}}({{\xb}})\right),\,{\sigma}_n(\ub):= \lambda\nabla\cdot\ub(\xb)+2\mu \dfrac{\partial u_n}{\partial {\mathbf{n}}}(\xb),$$
%with Taylor expansion and $h=O(\delta)$, we note that
%$$\left|T_p(\bar{\xb})-\hat{\sigma}_p(\ub)\right|=O(h),\quad\left|T_n(\bar{\xb})-\hat{\sigma}_n(\ub)\right|=O(h).$$
%Note that the area of $B_\delta(\xb)\cap\triangle\Omega$ has $|B_\delta(\xb)\cap\triangle\Omega|\leq Ch\delta$. 
%Without loss of generality, in the following we assume that $\Omega\subset\hat{\Omega}$, and the proof for the general situation can be similarly derived. 
From the definition of nonlocal dilitation in \eqref{eq:continuousNonlocDilitation3}, we note that $\theta$ is not influenced by the perturbations on traction loads and normal/tangential directions, i.e., $\hat{\theta}=\theta$. To obtain the bound for the nonlocal operator in \eqref{eq:newform1} we separate the truncation estimate into two parts, the part from perturbation on $\mathbf{T}$ and the part induced by the perturbation of $\mathbf{n}$ and $\mathbf{p}$:
$$\hat{\mathcal{L}}_{N\delta}\ub-\hat{\fb}_{N\delta}-{\mathcal{L}}_{N\delta}\ub+{\fb}_{N\delta}=H_1+H_2.$$
\begin{align*}
H_1:=&\frac{C_\beta}{m(\delta)} \int_{B_\delta (\mathbf{x}) \backslash {\Omega}} K\frac{[\left(\mathbf{y}-\mathbf{x} \right)\cdot \hat{\mathbf{n}}]}{\left|\mathbf{y}-\mathbf{x}\right|^2} 
      [\left(\mathbf{y}-\mathbf{x} \right)\cdot \hat{\mathbf{p}}]^2(T_p(\bar{\xb})-\hat{T}_p(\bar{\xb}))\hat{\mathbf{p}} d\mathbf{y}\\
\nonumber &+\frac{C_\beta}{2m(\delta)} \int_{B_\delta (\mathbf{x}) \backslash{\Omega}} K\frac{[\left(\mathbf{y}-\mathbf{x} \right)\cdot \hat{\mathbf{n}}]}{\left|\mathbf{y}-\mathbf{x}\right|^2} 
      \left([\left(\mathbf{y}-\mathbf{x} \right)\cdot \hat{\mathbf{n}}]^2-[\left(\mathbf{y}-\mathbf{x} \right)\cdot \hat{\mathbf{p}}]^2\right)(T_n(\bar{\xb})-\hat{T}_n(\bar{\xb}))\hat{\mathbf{n}} d\mathbf{y}\\
\end{align*}
\begin{align*}
H_2:=&\frac{C_\beta(\lambda+2\mu)\theta(\mathbf{x})}{2m(\delta)}\int_{B_\delta (\mathbf{x}) \backslash\Omega} \frac{K}{\left|\mathbf{y}-\mathbf{x}\right|^2} \left([\left(\mathbf{y}-\mathbf{x} \right)\cdot {\mathbf{n}}][\left(\mathbf{y}-\mathbf{x} \right)\cdot {\mathbf{p}}]^2
      {\mathbf{n}}-[\left(\mathbf{y}-\mathbf{x} \right)\cdot \hat{\mathbf{n}}][\left(\mathbf{y}-\mathbf{x} \right)\cdot \hat{\mathbf{p}}]^2\hat{\mathbf{n}}\right) d\mathbf{y} \\
 \nonumber&+\frac{C_\beta\lambda\theta(\mathbf{x})}{2m(\delta)} \int_{B_\delta (\mathbf{x}) \backslash\Omega}  \frac{K}{\left|\mathbf{y}-\mathbf{x}\right|^2}  \left({[\left(\mathbf{y}-\mathbf{x} \right)\cdot \hat{\mathbf{n}}]^3}\hat{\mathbf{n}}-{[\left(\mathbf{y}-\mathbf{x} \right)\cdot {\mathbf{n}}]^3}{\mathbf{n}}\right) d\mathbf{y}\\
& +\frac{C_\beta}{m(\delta)} \int_{B_\delta (\mathbf{x}) \backslash\Omega} \frac{K}{\left|\mathbf{y}-\mathbf{x}\right|^2} \left({[\left(\mathbf{y}-\mathbf{x} \right)\cdot {\mathbf{n}}]}[\left(\mathbf{y}-\mathbf{x} \right)\cdot {\mathbf{p}}]^2[T_p(\bar{\xb}){\mathbf{p}}]-{[\left(\mathbf{y}-\mathbf{x} \right)\cdot \hat{\mathbf{n}}]}[\left(\mathbf{y}-\mathbf{x} \right)\cdot \hat{\mathbf{p}}]^2[{T}_p(\bar{\xb})\hat{\mathbf{p}}]\right) d\mathbf{y}\\
\nonumber &+\frac{C_\beta}{2m(\delta)} \int_{B_\delta (\mathbf{x}) \backslash\Omega} \frac{K}{\left|\mathbf{y}-\mathbf{x}\right|^2} \left({[\left(\mathbf{y}-\mathbf{x} \right)\cdot {\mathbf{n}}]}\left([\left(\mathbf{y}-\mathbf{x} \right)\cdot {\mathbf{n}}]^2-[\left(\mathbf{y}-\mathbf{x} \right)\cdot {\mathbf{p}}]^2\right)[{T}_n(\bar{\xb}){\mathbf{n}}]\right. \\   
\nonumber &\left.-{[\left(\mathbf{y}-\mathbf{x} \right)\cdot \hat{\mathbf{n}}]}\left([\left(\mathbf{y}-\mathbf{x} \right)\cdot \hat{\mathbf{n}}]^2-[\left(\mathbf{y}-\mathbf{x} \right)\cdot \hat{\mathbf{p}}]^2\right)[{T}_n(\bar{\xb})\hat{\mathbf{n}}]\right) d\mathbf{y}.
\end{align*}
With a similar technique as in Thm. \ref{thm:truncation} and with the $\delta/h=M$ assumption in M-convergence tests, we can show that each term in $H_1$ and $H_2$ is bounded by $O(1)$.
%, with Taylor expansion and the approximated symmetry of operators in \eqref{eq:symassume} we have
\end{proof}

\begin{remark}
From the proof of Thm. \ref{thm:geo}, we can see that in M-convergence tests either an $O(h)$ error on the provided traction load or an $O(h)$ error on the approximated unit normal and tangential vectors will induce an $O(1)$ truncation estimate in \eqref{eq:newform1} for $\xb\in\omgi_N$, which is of the same order generated by the proposed traction-type boundary condition formulation as discussed in Thm. \ref{thm:truncation}. Therefore, when using meshfree formulations and the broken bond techniques to induce damage as in \cite{silling_2005_2}, the proposed nonlocal traction loading formulation is of optimal asymptotic M-convergence rate to its local limit.
\end{remark}

To sum up, with the meshfree discretization described above and the optimization-based quadrature weights $\omega_{j,i}$, $i=1,\cdots,N_p$, we solve for the displacement $\ub(\xb_i)$ and nonlocal dilatation $\theta_i$ from:
\begin{equation}\label{eqn:disc_static}
    \mathbb{K}\eta=\mathbf{F}.
\end{equation}
Here $\eta$ is the vector of unknowns organized as follows:
$$\eta=[u_1,\cdots,u_{DOF},v_1,\cdots,v_{DOF},\theta_1,\cdots,\theta_{DOF}]^T.$$
$DOF=\#{i:\xb_i\in\Omega}$ is the total number of material points to be solved, and $u$, $v$ are the components of displacement such that $\ub_i=[u_i,v_i]$. $\mathbb{K}$ is a $3DOF\times3DOF$ stiffness matrix. The right hand side $\mathbf{F}$ is organized following a similar way as for $\eta$.

 \section{Numerical Verification and Asymptotic Compatibility}\label{sec:num}
 
In this section we numerically verify the approach by investigating accuracy when recovering analytic solutions in the M-convergence limit with mixed boundary conditions. We consider: linear patch tests, smooth manufactured solutions, analytical solutions to curvilinear surface loading problems, and analytical solutions to linearly elastic composites. For each case, we consider various combinations of Dirichlet and traction-type boundary conditions, exploring also the effect of reduced regularity on the traction problem by considering both Lipschitz and smooth boundaries. For each case we consider refinements of both Cartesian grids with mesh spacing $h$, and nonuniform grids generated by perturbing the Cartesian grids with a uniformly distributed random vector field  $(\Delta x,\Delta y)$, $\Delta x, \Delta y \sim \mathcal{U}[-0.2h,0.2h]$. For the sake of brevity we report the formal convergence study in \ref{app:convergence}, but summarize the setup and main conclusions for each case below, particularly focusing on whether optimal first order convergence in $\delta$ is realized as $\delta \rightarrow 0$, or if a lack of boundary regularity leads to suboptimal convergence. In all cases considered, the scheme does provide AC convergence as $\delta,h\rightarrow 0$.

\subsection{Linear Patch Test}\label{sec:patch}

We consider as linear patch test the displacement
$$\ub(x,y)=(3x+2y,-x+2y)$$
on a square domain $\Omega=[-\pi/2,\pi/2] \times [-\pi/2,\pi/2]$, with three different boundary conditions:
\begin{enumerate}
    \item Full Dirichlet-type boundary conditions: $\partial\Omega_D=\partial\Omega$;
    \item Mixed boundary conditions with traction loads applied on a straight line: $\partial\Omega_N=\{(x,\pi/2)|x\in[-\pi/2,\pi/2]\}$ and $\partial\Omega_D=\partial\Omega\backslash\partial\Omega_N$;
    \item Mixed boundary conditions with traction loads applied on corner: $\partial\Omega_N=\{(x,\pi/2)|x\in[-\pi/2,\pi/2]\}\cup\{(\pi/2,y)|y\in[-\pi/2,\pi/2]\}$ and $\partial\Omega_D=\partial\Omega\backslash\partial\Omega_N$.
\end{enumerate}
Note that in linear patch tests, the local and nonlocal solutions coincide. On settings 2 and 3, a traction-type boundary condition
\begin{equation}\label{eqn:T}
\mathbf{T}=\left[\begin{array}{cc}
5\lambda+6\mu&\mu\\
\mu&5\lambda+4\mu\\
\end{array}\right]\mathbf{n}
\end{equation}
is applied on the interface $\partial\Omega_N$, with material parameters following the plane strain assumption:
$$\lambda=K\nu/((1+\nu)(1-2\nu)),\mu=K/(2(1+\nu)),$$
for Young's modulus $K=1$. Two values of Poisson ratio $\nu=0.3$ and $0.49$ are investigated which correspond to compressible and nearly-incompressible materials, respectively. To demonstrate independence of $M$-convergence rate to choice of $M$, we consider both $\delta=3.5h$ and $\delta=3.9h$. Note that in problems with the boundary condition setting 3, when $\xb$ is close to the corner the projection point $\overline{\xb}$ is possibly ill-defined and therefore induces ambiguity of definition on $\mathbf{T}(\overline{\xb})$. To resolve this possible issue, in setting 3 we define $\mathbf{T}(\overline{\xb})$ following \eqref{eqn:T} where $\nb$ is the numerical approximation of normal direction following \eqref{eqn:approx_n}.

\textbf{Uniform discretization}: With settings 1 and 2, we observe that the numerical solution passes the patch test to within machine precision. Note that in setting 2, $\partial\Omega_N$ consists of a straight line and therefore $B_\delta(\xb)\backslash\Omega$ is symmetric with respect to $\mathbf{n}$, and the numerical result is consistent with Thm. \ref{thm:linearpatch}. In setting 3, $B_\delta(\xb)\backslash\Omega$ is not symmetric when $\xb$ is close to the corner, and the numerical solution only passes the linear patch test approximately. In Figure \ref{fig:patch_corner_s}, we compare $L^2(\omg)$ error vs. $h$ for displacement and dilitation to demonstrate first order AC convergence for both $\ub$ and $\theta$, independent of $\delta/h$ and $\nu$.

\textbf{Non-uniform discretization}: For a randomly perturbed grid, machine precision accuracy is again observed for setting 1 imposing full Dirichlet-type boundary conditions. With setting 2, the patch test is no longer satisfied as $B_\delta(\xb)\backslash\Omega$ is generally asymmetric with respect to the background grid. We plot the $L^2(\omg)$ errors of $\ub$ and $\theta$ vs. $h$ in Figure \ref{fig:patch_linear_un}. To investigate the impact of error in calculation of boundary normals, we present $L^2(\omg)$ errors either approximately using the estimate from \eqref{eqn:approx_n}, or using the exact normals. From Figure \ref{fig:patch_linear_un}, we observe an $O(h)$ convergence for the $L^2(\omg)$ error of $\ub$, and a deteriorated convergence rate for $\theta$. When comparing the numerical results from approximated $\mathbf{n}$ and exact $\mathbf{n}$, we surprisingly observed smaller numerical errors from the cases with approximate normal direction $\mathbf{n}$. %This is possibly due to the fact that with non-uniform discretizations, the approximated boundary $\partial\Omega_N$ via breaking bonds is generally not a straight line, which makes the analytical normal direction $\mathbf{n}$ less likely to pass the symmetry requirement in Thm. \ref{thm:linearpatch}. On the other hand, since \eqref{eqn:approx_n} calculates $\mathbf{n}$ from the current non-uniform discretization, it numerically provides an approximated $\mathbf{n}$ for each $\xb\in\omgi_N$ with better approximated symmetry for $B_\delta(\xb)\backslash\Omega$. 

In Figure \ref{fig:patch_corner_un} we consider boundary condition setting 3. Since there is no analytical normal direction defined on the corner point, we only investigate the results from approximated normal unit vector through formulation \eqref{eqn:approx_n}. Comparing with the $O(h)$ convergence rate in the uniform discretization cases, setting 3 converges with suboptimal $0.75$-th order convergence for $\ub$ and $0.5$-th order for $\theta$ on non-uniform grids.

\subsection{Manufactured solution test}\label{sec:mconv}

%\subsection{Straight Line Interface: Mixed BC}

To study the rate of convergence to the AC limit, we manufacture the local solution
$$\ub_0(x,y) = [\sin(Ax)\sin(Ay),-\cos(Ax)\cos(Ay)]$$
by imposing forcing consistent with the local operator
$$\fb(x,y)=[2(\lambda+2\mu)A^2\sin(Ax)\sin(Ay),-2(\lambda+2\mu)A^2\cos(Ax)\cos(Ay)].$$
A square domain
$\Omega=[-\pi/2,\pi/2] \times [-\pi/2,\pi/2]$ and the three boundary condition settings described in the previous Section \ref{sec:patch} are applied, to again consider the effect of boundary regularity.  
On $\omgbb_D$, the Dirichlet boundary condition $\ub_D(x,y)=\ub_0(x,y)$ is applied, while on the Neumann boundary we apply the traction condition
\begin{displaymath}
\mathbf{T}(x,y)=\left[\begin{array}{cc}
2A(\lambda+\mu)\cos(A\pi/2)\sin(Ay)&2A\mu\sin(A\pi/2)\cos(Ay)\\
2A\mu\sin(Ax)\cos(A\pi/2)&2A(\lambda+\mu)\cos(Ax)\sin(A\pi/2)\\
\end{array}\right]\mathbf{n}.
\end{displaymath}
We adopt material parameters under plane strain assumptions:
$$K=1,\lambda=K\nu/((1+\nu)(1-2\nu)),\mu=K/(2(1+\nu)),$$
and compare Poisson ratios $\nu=0.3 \text{ or }0.49$ again corresponding to compressible/near-incompressible limits. The parameter $A$ is taken as $0.4$. For the possible ambiguity of projection point in setting 3, we set $\mathbf{T}$ following a similar convention as in the linear patch test: for $\xb\in\omgi_N$ close to the corner point $(\pi/2,\pi/2)$, we set
$\mathbf{T}(\overline{\xb})\approx\mathbf{T}({\xb})$ where $\nb$ is numerically approximated with \eqref{eqn:approx_n}.

\textbf{Uniform discretization}: For Dirichlet boundary condition second-order convergence is achieved, consistent with the analysis in Thm. \ref{thm:Diri} and the $L^2(\omg)$ convergence results are presented in Figure \ref{fig:sin_Diri_s}. For traction loadings on straight and corner boundaries (Settings 2 and 3), we present $L^2(\omg)$ convergence results in Figure \ref{fig:sin_linear_s} and Figure \ref{fig:sin_corner_s}, respectively. In these settings first-order convergence is observed for both $\ub$ and $\theta$.

\textbf{Non-uniform discretization}: With non-uniform particle distribution, Figure \ref{fig:sin_Diri_un} demonstrates second-order $L^2(\omg)$ convergence for both $\ub$ and $\theta$ under Setting 1. Under Setting 2, Figure \ref{fig:sin_linear_un} demonstrates again first-order convergence. Again, somewhat surprisingly, when using the estimated normals one obtains improved accuracy, albeit with the same convergence rates. In Figure \ref{fig:sin_corner_un}, we further consider Setting 3 where $\partial\Omega_N$ includes a corner. Comparing with the results from uniform discretizations as shown in Figure \ref{fig:sin_corner_s}, a similar convergence rate ($O(h)$) is obtained for both $\ub$ and $\theta$ on non-uniform discretizations with setting 3.

\subsection{Traction loading on curvilinear surfaces}\label{ssec:isoSol}

 \begin{figure}[t!]
   \centering
   \includegraphics[width=0.9\textwidth]{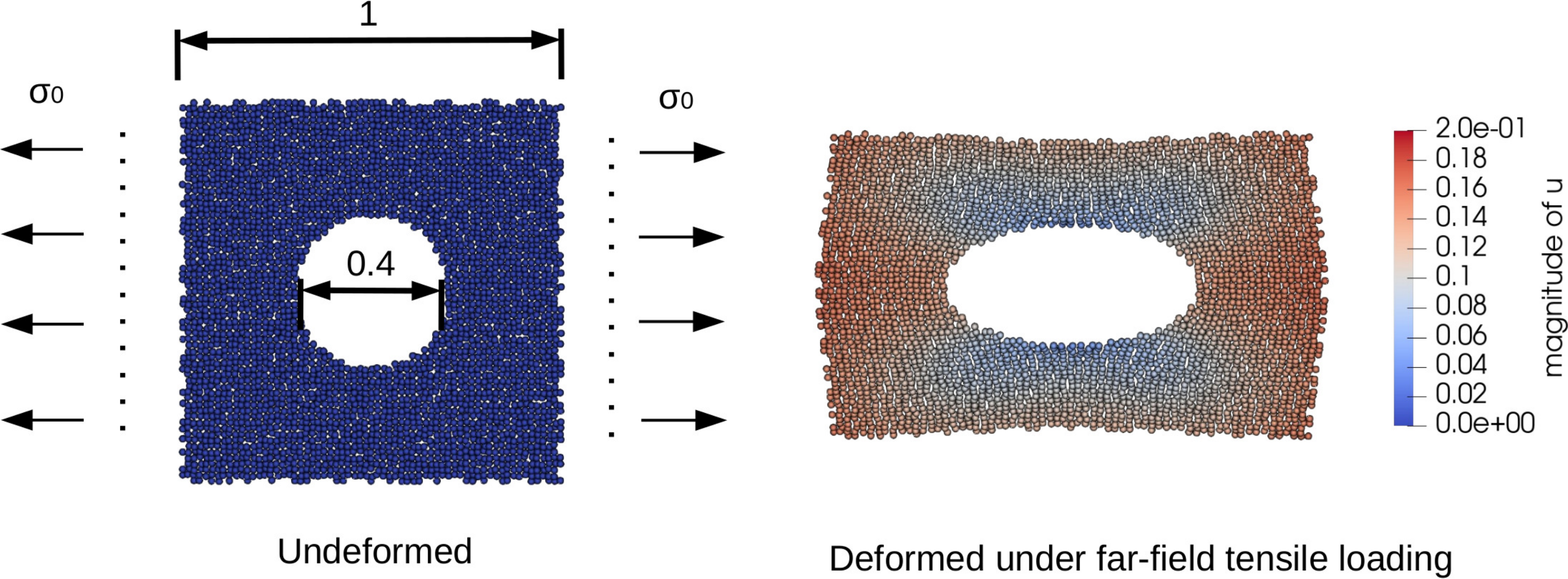}
   \caption{{Left:} Problem settings for circular hole in an infinite solid under remote loading with spheres representing a non-uniformed discretization. {Right:} Final deformed object when taking the far-field tensile stress $\sigma_0=0.3$.}
   \label{fig:holeInPlate}
 \end{figure}

 \begin{figure}[t!]
   \centering
   \includegraphics[width=0.8\textwidth]{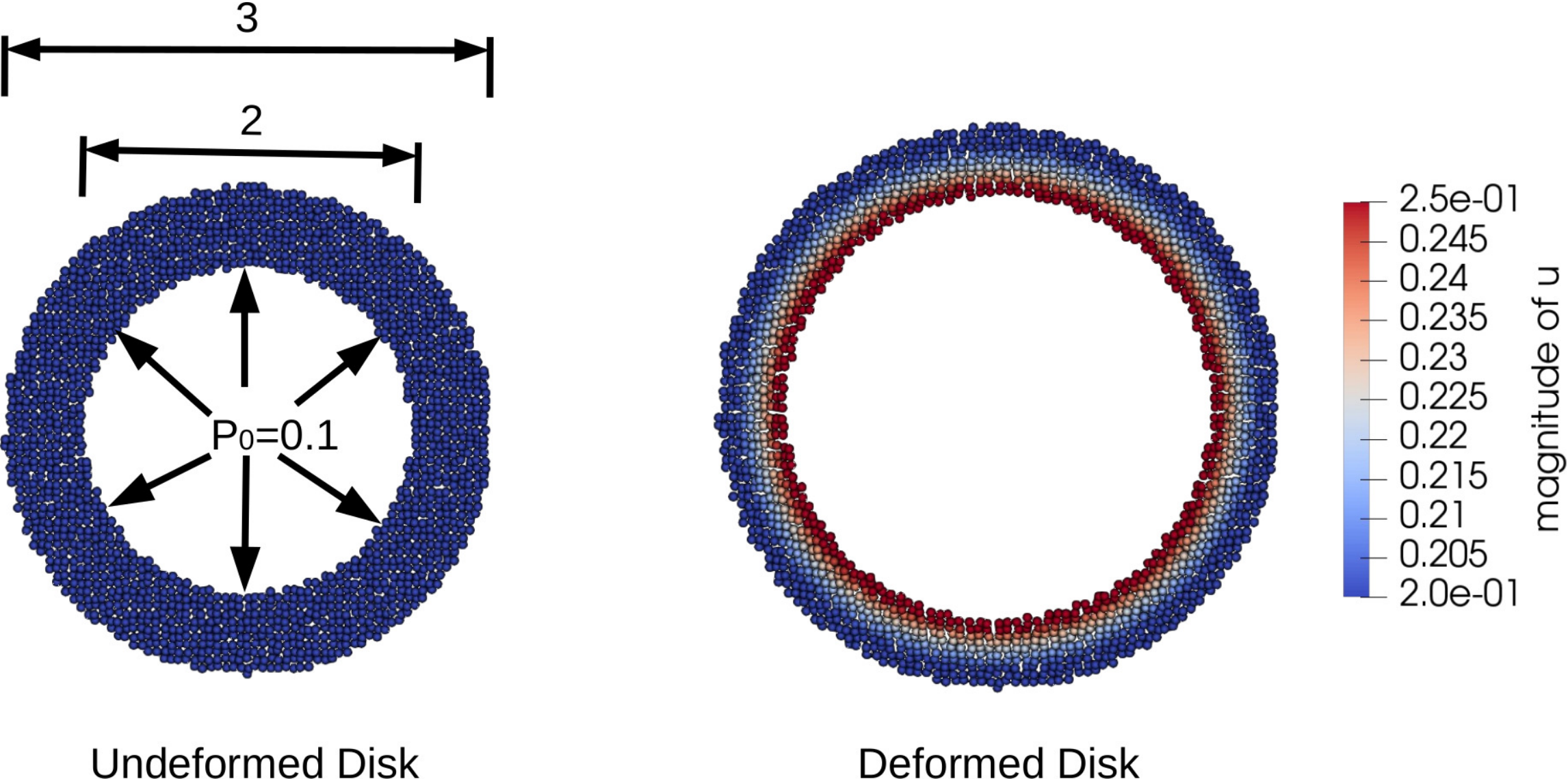}
   \caption{{Left:} Problem settings for a hollow disk under internal pressure with spheres representing a non-uniformed discretization. {Right:} Final deformed disk when taking the internal pressure $p_0=0.1$.}
   \label{fig:disk}
 \end{figure}

We next consider more physical settings corresponding to homogeneous and inhomogeneous traction loadings on a curvilinear surface. Two different problems are considered:
\begin{enumerate}
    \item We consider a free-surface circular hole of radius $a$ in an infinite solid under remote loading $\sigma_0$, as illustrated in Figure \ref{fig:holeInPlate}. Under a plane strain assumption the classical linear elasticity model yields the displacement field 
\begin{equation}\label{eq:holeSol}
\mathbf{u}_0(r,\theta) =\left( \begin{array}{c}
\frac{\sigma_0 a}{8 \mu} \left[ \frac{r}{a} (\kappa + 1) \cos \theta + \frac{2 a}{r}\left((1+\kappa) \cos \theta+\cos 3\theta \right) - \frac{2 a^3}{r^3} \cos 3 \theta \right]\\
\frac{\sigma_0 a}{8 \mu} \left[ \frac{r}{a} (\kappa - 3) \sin \theta + \frac{2 a}{r}\left((1-\kappa) \sin \theta+\sin 3\theta \right) - \frac{2 a^3}{r^3} \sin 3 \theta \right]\\
\end{array}\right).
\end{equation}
where $\kappa = 3-4\nu$ and $(r,\theta)$ are the radial distance and azimuthal angle in cylindrical coordinates. To set the problem up, we impose the analytic local solution $\ub_0$ as a Dirichlet-type condition on the nonlocal collar around the perimeter of a unit square. We then break all bonds crossing the circle of radius $a = 0.2$,  and apply $\mathbf{T}=0$ on the sharp interface $\partial\Omega_N=\{(x,y)|x^2+y^2=a^2\}$, fixing $\sigma_0=1$.
    \item As an example of imposing non-zero traction loads, we consider the deformation of a hollow cylinder under an internal pressure $p_0$, as illustrated in Figure \ref{fig:disk}. Under a plane strain assumption the classical linear elasticity model predicts displacements given by
    $$\ub_0(x,y) = \left[Ax+\dfrac{Bx}{x^2+y^2},Ay+\dfrac{By}{x^2+y^2}\right]$$
    where
    $$A=\dfrac{(1+\nu)(1-2\nu)p_0R_0^2}{K(R_1^2-R_0^2)},\;B=\dfrac{(1+\nu)p_0R_0^2R_1^2}{K(R_1^2-R_0^2)},$$
    $R_0$ and $R_1$ are the interior and exterior radius of the hollow disk. Here we take $R_0=1$ and $R_1=1.5$. We impose the analytic local solution $\ub_0$ as a Dirichlet-type condition on the nonlocal collar around the exterior boundary $\partial\Omega_D=\{(x,y)|x^2+y^2=R_1^2\}$ and break all bonds crossing the inner circle of radius $R_0 = 1$, and apply $\mathbf{T}=p_0\mathbf{n}$ on the interface $\partial\Omega_N=\{(x,y)|x^2+y^2=R_0^2\}$. In all M-convergence tests we take $p_0=0.1$.
\end{enumerate}
In the following tests we take the Young's modulus $K=1$ and test with different Poisson ratios $\nu=0.3 \text{ and }0.49$. To investigate the asymptotic compatibility when $\delta/h=C$, we employ uniform and non-uniform discretizations and refine $\delta$ and $h$ simultaneously while keeping the ratio $\delta/h$ a constant. In uniform discretizations, we take collocation points $X_h=\{(ah,bh),a,b\in\mathbb{N}\}$, and in non-uniform discretizations the uniform grid points are perturbed with $(\Delta x,\Delta y)$, $\Delta x, \Delta y \sim \mathcal{U}[-0.2h,0.2h]$. Note here even with uniform discretizations, the collocation points don't align with $\partial\Omega_N$ since the later is a circular curve, which introduces numerical errors as discussed in Thm. \ref{thm:geo}. In both settings we also investigate performances of the proposed formulation with the approximated normal unit vector formulation \eqref{eqn:approx_n} and  the analytical normal direction.

\textbf{Setting 1, Free circular surface}: In Figure \ref{fig:hole_s} and Figure \ref{fig:hole_un}, we demonstrate AC convergence for uniform and non-uniform particle distributions, respectively, for both compressible ($\nu=0.3$) and incompressible ($\nu=0.49$) materials. From the results, we observe first order convergence in the $L^2(\omg)$ norm for the displacements in all cases. A deteriorated order of convergence is observed for $\theta$, where a roughly $0.7$-th order is achieved for both uniform and non-uniform discretizations. When comparing the results with approximated $\mathbf{n}$ and the results with analytical $\mathbf{n}$, the formulation with analytical $\mathbf{n}$ performs similar or sometime slightly better than the results with approximated $\mathbf{n}$. Therefore, the formulation \eqref{eqn:approx_n} still provides a reasonable numerical approximation for the normal unit vector $\mathbf{n}$.

\textbf{Setting 2, Hollow disk under pressure}: In Figure \ref{fig:disk_s} and Figure \ref{fig:disk_un}, we present AC convergence for uniform and non-uniform particle configurations. We observe nearly first-order $L^2(\omg)$-norm convergence for displacements in both compressible and nearly-incompressible materials. Surprisingly, for nearly-incompressible materials, $O(h)$ order convergence is achieved for $\theta$, but for compressible material a reduced order (around $0.65$) convergence is observed. Again, the convergence rates are nearly identical for analytical $\nb$ and approximated $\nb$. Therefore, the approximation formulation \eqref{eqn:approx_n} again provides a good practical estimate of $\nb$.

\subsection{Composite materials with discontinuous material properties}\label{ssec:anisoSol}

\begin{figure}[t!]   \label{fig:inclusion}
   \centering
   \includegraphics[width=0.49\textwidth]{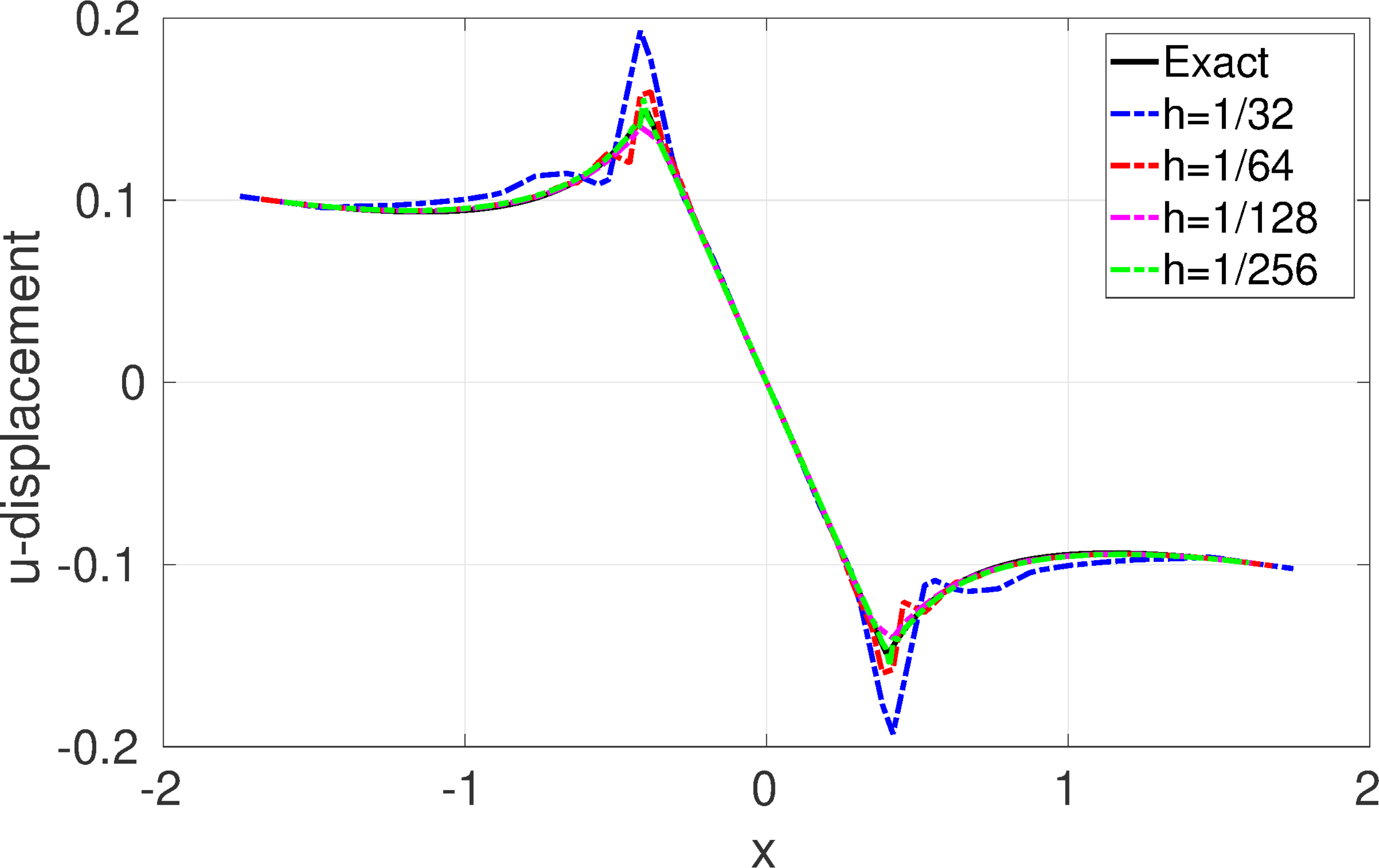}\;
   \includegraphics[width=0.49\textwidth]{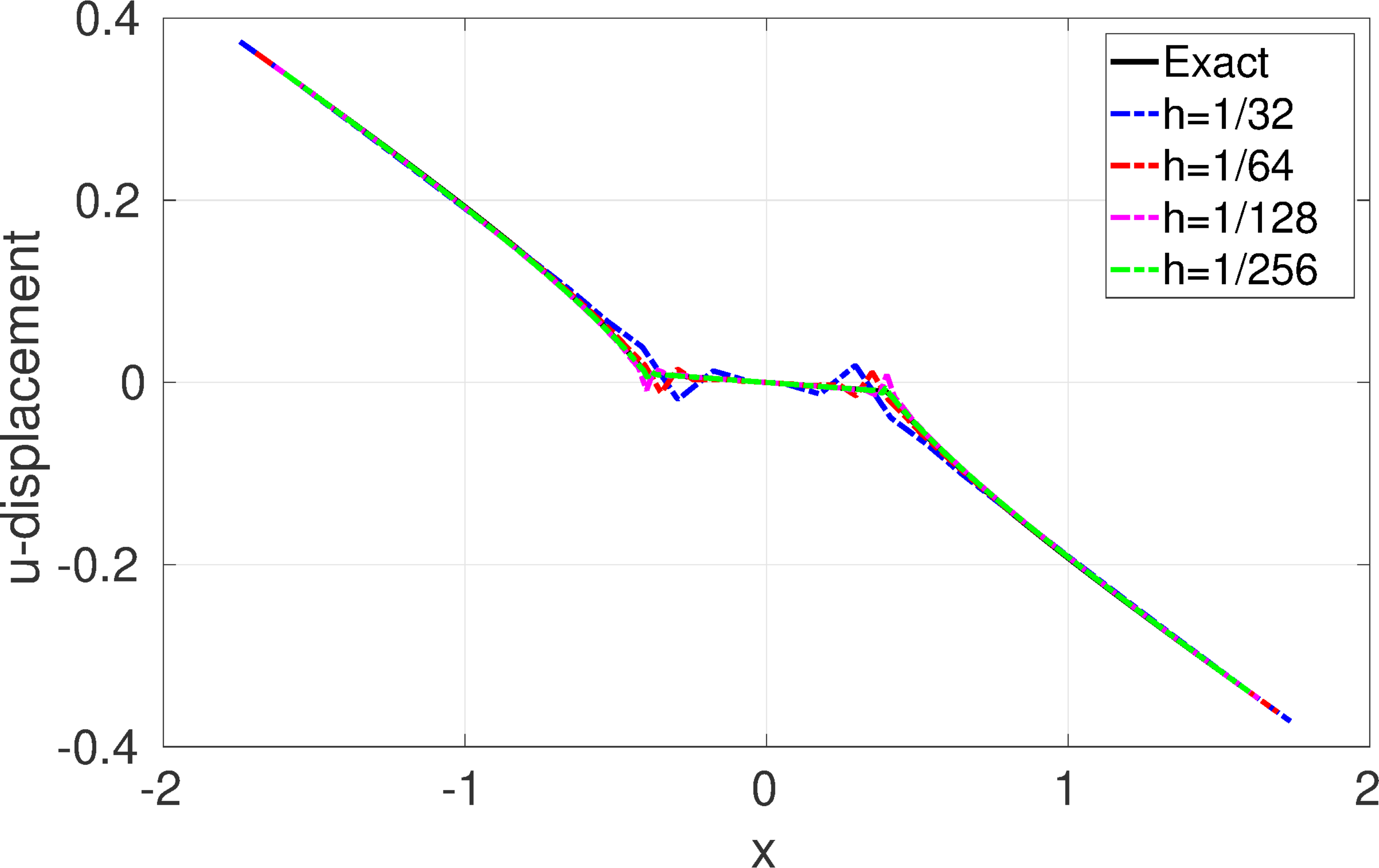}\;
   \includegraphics[width=0.49\textwidth]{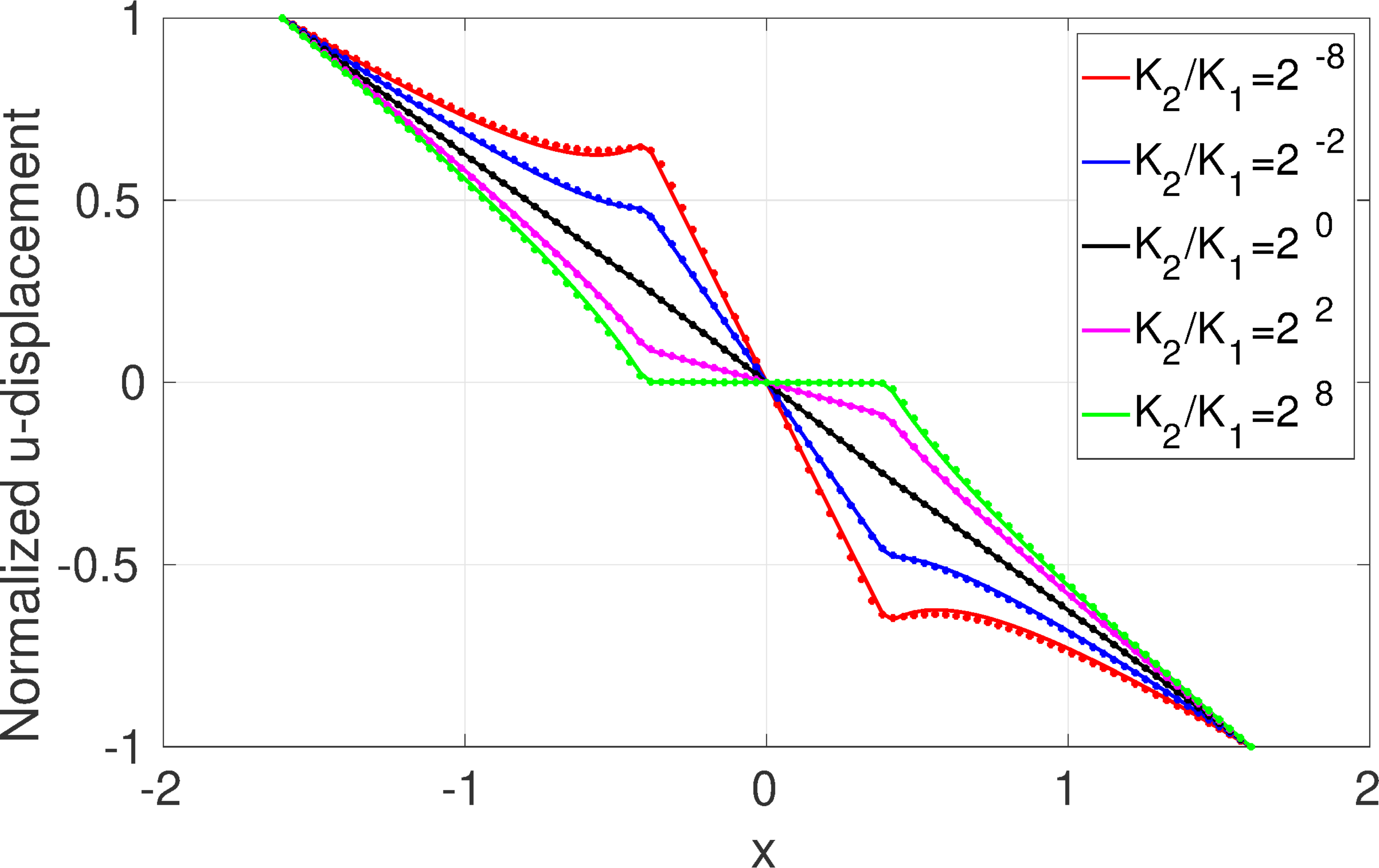}
   \caption{Displacements across the $y=0$ line, comparing prediction to analytic solution for composite problem. {Upper left:} Convergence to analytic solution for a soft inclusion ($K_2/K_1 = 1/64$). {Upper right:} Convergence to analytic solution for a stiff inclusion ($K_2/K_1 = 64$). {Bottom:} For a fixed resolution of $64^2$ points, reproduction of analytic solution for a wide range of $K_2/K_1 \in \left(2^{-8},2^8\right)$. Solid line corresponds to analytic solution, while dots correspond to numerical result.}
 \end{figure}

We now further consider an extension of the state-based peridynamics formulation \eqref{eq:nonlocElasticity} to composite materials constituted of $n$ phases, so that the domain may be partitioned into disjoint subdomains with piecewise constant material properties, i.e. $\Omega = \underset{k}{\cup}{\Omega_k}$, $\Omega_k\cap \Omega_l = \emptyset$, and $\lambda(\mathbf{x}) = \lambda_k$, $\mu(\mathbf{x}) = \mu_k$ for $\mathbf{x} \in \Omega_k$. Discussions on the mathematical properties of this heterogeneous system can be found in, e.g., \cite{capodaglio2020energy}. Specifically, when $\lambda(\xb)$ and $\mu(\xb)$ may vary for each material point $\xb$, we propose the following formulation:
\begin{equation}\label{eq:nonlocElasticity_comp}
    \mcL_\delta \mathbf{u}:=-\frac{C_\alpha}{m(\delta)}  \int_{B_\delta (\mathbf{x})} \left(\lambda(\xb,\yb) - \mu(\xb,\yb)\right) K(\left|\mathbf{y}-\mathbf{x}\right|) \left(\mathbf{y}-\mathbf{x} \right)\left(\theta(\mathbf{x}) + \theta(\mathbf{y}) \right) d\mathbf{y}
    \end{equation}
    \begin{equation*}
  -  \frac{C_\beta}{m(\delta)}\int_{B_\delta (\mathbf{x})} \mu(\xb,\yb) K(\left|\mathbf{y}-\mathbf{x}\right|)\frac{\left(\mathbf{y}-\mathbf{x}\right)\otimes\left(\mathbf{y}-\mathbf{x}\right)}{\left|\mathbf{y}-\mathbf{x}\right|^2}  \left(\mathbf{u}(\mathbf{y}) - \mathbf{u}(\mathbf{x}) \right) d\mathbf{y} = \mathbf{f}(\xb),
\end{equation*}
where the two-point functions $\mu(\cdot,\cdot)$, $\lambda(\cdot,\cdot)$ denote averaged material properties satisfying $\mu(\zb,\zb) = \mu(\zb)$ and $\lambda(\zb,\zb) = \lambda(\zb)$. We will consider for the purposes of this work the harmonic mean 
\begin{equation}
   \frac{2}{ \mu(\xb,\yb)} =\frac{1}{\mu(\xb)} + \frac{1}{\mu(\yb)},\quad\frac{2}{ \lambda(\xb,\yb)} =\frac{1}{\lambda(\xb)} + \frac{1}{\lambda(\yb)}.
\end{equation}
Correspondingly, to evaluate the above formulation, we modify the meshfree formulation with optimization-based quadrature weights in \eqref{eq:discreteNonlocElasticity} as follows:
\begin{align}
    \nonumber&-\frac{C_\alpha}{m(\delta)} \sum_{\xb_j \in B_\delta (\mathbf{x}_i)} \left(\lambda_{ij} - \mu_{ij}\right) K_{ij} \left(\mathbf{x}_j-\mathbf{x}_i\right)\left(\theta_i + \theta_j \right) \omega_{j,i}\\
  &-  \frac{C_\beta}{m(\delta)}\sum_{\xb_j \in B_\delta (\mathbf{x}_i)} \mu_{ij} K_{ij}\frac{\left(\mathbf{x}_j-\mathbf{x}_i\right)\otimes\left(\mathbf{x}_j-\mathbf{x}_i\right)}{\left|\mathbf{x}_j-\mathbf{x}_i\right|^2} \cdot \left(\mathbf{u}_j- \mathbf{u}_i \right) \omega_{j,i} = \mathbf{f}_i,\label{eq:discreteNonlocElasticity_comp}
\end{align}
where $\lambda_{ij}:=\lambda(\xb_i,\xb_j)$ and $\mu_{ij}:=\mu(\xb_i,\xb_j)$.

We numerically investigate the AC convergence of the nonlocal solution for a hydrostatically loaded cylindrical inclusion of radius $a$ in an infinite plate. We denote the interior of the inclusion as $\Omega_1$ and the exterior as $\Omega_2$, with corresponding constant material properties $\left(\mu_1,\lambda_1\right)$ and $\left(\mu_2,\lambda_2\right)$. Assuming a far-field hydrostatic stress $P_\infty$ and plane strain conditions, we define the coefficients
\begin{align*}
C_A &=   \frac{P_\infty}{2(\lambda_1 + \mu_1)},\\
C_B &=   \frac{P_\infty\left(\lambda_1 + \mu_1 +\mu_2\right)}{2(\lambda_1 + \mu_1) (\lambda_2 + 2 \mu_2)},\\
C_C &= - \frac{P_\infty a^2\left(\lambda_1-\lambda_2 + \mu_1 -\mu_2\right)}{2(\lambda_1 + \mu_1) (\lambda_2 + 2 \mu_2)},
\end{align*}
and the analytic local solution for the displacement field in cylindrical coordinates is given by 
\begin{align*}
u_r &=
\begin{cases}
C_A r         \quad \mathbf{x} \in \Omega_1,\\
C_B r + C_C/r \quad \mathbf{x} \in \Omega_2,
\end{cases}\qquad u_\theta = 0.
\end{align*}
We use this solution to assess the stability of the method in the vicinity of large jumps in material properties - for such scenarios high-order meshfree reconstructions have been shown to demonstrate unphysical oscillations near material interfaces \cite{trask2017high}. Note that the consistency conditions derived only guarantee asymptotic compatibility under the assumption of an isotropic material; this benchmark thus explores the applicability of the approach beyond the guarantees of the approximation theory in Thm. \ref{thm:Diri}.

We first investigate whether the discretization is AC. We take $a = 0.2$, impose a jump in the bulk modulus ($K_1 = 2$, $K_2 = 1$), and apply the analytic local solution $\ub_0$ as a Dirichlet-type condition on the nonlocal collar around the perimeter of a unit square. We consider three scenarios corresponding to different material compressibilities: (1) $\nu_1 = \nu_2 = 0.25$, (2) $\nu_1 =\nu_2= 0.49$ and (3) $\nu_1=0.49, \nu_2 = 0.25$. We present $L^2(\omg)$ convergence for both uniform and randomly perturbed particle distributions in Figure \ref{fig:hetero_s} and Figure \ref{fig:hetero_un} respectively. For both scenarios, for all three Poisson ratio combinations we obtain first- and half- order convergence for the displacement and dilitation, respectively.

We next investigate the stability of the approach over a large range of material parameters. To do this, we set $\mu_1=1$ and fix the Poisson ratio in both phases to $\nu_1 = \nu_2 = \frac14$ and impose a jump in the Young's modulus of $K_2/K_1 = Q$, for $Q \in \left\{2^{-8}, 2^8\right\}$. In these tests we employ a square domain $\Omega:=[-\pi/2,\pi/2]\times[-\pi/2,\pi/2]$ and apply the analytic local solution $\ub_0$ as a Dirichlet-type condition on the nonlocal collar around the perimeter of $\omg$. In Figure \ref{fig:inclusion}, we plot a profile of the x-component of displacement along the $y = 0$ line to provide a qualitative assessment of the solution. We demonstrate convergence for both a stiff inclusion ($K_1 = 64 K_2$), a soft inclusion ($K_2 = 64 K_1$), and then illustrate that we reproduce well the displacement for a wide range of parameters.

 \section{Fracture dynamics for brittle fracture experiments}\label{sec:exp}

The previous sections have established the ability of the scheme to recover local solutions of boundary value problems in elasticity governed by traction loadings and ensured that the breaking bonds treatment does not impair the AC convergence of the quadrature treatment. Of course, the main appeal of peridynamic discretizations is to handle fracture problems, and we devote the remainder of the paper to demonstrating how the scheme prescribed previously adapts to practical engineering settings, where now free surfaces are associated with the time evolution of a fracture surface. We specifically consider brittle fracture mechanics in linearly elastic materials and provide validation against experiment and existing numerical results. The main objective of this section is to provide a proof-of-principle demonstration that the framework introduced thus far applies to realistic settings, however overall the provided preliminary validation provides good agreement.

%The previous sections have endeavored to verify that the state-based peridynamics formulations described in this paper can accurately approximate solutions to the nonlocal model and achieve asymptotic compatibility to the correct local limits in linear elasticity equations \eqref{eqn:local}. However, this model is not guaranteed a priori to describe the dynamics of material crack initiation and propagation. The present section serves both to further illustrate the application of state-based peridynamics and the meshfree discretization to realistic problems and to argue that the modeling assumptions from mixed boundary conditions in \eqref{eq:newform1} can represent the dynamics of boundary tensile loadings and new interfaces created by cracks. Specifically, we numerically simulate the cracks in brittle solids and provide preliminary experimental validation of the model by quantitatively comparing our numerical simulation outputs with the results of established numerical simulations and/or experiments.

In this section we introduce an inertial term to handle dynamics
{\small\begin{equation}\label{eqn:probdyn}
\left\{\begin{array}{ll}
    \rho \frac{\partial^2{\ub}(\xb,t)}{\partial t^2}+\mcL_\delta\ub(\xb,t) = \mathbf{f}(\xb,t),& \text{ for }(\xb,t)\in\Omega\backslash\omgi_N\times[0,T],\\
\rho \frac{\partial^2{\ub}(\xb,t)}{\partial t^2}+\mcL_{N\delta}\ub(\xb,t) = \mathbf{f}_{N\delta}(\xb,t),& \text{ for }(\xb,t)\in\omgi_N\times[0,T],\\
%%%%%%%%%%%%%%%%%%%%%%%%%%%%%%% 
\theta(\mathbf{x},t)=\dfrac{d}{m(\delta)}\int_{B_\delta (\mathbf{x})} K(\left|\mathbf{y}-\mathbf{x}\right|) (\mathbf{y}-\mathbf{x})^T \left(\mathbf{u}(\mathbf{y},t) - \mathbf{u}(\mathbf{x},t) \right)d\mathbf{y},& \text{ for }(\xb,t)\in(\omg\cup\omgb_D\backslash\omgi_N)\times[0,T],\\
\theta(\mathbf{x},t)=\dfrac{d}{m(\delta)}\int_{B_\delta (\mathbf{x})\cap \Omega} K(\left|\mathbf{y}-\mathbf{x}\right|) (\mathbf{y}-\mathbf{x})^T \mathbf{M}(\mathbf{x})\left(\mathbf{u}(\mathbf{y},t) - \mathbf{u}(\mathbf{x},t) \right)d\mathbf{y},& \text{ for }(\xb,t)\in\omgi_N\times[0,T],\\
\mathbf{u}(\mathbf{x},t)=\mathbf{u}_D(\mathbf{x},t), & \text{ for }(\xb,t)\in\omgbb_D\times[0,T],\\
\mathbf{u}(\mathbf{x},0)=\mathbf{u}_{IC}(\mathbf{x}),\;\frac{\partial\mathbf{u}(\mathbf{x},0)}{\partial t}=\mathbf{v}_{IC}(\mathbf{x}),\;\frac{\partial^2\mathbf{u}(\mathbf{x},0)}{\partial t^2}=\mathbf{w}_{IC}(\mathbf{x}), & \text{ for }\xb\in\omg\cup\omgbb_D,
\end{array}\right.
\end{equation}}
where $\rho$ is the material density and $\ub_{IC}$, $\vb_{IC}$, $\wb_{IC}$ are the initial displacement, velocity and acceleration fields, respectively. To model brittle fracture, for $\xb_j\in B_\delta(\xb_i)$ we break the bond between $\xb_i$ and $\xb_j$ when the associated strain exceeds a critical strain criteria $s_0$:
\begin{equation}\label{eq:strain}
    s_{ij}: = \frac{||\mathbf{u}_j-\mathbf{u}_i + \xb_j-\xb_i||-||\xb_j - \xb_i||}{||\xb_j - \xb_i||}>s_0.   
\end{equation}
We employ the criteria derived in \cite{zhang2018state} relating $s_0$ to material parameters:
\begin{equation}\label{eq:criteria}
    s_0 = \sqrt{\frac{G_0}{4(\lambda - \mu)\beta' + 8\mu \beta}},
\end{equation}
where $\beta = \frac{3\delta}{4\pi}$, $\beta' = 0.23873 \delta$, and $G_0$ is the critical energy release rate or fracture energy. For $\xb_j\in B_\delta(\xb_i)$, this damage criterion can be implemented by replacing the static state weight $\gamma_{j,i}$ in \eqref{eqn:gamma} with a history-dependent scalar boolean state function $\gamma_{j,i}(t)$:
\begin{align}
    \nonumber \gamma_{j,i}(t) &= \begin{cases}
    1, \quad \text{if }s_{ij}(\tau)\leq s_0,\;\forall\tau\leq t, \text{ and }\xb_j\in B_\delta(\xb_i)\Omega,\\
    0, \quad \text{otherwise}, \\
    \end{cases}
\end{align}
such that $\tilde{\omega}_{j,i}=\omega_{j,i}\gamma_{j,i}(t)$, $\hat{\omega}_{j,i}=\omega_{j,i}(1-\gamma_{j,i}(t))$. To postprocess fracture evolution and identify cracks, we define the damage as
\begin{equation}
    \phi_i(t) = \frac{\sum\limits_{\xb_j \in B(\xb_i)} (1-\gamma_{j,i}(t))}{\sum\limits_{\xb_j \in B(\xb_i)}1}.
\end{equation}
For the purposes of fracture identification, we say that a crack occurs at $\xb_i$ if $\phi_i$ exceeds $0.35$. This threshold is somewhat arbitrary, but necessary to, e.g., postprocess the crack propagation velocity at which a crack grows.

To discretize we apply the Newmark scheme together with the meshfree quadrature established previously
\begin{equation}\label{eqn:probdis}
\left\{\begin{array}{ll}
    \frac{4\rho}{\Delta t^2} \ddot{\mathbf{u}}_i^{n+1}+(\mcL^h_\delta
    \mathbf{U})_i^{n+1}= \mathbf{f}_{i}^{n+1}
    +\frac{4\rho}{\Delta t^2}(\mathbf{u}_{i}^n + \Delta t \dot{\mathbf{u}}_i^n + \frac{\Delta t^2}{4}
    \ddot{\mathbf{u}}_i^n)
    ,& \quad \text{for }\xb_i \text{ in }\Omega\backslash\omgi_N^{n+1},\\
\frac{4\rho}{\Delta t^2} \ddot{\mathbf{u}}_i^{n+1}+(\mcL^h_{N\delta}\mathbf{u})_i^{n+1} = (\mathbf{f}^h_{N\delta})_i^{n+1}+\frac{4\rho}{\Delta t^2}(\mathbf{u}_{i}^n + \Delta t \dot{\mathbf{u}}_i^n + \frac{\Delta t^2}{4}\ddot{\mathbf{u}}_i^n),&\quad\text{for } \xb_i \text{ in }\omgi_N^{n+1},\\
%%%%%%%%%%%%%%%%%%%%%%%%%%%%%%%
\theta_i^{n+1}=\dfrac{d}{m(\delta)} \sum\limits_{\xb_j \in B_{\delta}(\xb_i)} K_{ij} (\mathbf{x}_j-\mathbf{x}_i)^T \left(\mathbf{u}_j^{n+1} - \mathbf{u}_i^{n+1} \right)\omega_{j,i},&\quad \text{for }\xb_i \text{ in }\omg\cup\omgb_D\backslash\omgi_N^{n+1},\\
\theta_i^{n+1}=\dfrac{d}{m(\delta)} \sum\limits_{\xb_j \in B_{\delta}(\xb_i)} K_{ij} (\mathbf{x}_j-\mathbf{x}_i)^T \mathbf{M}_i\left(\mathbf{u}_j^{n+1} - \mathbf{u}_i^{n+1} \right)\tilde{\omega}_{j,i},&\quad \text{for }\xb_i \text{ in }\omgi_N^{n+1},\\
\mathbf{u}_i^{n+1}=\mathbf{u}_D(\mathbf{x}_i), & \quad \text{for }\xb_i \text{ in }\omgbb_D,\\
%\mathbf{u}_i^0=\mathbf{u}_{IC}(\mathbf{x}_i),\;\dot{\ub}_i^0=\mathbf{v}_{IC}(\mathbf{x}_i),\;\ddot{\mathbf{u}}_i^0=\mathbf{w}_{IC}(\mathbf{x}_i), &\quad \text{ for }\xb_i \in\omg\cup\omgbb_D,\\
\end{array}\right.
\end{equation}
where $\Delta t$ is the time step size, $\fb_i^{n+1}:=\fb(\xb_i,t^{n+1})$, $\mcL^h_{\delta}$ and $\mcL^h_{N\delta}$ are the discretized nonlocal operators as defined in \eqref{eq:discreteNonlocElasticity} and \eqref{eq:discreteNonlocElasticity2}, respectively, and $\fb_{N\delta}^h$ is also as defined in \eqref{eq:discreteNonlocElasticity2}. The acceleration and velocity at the $n+1$-th time step are then calculated as follows:
$$\ddot{\mathbf{u}}_i^{n+1} 
:= \frac{4}{\Delta t^2}(\mathbf{u}_i^{n+1}
-\mathbf{u}_i^{n}-\Delta t \mathbf{u}_i^n) - \ddot{\mathbf{u}}_i^n, \quad
\dot{\mathbf{u}}_i^{n+1} := \dot{\mathbf{u}}_i^n
+ \frac{\Delta t}{2} (\ddot{\mathbf{u}}_i^n
+ \ddot{\mathbf{u}}_i^{n+1}).$$
Note that because the evolving fracture creates new free surfaces, $\partial\omg_N$ and $\omgi_N$ alter with time. To capture the implicit coupling between the material response and the evolving geometry due to fracture evolution, we employ subiterations at each time step 
%Particularly, based on the displacement field, we evaluate the damage criteria \eqref{eq:strain} for each bond. If a bond is broken, we subcycle to repeat the current step until no new broken bonds occur. 
%At each simulation step, a subiteration technique is used. 
{as follows. We first assume no new bonds have been broken at the current time step and solve for the displacement field. Based on the displacement field, we evaluate the damage criteria \eqref{eq:strain} for each bond. If any bond meets the criteria of breaking, we break all these bonds, update the corresponding state functions $\gamma_{j,i}$ and quadrature weights $\hat{\omega}_{j,i}$ and $\tilde{\omega}_{j,i}$, then solve for the displacement field again with new free surfaces. We repeat this procedure until no new broken bonds are detected, and finally proceed to the next time step.}
% To summarize:
% \begin{enumerate}
%     \item Set initial conditions for $\mathbf{u}_i^0$, $\dot{\mathbf{u}}_i^0$, and $\ddot{\mathbf{u}}_i^0$.
%     \item for $n = 0,1,2,\dots, T/\Delta t$, do 
%     \begin{enumerate}[a.]
%         \item Solve for the displacement field via \eqref{eqn:probdis}.
%         \item Check the damage criteria \eqref{eq:strain} based on the predicted displacement field, then update the quadrature weights $\tilde{\omega}_{j,i}$ and $\hat{\omega}_{j,i}$ if $s_{ij}>s_0$. Go to step c if no new broken bond is found, otherwise repeat step a.
%         \item Update the displacement, velocity and acceleration fields at time step $n+1$, and continue to the $n+1$-th iteration.
%     \end{enumerate}
% \end{enumerate}

We consider three benchmark problems involving material damage. In Sections \ref{sec:glass1}-\ref{sec:glass2} we study dynamic crack propagation and branching in glass. In Section \ref{sec:glass1}, we adopt the benchmark problem from \cite{bobaru2015cracks} and simulate a pre-cracked glass plate under sudden tensile loading. In Section \ref{sec:glass2}, we reproduce a recent experiment considering V-notched glass samples impacted by a striker \cite{dondeti2020comparative}. In Section \ref{sec:ring}, we simulate the material fragmentation of a cylinder under internal pressure and identify the number of fragments. %Although the main aim of this section is to demonstrate the capability of the proposed formulation on realistic problem settings, we provide preliminary validation results by comparing our numerical results with available numerical simulations and experimental measurements. 

\subsection{Dynamic brittle fracture I: Pre-cracked glass under tensile loading}\label{sec:glass1}
 
  \begin{figure}[t!]
   \centering
   \includegraphics[width=0.4\textwidth]{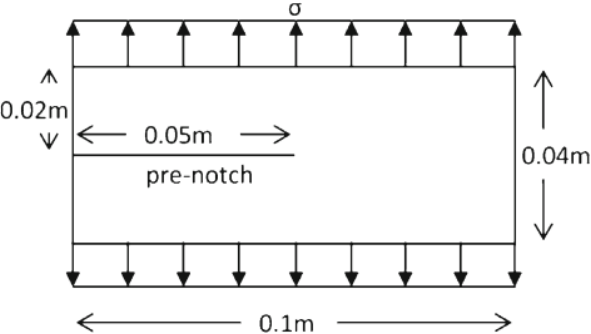}
   \caption{Problem setup for pre-cracked soda-lime glass experiment under tensile loading, following \cite{ha2010studies}.}
   \label{fig:glasssetup}
 \end{figure}
 
\begin{table}[]
\center
\begin{tabular}{|cccc|}
\hline
Young's modulus $E$ & Poisson ratio $\nu$ & Density $\rho$&Fracture energy $G_0$\\
\hline
72 $GPa$&0.23&2440 $kg/m^3$&3.8 $J/m^2$\\
\hline
\end{tabular}
\caption{Matieral properties used in pre-cracked soda-lime glass experiment.}
\label{tab:glass}
\end{table}

We first investigate the crack propagation and branching of soda-lime glass as a prototypical brittle fracture exemplar, whereby a pre-notched thin rectangular plate is subject to tensile loads on its top and bottom (Figure \ref{fig:glasssetup}). Following the setup in \cite{ha2010studies}, we consider plate dimensions of $0.1 m$ by $0.04 m$ with an initial crack of length $0.05m$, and a constant tensile load $\sigma=2MPa$ applied on the top and bottom of the sample starting at $T=0$. All other boundaries, including the new boundaries created by cracks, are treated as free surfaces. The mechanical properties of soda-lime glass are listed in Table \ref{tab:glass}. This problem was studied in several numerical studies on bond-based peridynamics \cite{ha2010studies,gu2017voronoi,du2017peridynamic} and non-ordinary state-based peridynamics \cite{zhou2016numerical} (see \cite{diehl2019review} for a review). Experimentally the crack propagation speed is fairly reproducible and was reported as 1580 $m/s$ in \cite{bowden1967controlled}. To validate our scheme's ability to reproduce crack propagation speed and branching location we compare against available numerical and experimental results from \cite{ha2010studies,bobaru2015cracks,bowden1967controlled}.

%, and the numerical results were quantitatively compared with experimental results published on the crack propagation velocity in the region of branching \cite{field1971brittle} or the maximum propagation velocity measured \cite{bowden1967controlled}.

  \begin{figure}[t!]
   \centering
   \includegraphics[width=0.95\textwidth]{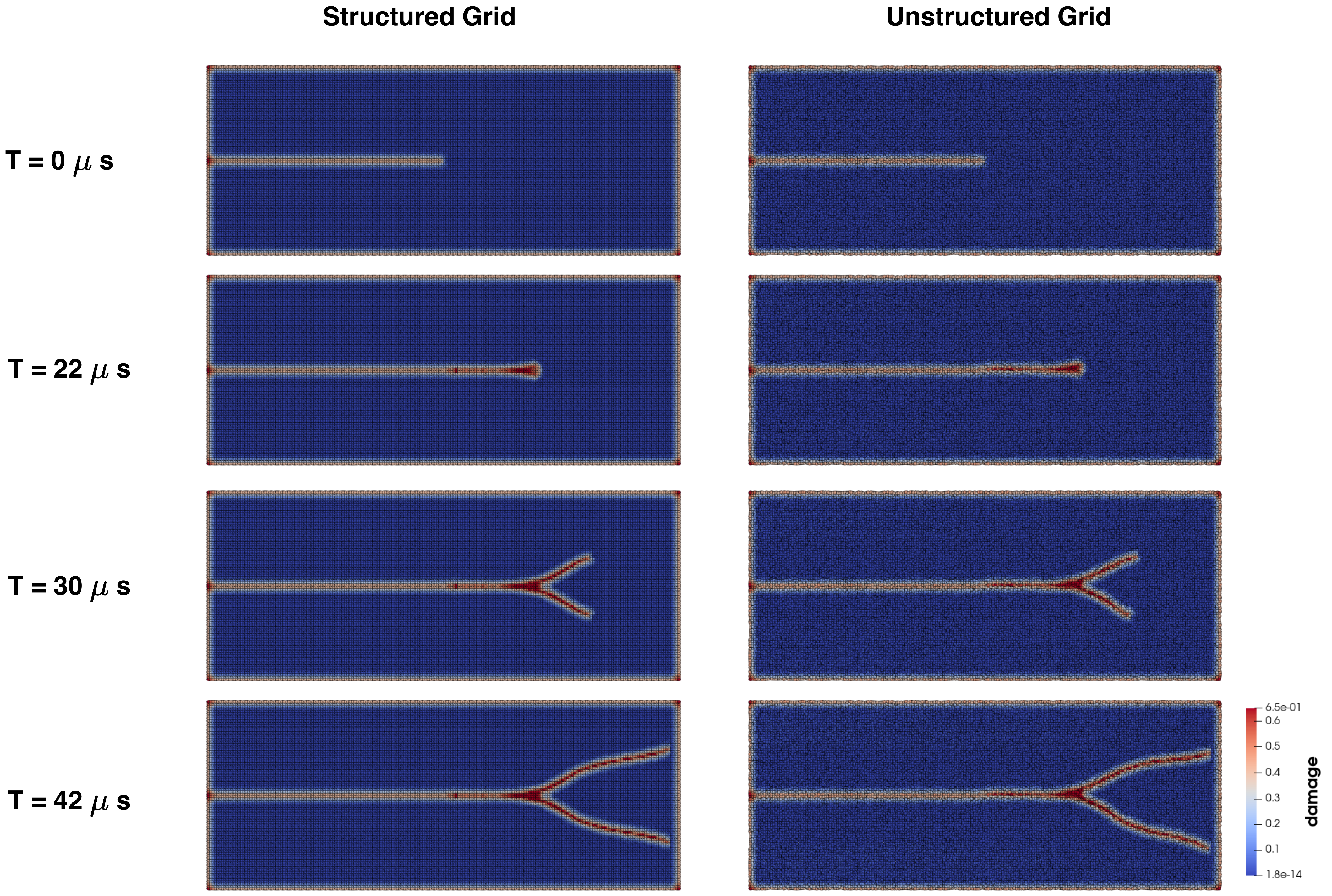}
   \caption{Simulated damage evolution of pre-cracked soda-lime glass crack branching study, using resolution $h=0.0005m$ and $\delta=0.002m$. \textit{Left:} uniform discretization, \textit{Right:} non-uniform discretization.}
   \label{fig:glassdamage}
 \end{figure}
 
    \begin{figure}[t!]
   \centering
   \includegraphics[width=0.95\textwidth]{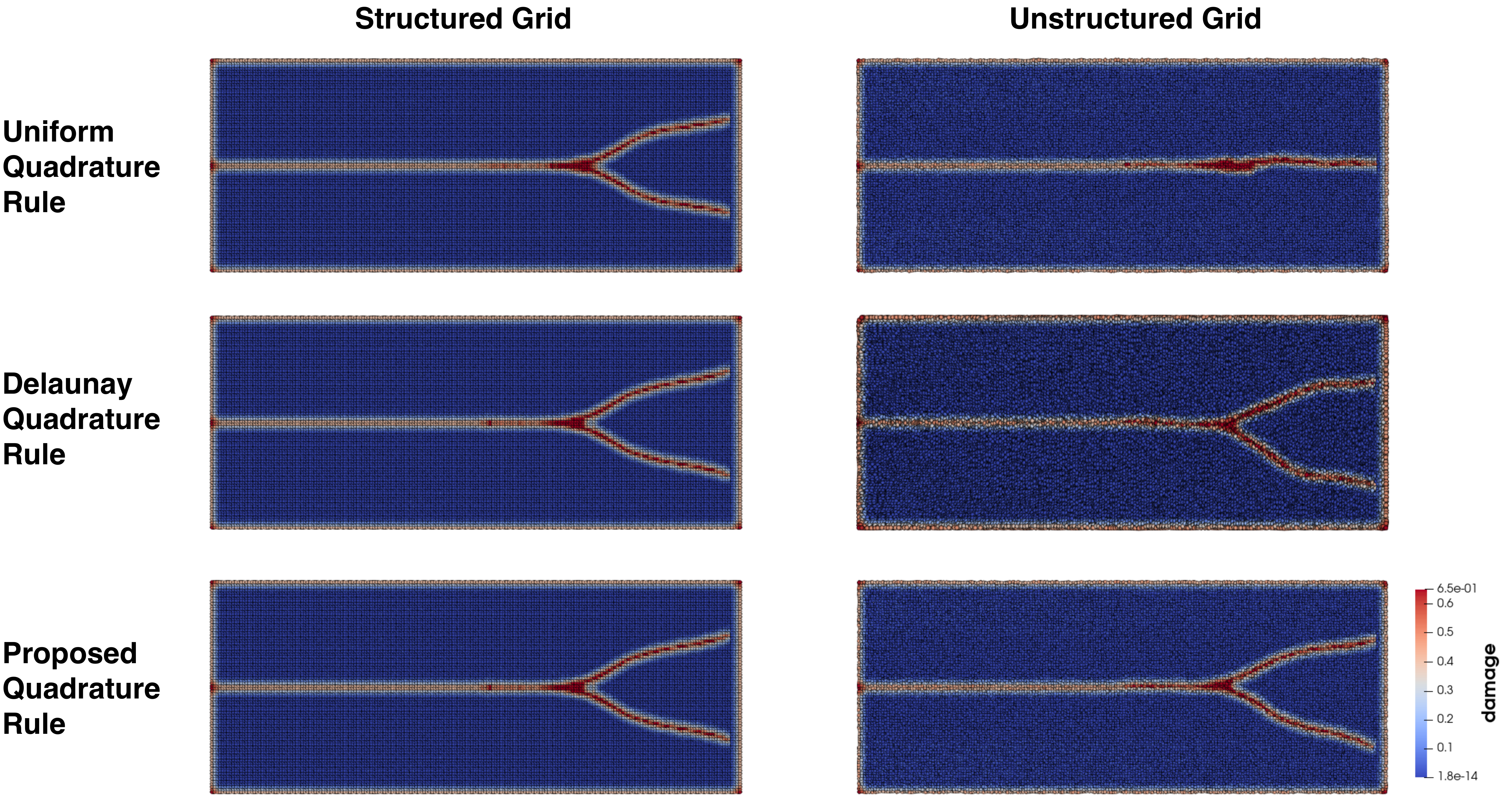}
   \caption{Simulated fracture pattern at $T=42\mu s$ for pre-cracked soda-lime glass crack branching study, using various quadrature strategies and resolution $h=0.0005m,\,\delta=0.002m$. \textit{Left:} uniform discretization, \textit{Right:} non-uniform discretization. For non-uniform grids, the proposed meshfree approach predicts branching in line with a Delaunay mesh-based discretization \cite{gu2017voronoi,bussler2017visualization}.}
   \label{fig:glassdamage_mtd}
 \end{figure}

We first plot in Figure \ref{fig:glassdamage} the fracture evolution based on a uniform grid spacing with $h=5\times 10^{-4}m$, $\delta=4h$ and $\Delta t=6.25 \times 10^{-2} \mu s$. A qualitative comparison to the results from \cite[Figure 5]{bobaru2015cracks} shows that we qualitatively recover the same dynamics as existing simulations, independent of particle distribution. 

For engineering applications, non-uniform discretizations are desirable to handle complex geometries and establish grid independence. For many discretizations, so-called grid-imprinting may qualitatively numerically skew fracture patterns so that they correlate with mesh orientation and special care is often required in numerical methods \cite{bobaru2015cracks,bobaru2011adaptive}. To this end, we compare the effect of particle anisotropy on the resulting fracture, comparing our approach to a popular meshfree quadrature rule from designed for Cartesian particle distributions \cite{parks2008pdlammps,Yu2018paper}. We also compare against a mesh-based approach, building a Delaunay mesh on a Cartesian grid with nodal spacing $h$ and assigning particles at cell centroids with quadrature weight equal to the cell measure \cite{gu2017voronoi,bussler2017visualization}. We generate non-uniform discretizations by perturbing either the particle locations or Delaunay nodes by $0.2h$, and consider $h=5\times 10^{-4}m$. Ideally, we would hope to recover results comparable to the mesh-based approach on non-uniform discretizations, without the need to introduce a mesh into the problem. By comparing the corresponding fracture patterns in Figure \ref{fig:glassdamage} and Figure \ref{fig:glassdamage_mtd}, we can see that with our proposed meshfree scheme is more robust to particle anisotropy than traditional meshfree quadrature, providing nearly identical results on uniform or nonuniform grids.

 \begin{table}[]
\center
\begin{tabular}{|c|ccc|cc|}
\hline
Quantity & Exp & BB ($\delta=$2e-3m) & BB ($\delta=$5e-4m) & SB ($\delta=$2e-3m) & SB ($\delta=$1e-3m)\\                       
\hline
Branching Location ($m$) & -- & 0.065 & 0.068 & 0.070 & 0.068\\
Branching Time ($\mu s$)& -- & 23.0 & 21.5 & 21.8 & 22.5\\
%Averaged Prop Speed ($m/s$) &1500 \cite{field1971brittle}&--&--& 1466 & 1444 \\
Max Prop Speed ($m/s$) & 1580 & 2000& 1679 & 2250 & 2000\\
\hline
\end{tabular}
\caption{Quantitative comparison of crack dynamics to existing experimental and numerical works. Here ``Exp'' stands for experimental results from \cite{bowden1967controlled}, ``BB'' is the estimated result with the bond-based peridynamics measured from \cite{ha2010studies}, and ``SB'' corresponds to the current approach with the state-based peridynamics.}
\label{tab:glasssamage}
\end{table}

\begin{figure}
\centering 
  \begin{subfigure}
  \centering
    \includegraphics[width=.49\linewidth]{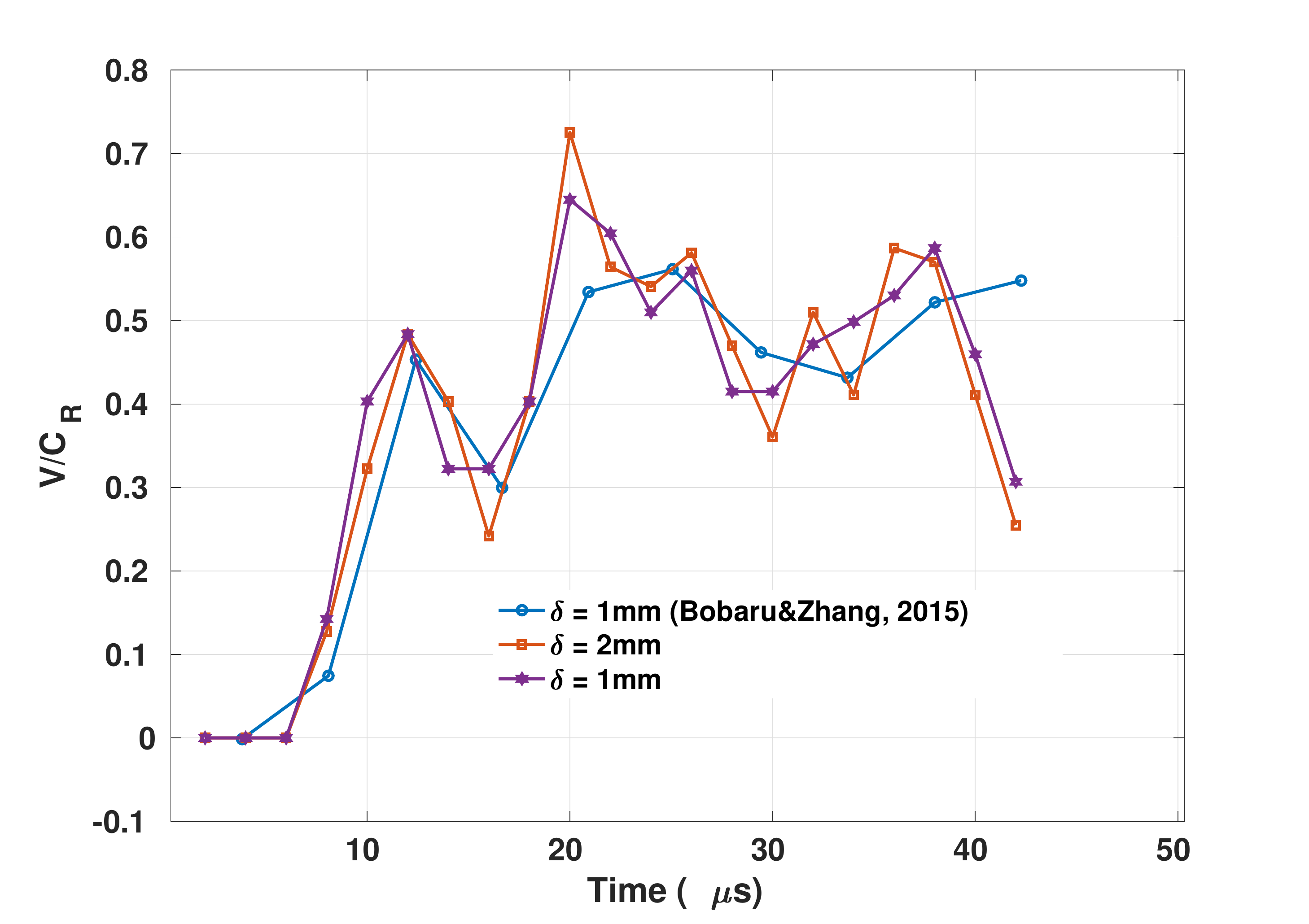}
  \end{subfigure}
%  \begin{subfigure}
%  \centering
%    \includegraphics[width=.49\linewidth]{pics/v_mid.eps}
%  \end{subfigure}
  \caption{Comparison of the soda-lime glass (normalized) crack propagation speed: proposed formulation versus numerical results reported in \cite{bobaru2015cracks}.}
  \label{fig:soda_velocity}
  \end{figure}

To quantitatively validate our simulation results, we validate the time and location of crack branching and the crack propagation speed and compare against \cite{ha2010studies,bobaru2015cracks,bowden1967controlled}. In all experiments, we kept a fixed time step size $\Delta t=6.25 \times 10^{-2} \mu s$ and a fixed ratio $\delta/h=4$. Theoretically, the nonlocal length scale in state-based peridynamics should be smaller than geometrical features to prevent unrealistic nonlocal interactions. Therefore, we also investigate the M-convergence test by decreasing $h$ and $\delta$ simultaneous to see if crack propagation features converge, considering $h = 5\times10^{-4}m$ and $h=2.5\times 10^{-4}m$. In Table \ref{tab:glasssamage}, we compare these quantities of interest against numerical and experimental data. We obtain good agreement for the branching time and location, but overestimate the maximum speed. This may be a result of under-resolution, as the overestimation is reduced under refinement. However, we note that several other methods \cite{diehl2019review,zhou2016numerical,gu2017voronoi} achieve similar results. To conclusively establish an improvement in the current formulation regarding this quantity of interest, we defer a deeper investigation of this discrepancy to an upcoming work involving a parallel implementation of the current scheme allowing a more involved refinement study.

In Figure \ref{fig:soda_velocity} we plot the predicted crack propagation speeds under different $\delta$ as functions of time, and compare them with the numerical results from \cite{ha2010studies}. All results are normalized by the Rayleigh wave speed $c_R = 3102 m/s$. We can observe that the numerical simulation show a similar trend: prior to the crack entering the branching phase, the speed gradually decreases, and then rapidly increases after branching. These trends are also observed in experiments \cite{sundaram2018dynamic,dondeti2020comparative}. 
 
 \subsection{Dynamic brittle fracture II: V-notched Glass Under Impact}\label{sec:glass2}
 
 \begin{figure}
 \centering
 \begin{subfigure}
\centering
\includegraphics[width=.38\linewidth]{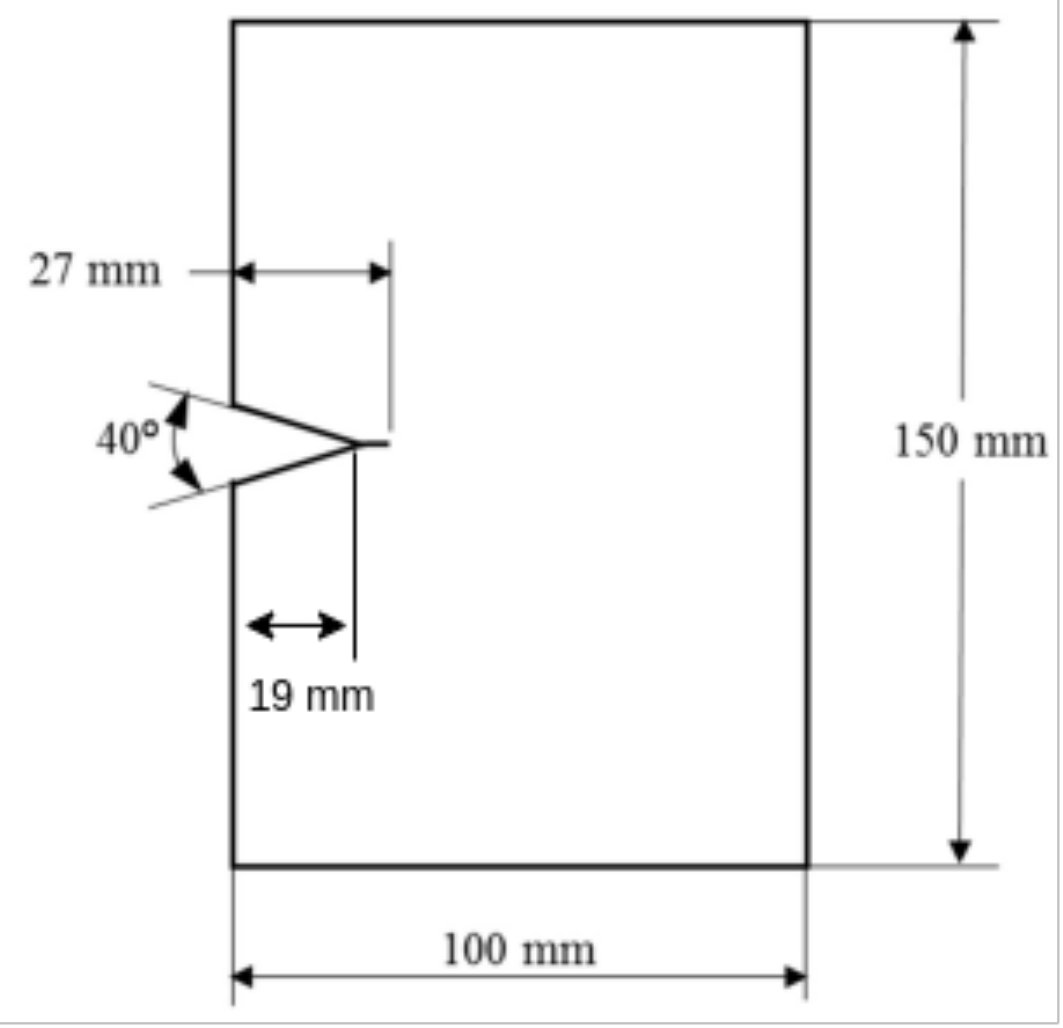}
\end{subfigure}
  \caption{Problem setup for V-notched soda-lime glass specimen under impact, following \cite{dondeti2020comparative}.}
  \label{fig:vnotch_setup}
\end{figure}

\begin{figure}
\centering
  \begin{subfigure}
  \centering
    \includegraphics[width=.22\linewidth]{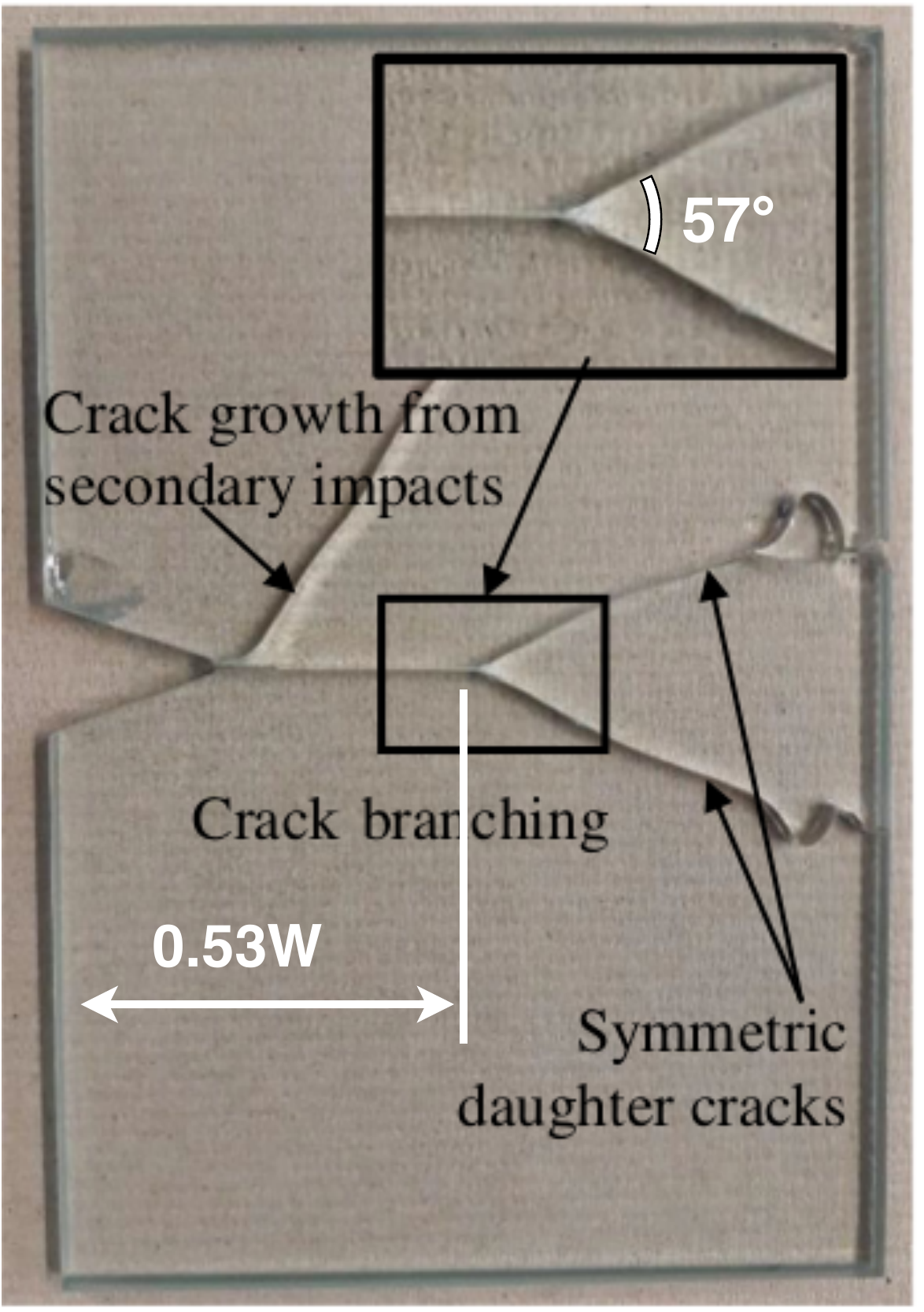}
  \end{subfigure}
  \begin{subfigure}
  \centering
    \includegraphics[width=.22\linewidth]{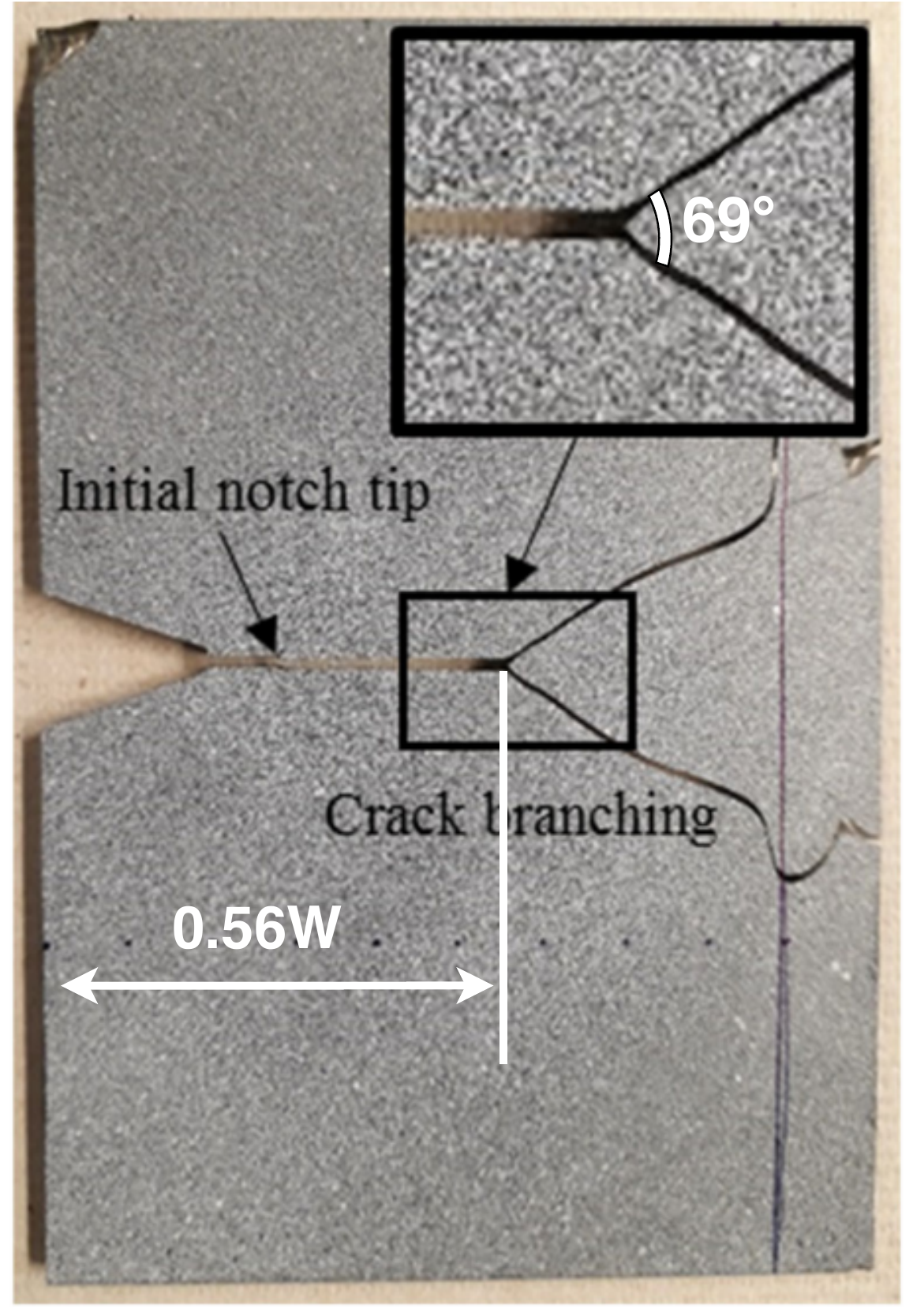}
  \end{subfigure}
  \begin{subfigure}
  \centering
    \includegraphics[width=.22\linewidth]{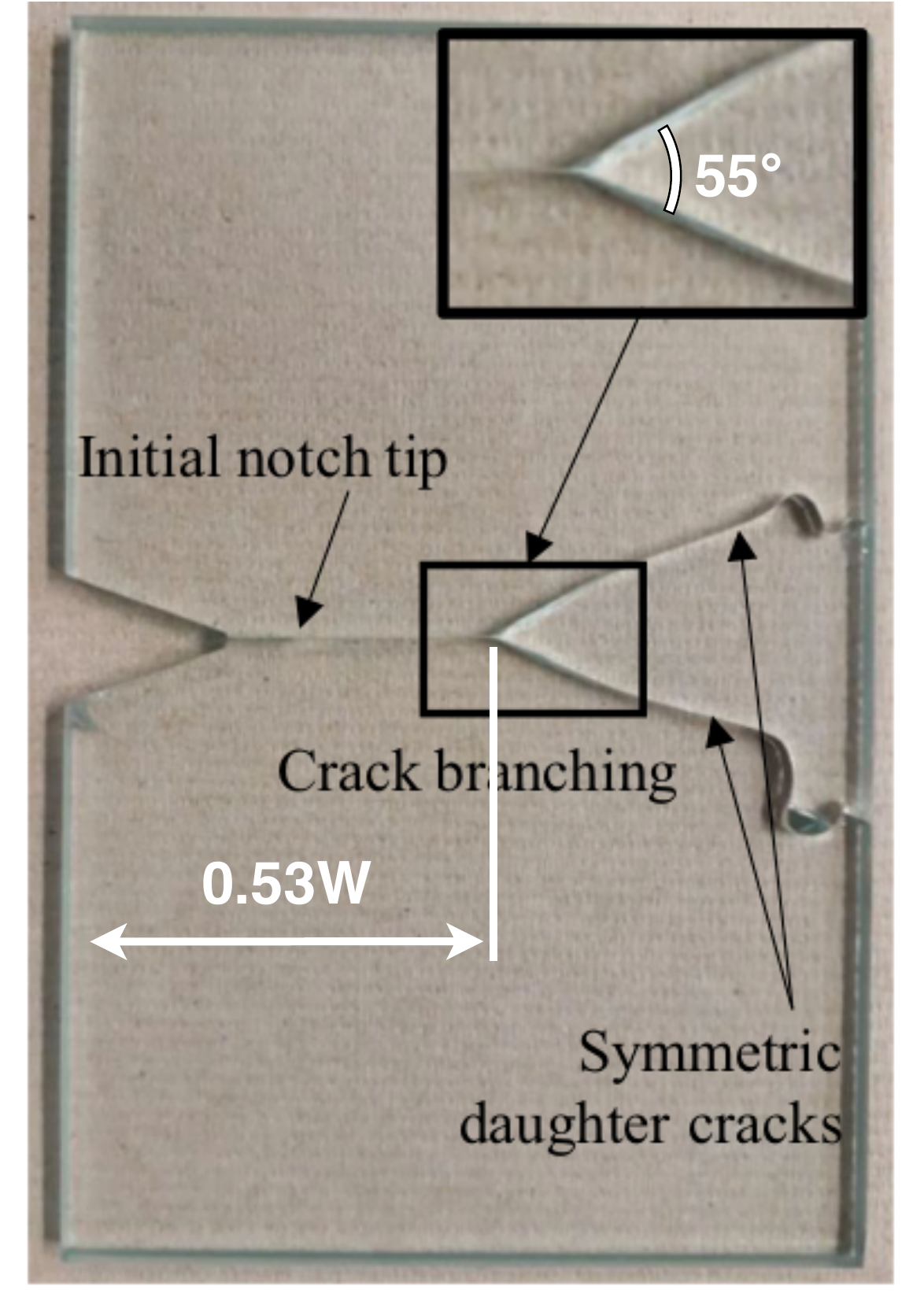}
  \end{subfigure}
  \begin{subfigure}
\centering 
    \includegraphics[width=.22\linewidth]{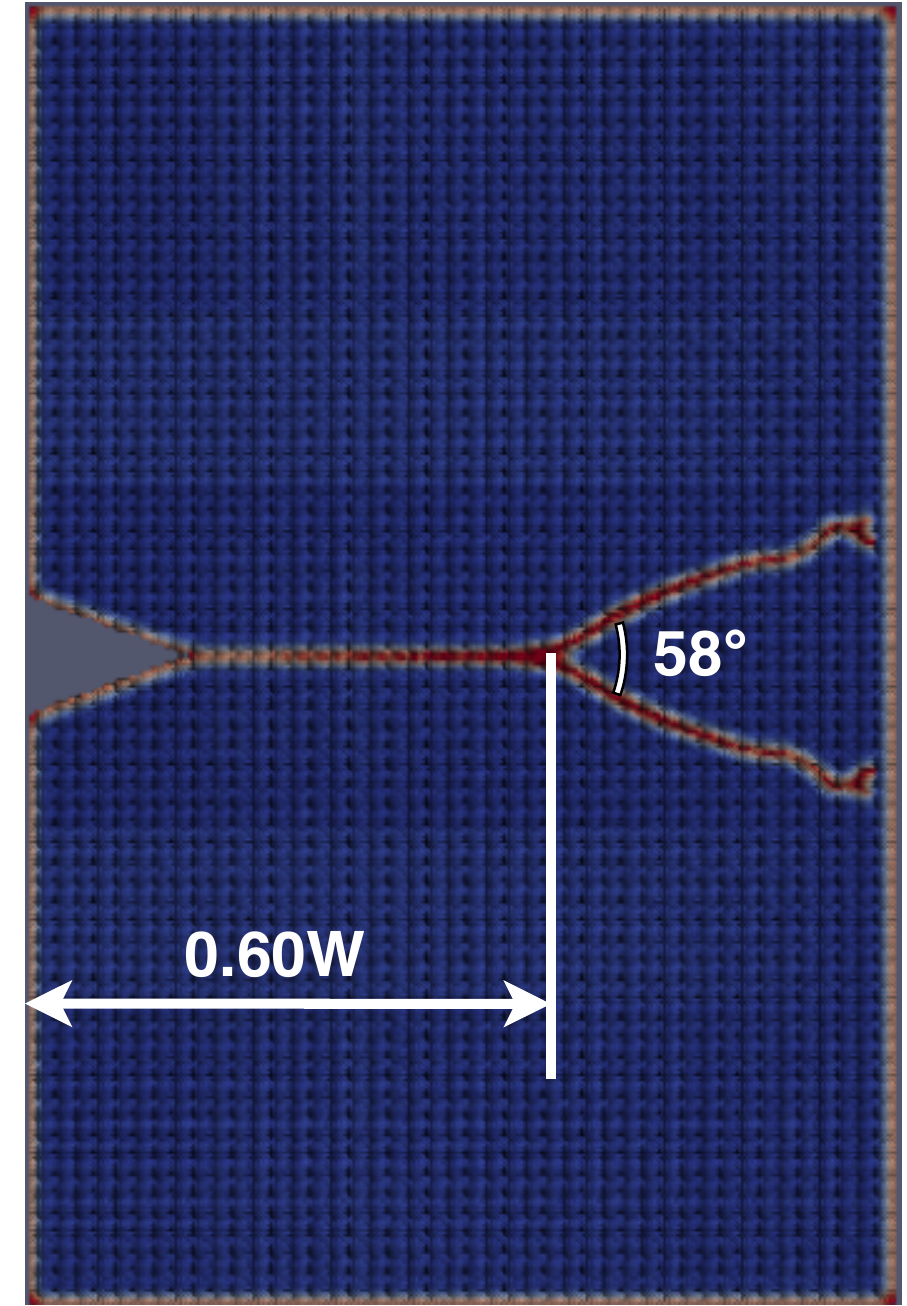}
  \end{subfigure}
  \caption{Experimental fracture patterns and comparison to numerical prediction for V-notch case. From left to right: experimental results from photoelasticity in \cite{dondeti2020comparative}; experimental results from Digital Image Correlation (DIC) in \cite{dondeti2020comparative}; experimental results from Digital Gradient Sensing (DGS) in \cite{dondeti2020comparative}; numerical simulation results from the proposed approach.}
  \label{fig:Vnotch_exp}
  \end{figure}

Recently, Dondeti and Tippur have studied impact-induced crack branching experiments on soda-lime glass by applying three prevalent optical techniques: transmission photoelasticity, 2D Digital Image Correlation (DIC) and transmission Digital Gradient Sensing (DGS) \cite{dondeti2020comparative}. Following the setup sketched in Figure \ref{fig:vnotch_setup}, a Hopkinson pressure bar was used to impart a impulse upon a V-notch and study the resulting fracture - we defer to \cite{dondeti2020comparative} for further details of the experimental setup. Three nominally identical but separate experiments were carried out to compare three different optical techniques in \cite{dondeti2020comparative}, and the experimental results are reproduced here in the first three plots of Figure \ref{fig:Vnotch_exp}. Although the branching location and the branching angles were not reported in \cite{dondeti2020comparative}, we used the photographs shown in Figure \ref{fig:Vnotch_exp} to measure branch locations and angles to serve as validation data. For the three specimens, crack branching was observed at $53\%$, $56\%$, and $53\%$ of the width, with branching angles as $57^\circ$, $69^\circ$, and $55^\circ$, respectively. Moreover, one can observe that in all specimens the crack path presents small oscillations near the far end of the sample due to wave reflections/spalling, which we aim to reproduce.
 \begin{table}[]
\center
\begin{tabular}{|cccc|}
\hline
Young's modulus $E$ & Poisson ratio $\nu$ & Density $\rho$& Fracture energy $G_0$\\
\hline
70 $GPa$&0.22&2500 $kg/m^3$& 8 $J/m^2$\\
\hline
\end{tabular}
\caption{Mechanical properties for V-notched soda-lime glass specimens taken from \cite{dondeti2020comparative}.}
\label{tab:glass_Vnotch}
\end{table} 

 \begin{figure}
\centering
\begin{subfigure}
\centering
\includegraphics[width=.49\linewidth]{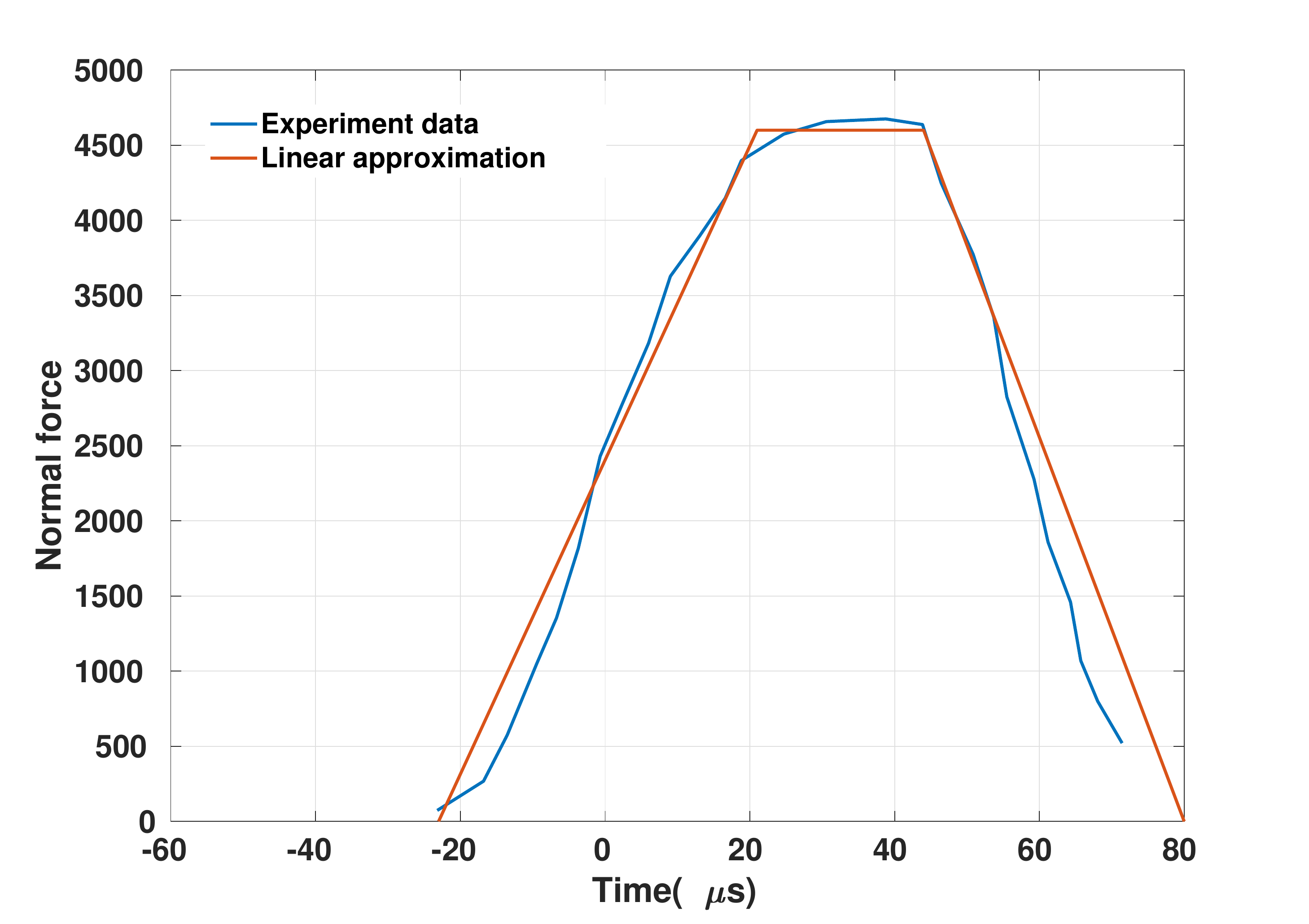}
\end{subfigure}
\begin{subfigure}
\centering
\includegraphics[width=.49\linewidth]{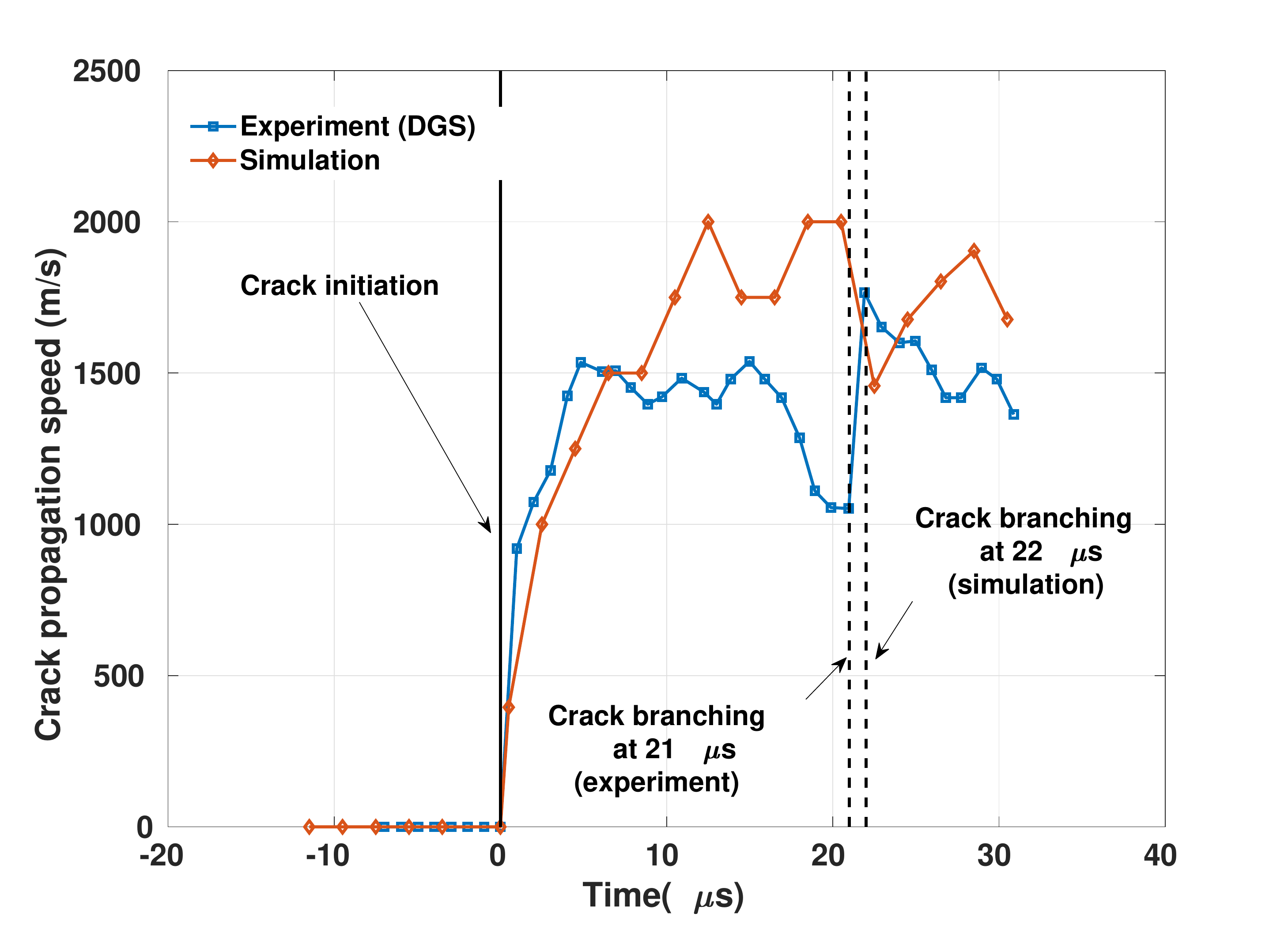}
\end{subfigure}
\caption{Left: Normal force loads applied to V-notch, demonstrating forces measured from experiment (DGS) (reproduced from \cite{dondeti2020comparative}) and approximate piecewise linear load applied in simulations. Right: Crack speed profiles, comparing experimental results (DGS) (reproduced from \cite{dondeti2020comparative}) to numerical prediction. 
%For numerical results, improved agreement is achieved if pressure is applied over a portion of the V-notch (13.1mm/19mm)\cite{mehrmashhadi2020comparison}. 
Time is denoted after the onset of branching.}
\label{fig:vnotch_forcespeed}
\end{figure} 

In this experiment, detailed information regarding contact force history, crack propagation speed, branching angle and point of branching are provided, allowing validation of numerical simulations against experiment using identical experimental loading conditions. In Table \ref{tab:glass_Vnotch} we list the material properties of soda-lime glass as provided in \cite{dondeti2020comparative}, where the fracture energy $G_0 = 8J/m^2$ is measured during the experiment using the DGS technique when crack initiates. The force histories on the V-notch faces of the specimen by the long-bar were evaluated with DGS, as reproduced in blue in the left plot of Figure \ref{fig:vnotch_forcespeed}. In our numerical simulations, a piecewise linear approximation of the applied normal force is applied uniformly over the V-notch surface as a time-varying traction load. Following the settings in \cite{dondeti2020comparative}, the frictional effect is neglected. Moreover, since the actual measurement of the bar tip shape was not provided in experiments, we assume that the full length of the V-notch is loaded, although we note that the predicted failure patterns might differ from the ones produced by the partial loading of the notch surfaces \cite{mehrmashhadi2020comparison}. Crack velocities were also estimated in \cite{dondeti2020comparative}, and the results indicate that both the photoelastic recording and the DGS method provided reliable velocity history profiles. %In \ref{fig:vnotch_forcespeed} we reproduce the experimental velocity history profile from DGS and compare our numerical estimates with it.

%recent experimental tests on dynamic fracture/crack branching in glass induced by impact. 
%We adopt the estimated normal force on the V-notch from the DGS approach, which can be found in the left figure of Figure \ref{fig:vnotch_forcespeed}. The frictional effect is neglected here. The normal force follows the 'ramp-up, plateau, ramp-down' pattern, more specifically, it hits the plateau right  after the crack branches. 
 %In the simulation, the crack initiation and crack
 %branching points are earlier than those in the experiment, %therefore we 
 %modify the normal force to keep it as a constant after branching before starting to decrease. A similar measure has been taken to approximate the loading conditions in \cite{mehrmashhadi2019uncovering}. 

To simulate the experiment, plane stress assumptions are adopted and traction loads are applied consistent with the experimentally measured normal force at the V-notch and free surfaces over the remainder of the boundary. A uniform discretization is employed with grid size $h=0.5mm$, and we select horizon $\delta = 4 h$, and time step $\Delta t = 0.125 \mu s$. The predicted fracture pattern and crack velocity profile is given in Figures \ref{fig:Vnotch_exp} and Figure \ref{fig:vnotch_forcespeed}, respectively. In Figure \ref{fig:Vnotch_exp}, results show that branching happens at location $60\%$ away from the left edge of the sample, with a branching angle of around $58^\circ$. While the branching angle matches very well within the range of angles from experimental measurements ($55^\circ$-$69^\circ$), the branching location is a little further than the measurements in experiment ($53\%-56\%$), in what follows we explore possible explanations. Oscillations in the fracture surface are reproduced near the back of the specimen. In Figure \ref{fig:vnotch_forcespeed} we provide comparisons of the crack speed as a function of time.

For the results provided, the results provide qualitative agreement sufficient for the purposes of this work. We do offer speculation regarding possible explanation and areas which may lead to improved quantitative agreement. Regarding the discrepancy in branching location, Mehrmashhadi et al. was able to achieve better agreement with experiment by applying the normal force loading over a subset of the full V-notch, to model the effect of reduced area under contact \cite{mehrmashhadi2020comparison}. We remark that we were able to achieve improved agreement in crack branching location with similar techniques. We omit any results along these lines however, as our focus is only to demonstrate our boundary treatment for a realistic problem and a careful analysis of physical modeling assumptions is beyond the scope of this work. We also note that Mehrmashhadi et al. was able to employ a finer mesh; again we defer a careful analysis of such effects to a future work where we introduce a scalable implementation.

 \begin{figure}
 \centering
 \begin{subfigure}
\centering
\includegraphics[width=.45\linewidth]{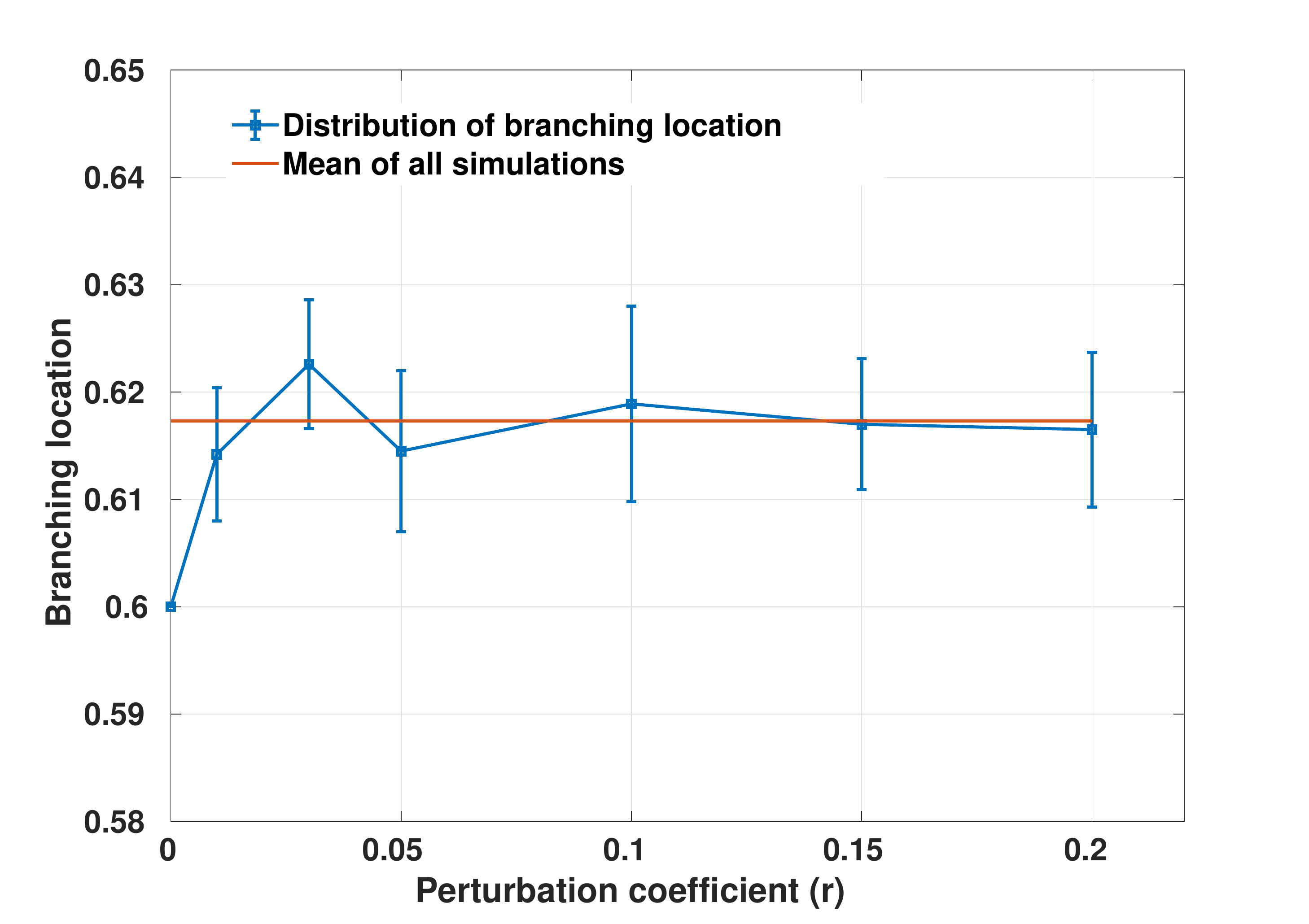}
\end{subfigure}
 \begin{subfigure}
\centering 
    \includegraphics[width=.45\linewidth]{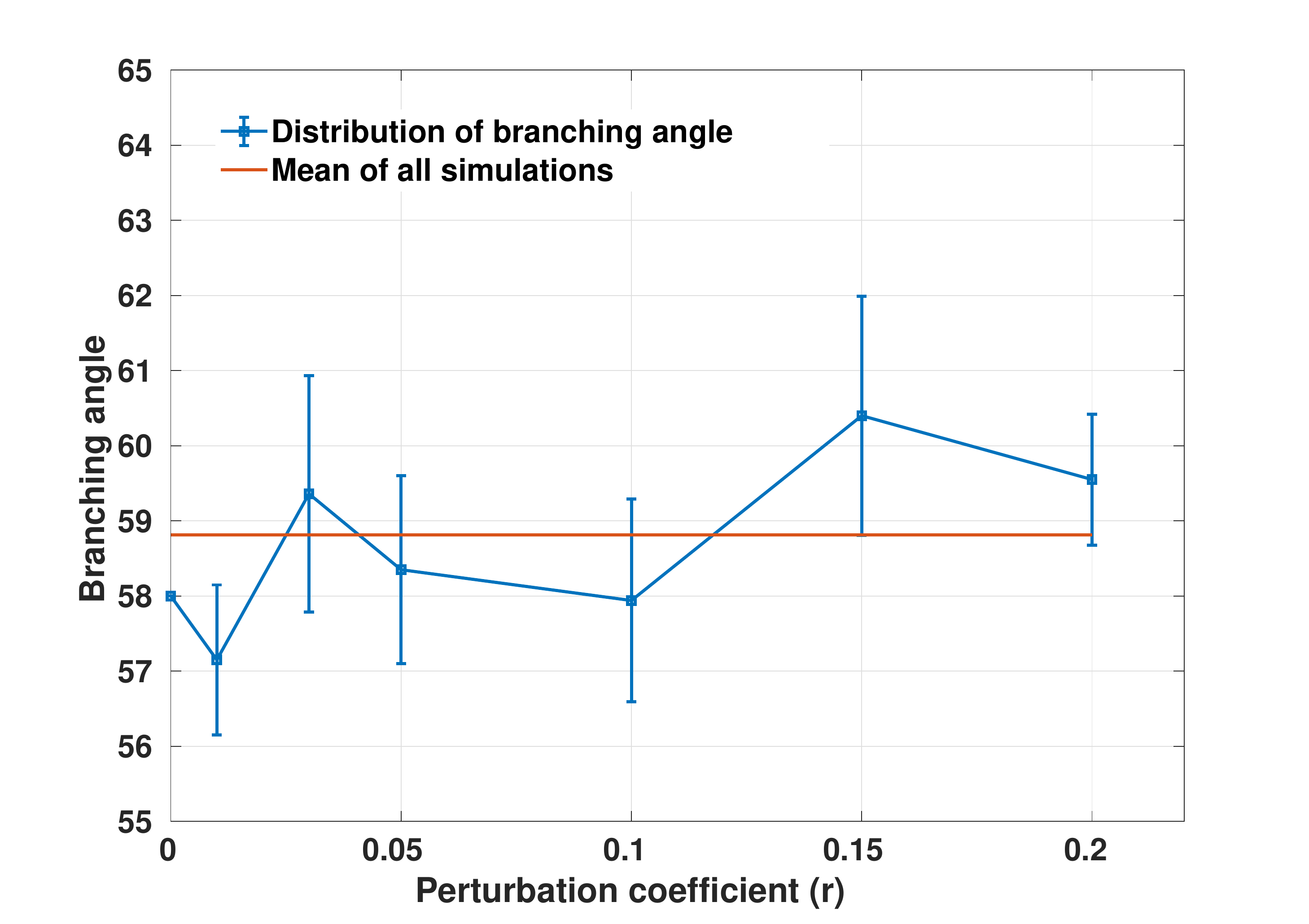}
  \end{subfigure}
  \caption{Reproducability and effect of anisotropy in spatial location of crack branching point. Non-uniform discretizations are employed by perturbing the uniform particle locations by $rh$, $r\in[0,0.2]$. Left: Mean and standard error of predicted branching location. Right: Mean and standard error of branching angles.}
  \label{fig:vnotch_unstructure}
\end{figure}

Next we characterize the reproducability of the predicted crack paths, considering in particular the effect of anisotropy in the underlying discretization. For an increasing magnitude of perturbation ratio $r\in(0,1)$, a quasi-uniform pointset is generated by perturbing every point in the uniform grid by a uniformly distributed random variable of magnitude $rh$. In this study we take $h = 0.5 mm$, $\Delta t = 0.125 \mu s$, $\delta=2mm$ and $r \in 
\{0.01, 0.03, 0.05,0.1,0.15,0.2\}$. For each $r$, we calculate solutions corresponding to 20 non-uniform particle distributions. To investigate the impact of non-uniform grids on crack features, we record the branching location and branching angle, and report their means and standard errors versus the grid perturbation ratio $r$ in Figure \ref{fig:vnotch_unstructure}. 

For the branching location, all simulations predict fairly consistent results: the crack starts to branch at around $60\%-62\%$ of the specimen width. Larger variations are observed on the branching angle when $r\geq 0.03$, which is possibly due to the fact that these estimates are sensitive to the placement of the branching points. Across all $r$, a mean angle with $57^\circ-60^\circ$ is predicted, which lies in the range observed from experiments ($55^\circ-69^\circ$). The numerical results indicate that these crack features are not overly sensitive to small perturbations in the discretization grids, demonstrating the suitability of the scheme to handle nontrivial problems without imparting grid anisotropy effects on the resulting fracture prediction.

\subsection{Fragmentation of Cylinder Expansion}\label{sec:ring}

% \begin{figure}[t!]
%   \centering
%   \includegraphics[width=0.95\textwidth]{pics/ring_all.eps}
%   \caption{Expanding Ring Fragments.}
%   \label{fig:glassdamage_mtd}
% \end{figure} 

 \begin{table}[]
\center
\begin{tabular}{|c|cccc|}
\hline
Material properties& Young's modulus $E$ & Poisson ratio $\nu$ & Density $\rho$& Fracture energy $G_0$\\
\hline
Value&$200 GPa$&$0.3$&$7800 kg/m^3$&  $1.125 \times 10^5 J/m^2$\\
\hline
\end{tabular}
\caption{Mechanical properties for the cylinder fragmentation under internal pressure example, following \cite{rabczuk2004cracking,abd2001fracture}.}
\label{tab:ring}
\end{table}

 \begin{figure}
     \centering
     \begin{subfigure}
     \centering
     \includegraphics[width = .30\linewidth]{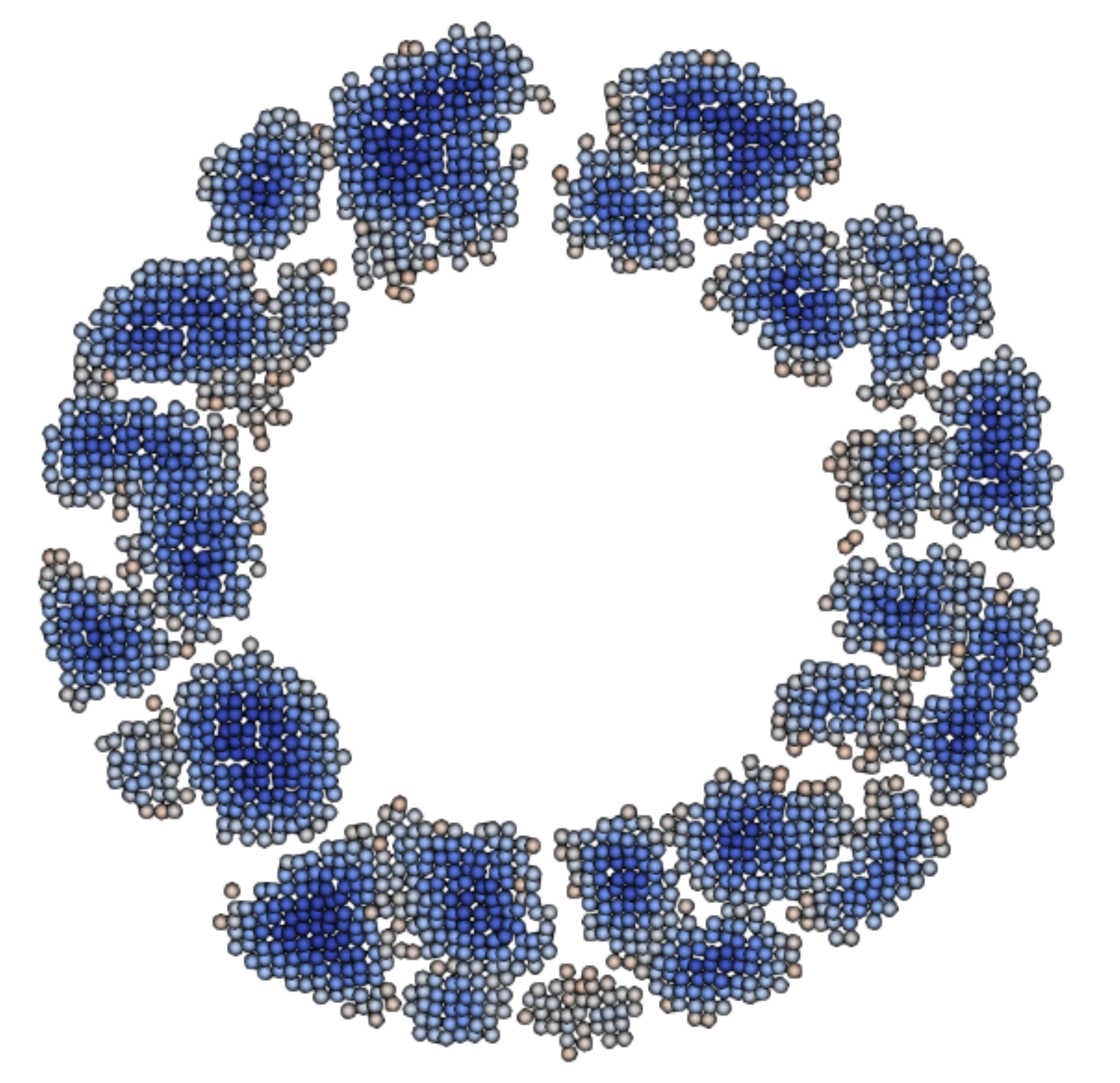}
     \end{subfigure}
     \qquad\qquad
     \begin{subfigure}
     \centering
     \includegraphics[width = .30\linewidth]{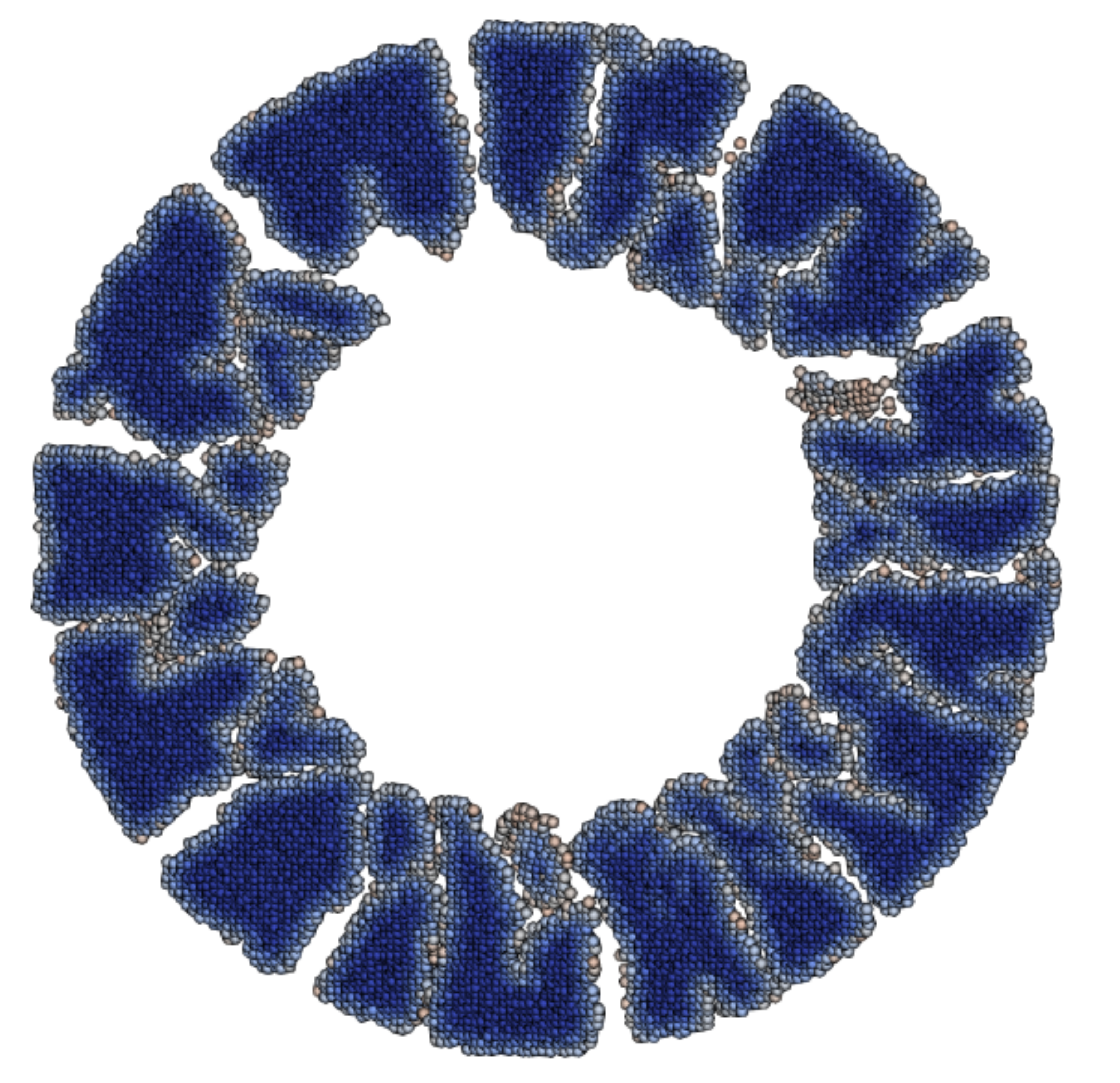}
     \end{subfigure}
     \begin{subfigure}
     \centering
     \includegraphics[width = .30\linewidth]{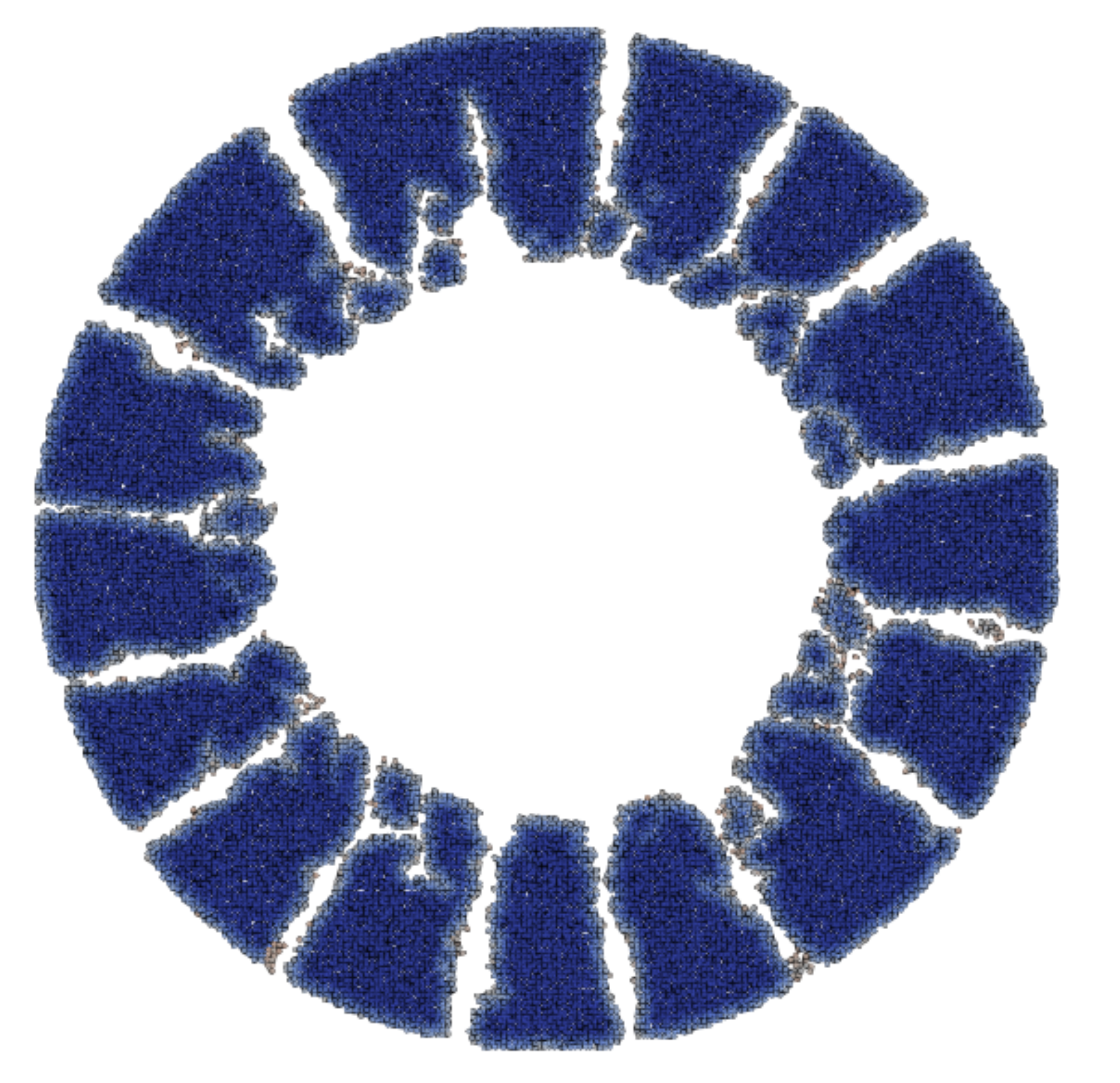}
     \end{subfigure}
     \qquad\qquad
     \begin{subfigure}
     \centering
     \includegraphics[width = .30\linewidth]{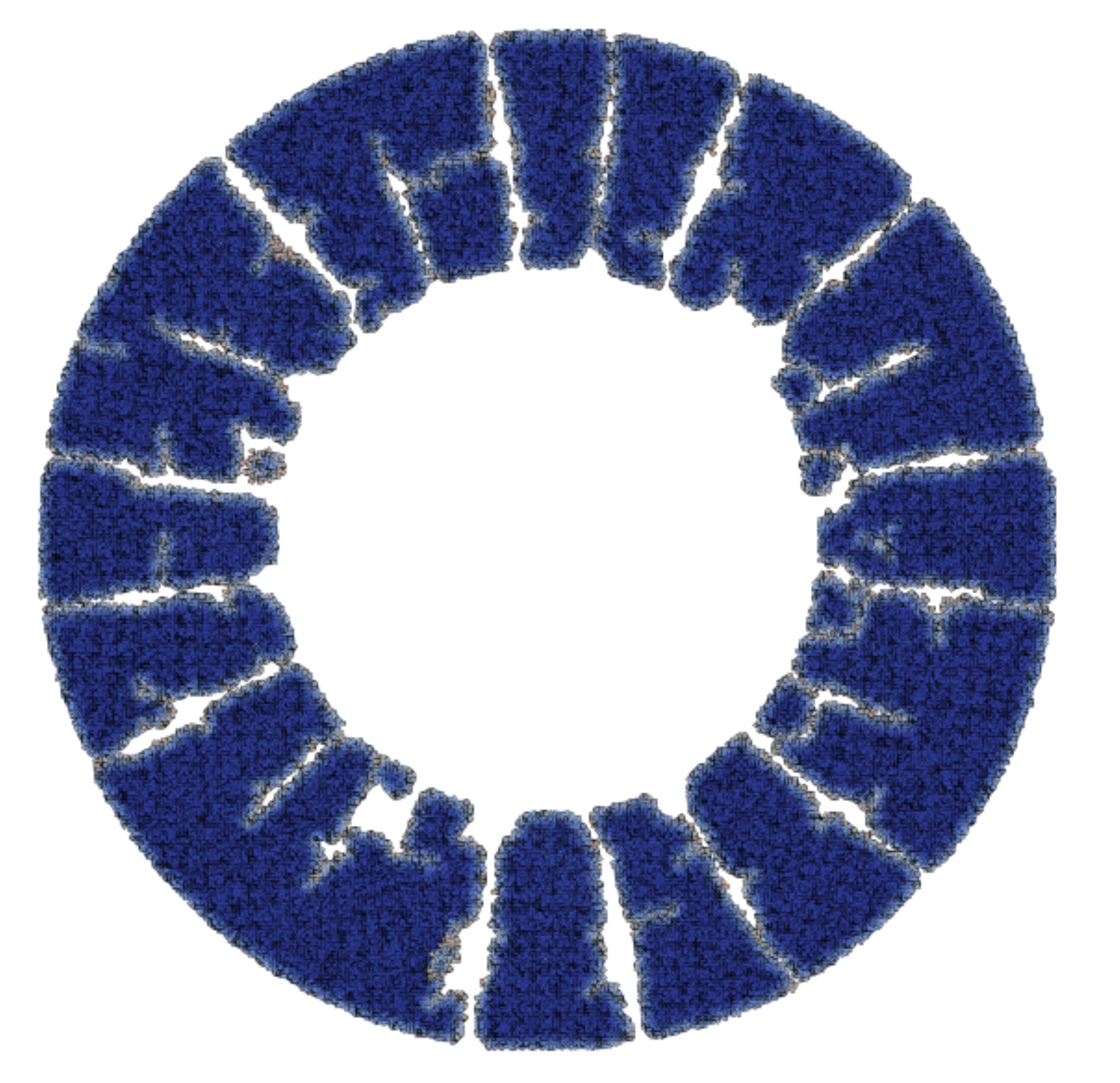}
     \end{subfigure}
     \caption{Predicted fragmentation of a cylinder under internal pressure at $T = 0.2 ms$, with different discretization resolutions. Left top: 3124 particles. Right top: 12587 particles. Left bottom: 22413 particles. Right bottom: 35035 particles.}
     \label{fig:ringex_simulation}
 \end{figure}
 
 \begin{table}
     \centering
     \begin{tabular}{ccc}
     \hline
        Number of particles  &  Number of large fragments & Number of small fragments \\
         \hline
         3124 &  12 & 2\\
         12587 & 15 & 3\\
         22413 & 15 & 9\\
         35035 & 14 & 4 \\
         \hline
     \end{tabular}
     \caption{Predicted number of large and small fragments in the cylinder under internal pressure simulation.}
     \label{tab:numer_frag}
 \end{table}

 In the last example simulation, we consider the fragmentation of a cylinder under internal pressure, so as to evaluate the proposed algorithm on handling multiple cracks and fragments. Following a similar setting as in \cite{rabczuk2004cracking}, a cylinder with inner radius $80 mm$ and outer radius $150 mm$ is employed, with material properties listed in Table \ref{tab:ring}. The cylinder is subject to a internal $p = p_0 e^{-t/t_0}$, where $p_0 = 2.5 {GPa}$, and $t_0 = 0.01 {ms}$. We run the simulation with $\delta=4h$, $\Delta t =0.05\mu s$, and four different levels of spatial resolution: 3124, 12587, 22413, and 35035 discretization points (particles). For each set of discretization points, we generate non-uniform grids by perturbing particle positions by a uniformly distributed perturbation of magnitude $0.2h$. In Figure \ref{fig:ringex_simulation}, we show simulation results at $T = 0.2 {ms}$, when the cylinder breaks into fragments. The number of large and small fragments are listed in Table \ref{tab:numer_frag}, where we can observe that the number of large fragments is generally consistent except for the case with the coarsest resolution. The number of small fragments generally increases when using finer discretizations but it's not monotonic. These observations as well as the number of large fragments are consistent with the simulation results in \cite{rabczuk2004cracking}, where a particle model was employed and 15-16 numbers of large fragments were predicted in numerical simulations with $12500-39000$ particles. This suggests the current scheme is appropriate for handling blast loading predictions, and provides consistent predictions as resolution is refined.

 \section{Conclusion and Future Work}\label{sec:conclusion}
 
 Peridynamics presents a flexible framework for modeling fracture mechanics. In particular, bond-based fracture models admit a sharp representation of fracture surfaces while avoiding the loss of mass associated with damage models and element death \cite{bathe1978some}. This flexibility comes with a cost however, as the free-surface introduced during fracture compounds traditional challenges in peridynamic models related to nonlocal boundary conditions. This work has presented a complete workflow demonstrating for linearly elastic material how quadrature, boundary and traction loading may be handled in such a way that one preserves a limit to the relevant local problem as resolution is increased. This is a major contribution to the field of peridynamics - while numerous works have demonstrated the flexibility of peridynamics in modeling a diverse set of physical phenomena, comparatively few have demonstrated rigorous notions of convergence and grid independence. Rigorous accuracy guarantees are fundamental to trusting predictions made by numerical models, and this work aims to provide an important first step toward putting peridynamics on the same footing as e.g. finite element methods for local mechanics.
 
 The primary focus of this work has been to establish schemes, quadrature rules, and boundary treatment and provide rigorous mathematical analysis. While numerical examples have been provided at a level appropriate for establishing the scheme's feasibility for practical problems, an important next step is to generate a performant parallel implementation allowing one to consider high-resolution predictions in two and three dimensions. For several of the validation studies provided here we were unable to reach the resolution used by other state-of-the-art peridynamic discretizations due to memory limitations of our serial implementation. The method itself is embarassingly parallelizable, as the generation of quadrature weights and dilitation corrections involves only the local construction and inversion of small linear matrices. In an upcoming work we will provide a clear demonstration of how the convergence guarantees provided by our approach translates to improved prediction accuracy for realistic problems. We will additionally consider the generalization of this approach to nonlinear elastoplasticity governing ductile failure.

\section*{Acknowledgements}
Sandia National Laboratories is a multi-mission laboratory managed and operated by National Technology and Engineering Solutions of Sandia, LLC., a wholly owned subsidiary of Honeywell International, Inc., for the U.S. Department of Energy’s National Nuclear Security Administration under contract DE-NA0003525. This paper describes objective technical results and analysis. Any subjective views or opinions that might be expressed in the paper do not necessarily represent the views of the U.S. Department of Energy or the United States Government.

H. You and Y. Yu are supported by the National Science Foundation under award DMS 1753031. N. Trask's work is supported under the Sandia National Laboratories Laboratory Directed Research and Development (LDRD) program, and by the U.S. Department of Energy, Office of Science, Office of Advanced Scientific Computing Research under the Collaboratory on Mathematics and Physics-Informed Learning Machines for Multiscale and Multiphysics Problems (PhILMs) project. SAND number: SAND2021-0063 O

\appendix
\section{Convergence studies}\label{app:convergence}

\subsection{Linear Patch Tests}\label{sec:patchapp}

\begin{figure}[!htb]\centering
 \subfigure{\includegraphics[width=0.49\textwidth]{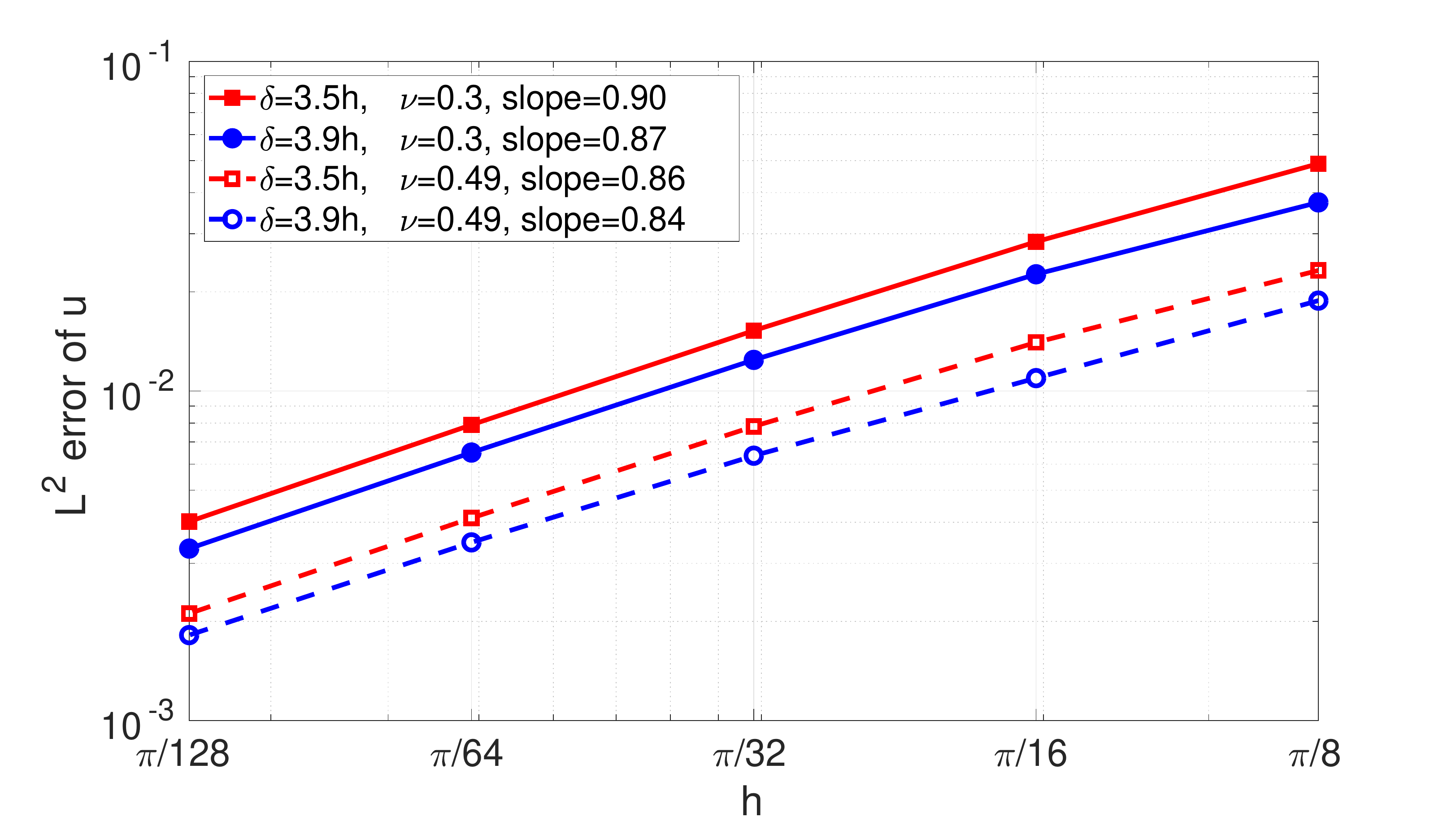}}
  \subfigure{\includegraphics[width=0.49\textwidth]{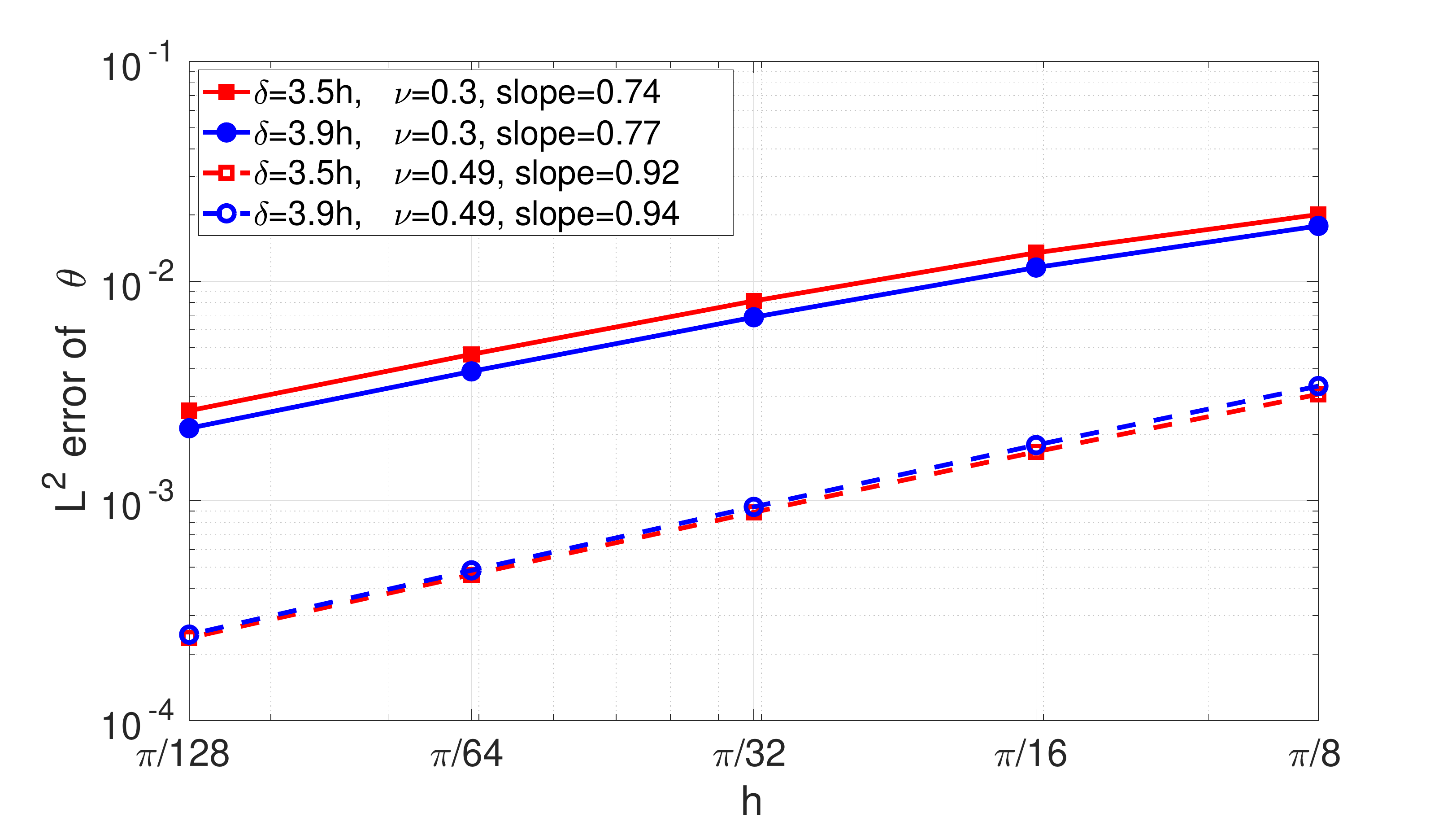}}
 \caption{Linear patch test with uniform discretizations: traction loads applied on boundary including a corner. Left: $L^2(\omg)$ errors of displacement $\ub$. Right: $L^2(\omg)$ errors of dilatation $\theta$.
 }
 \label{fig:patch_corner_s}
\end{figure}

 \begin{figure}[!htb]\centering
 \subfigure{\includegraphics[width=0.49\textwidth]{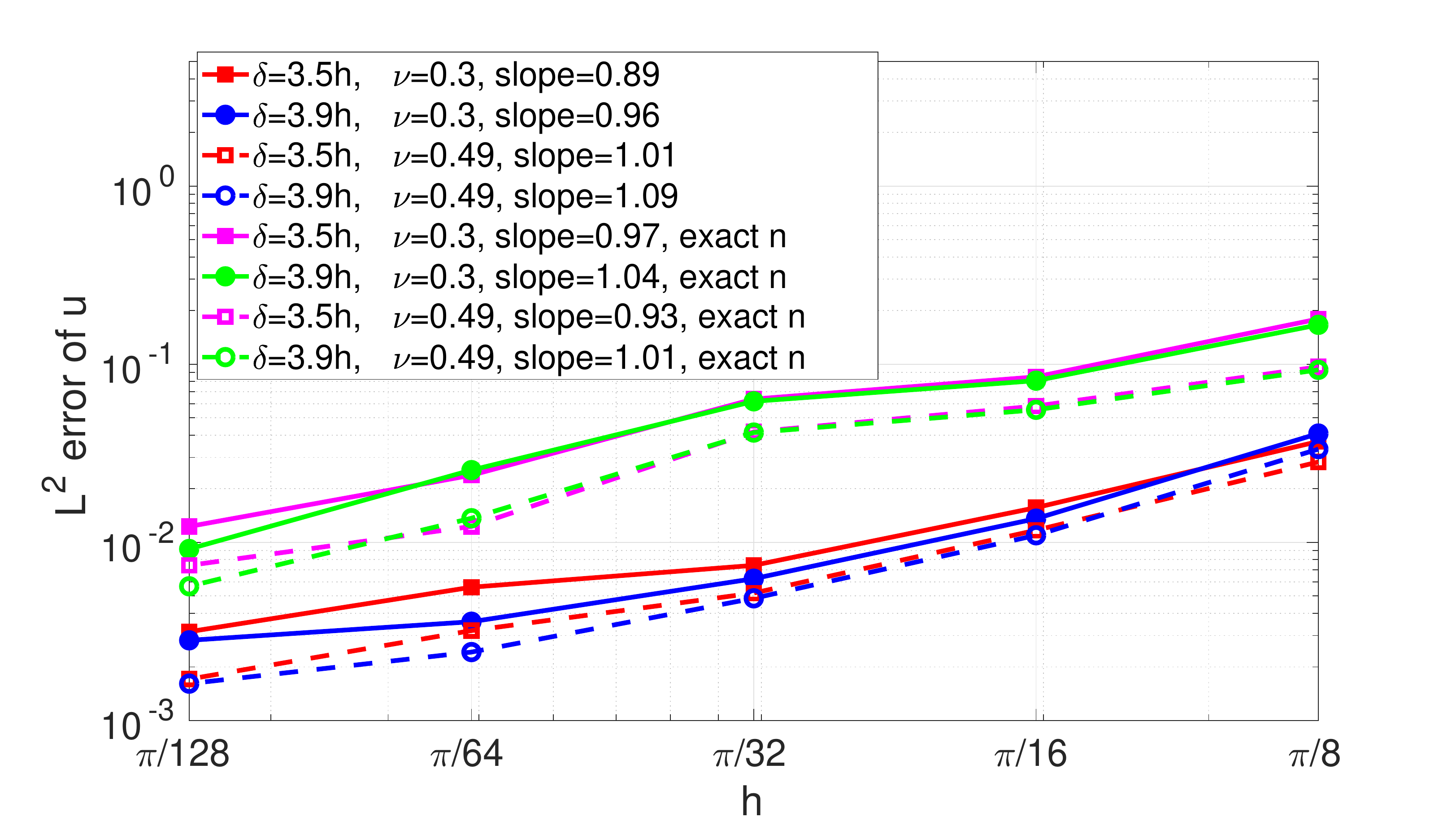}}
  \subfigure{\includegraphics[width=0.49\textwidth]{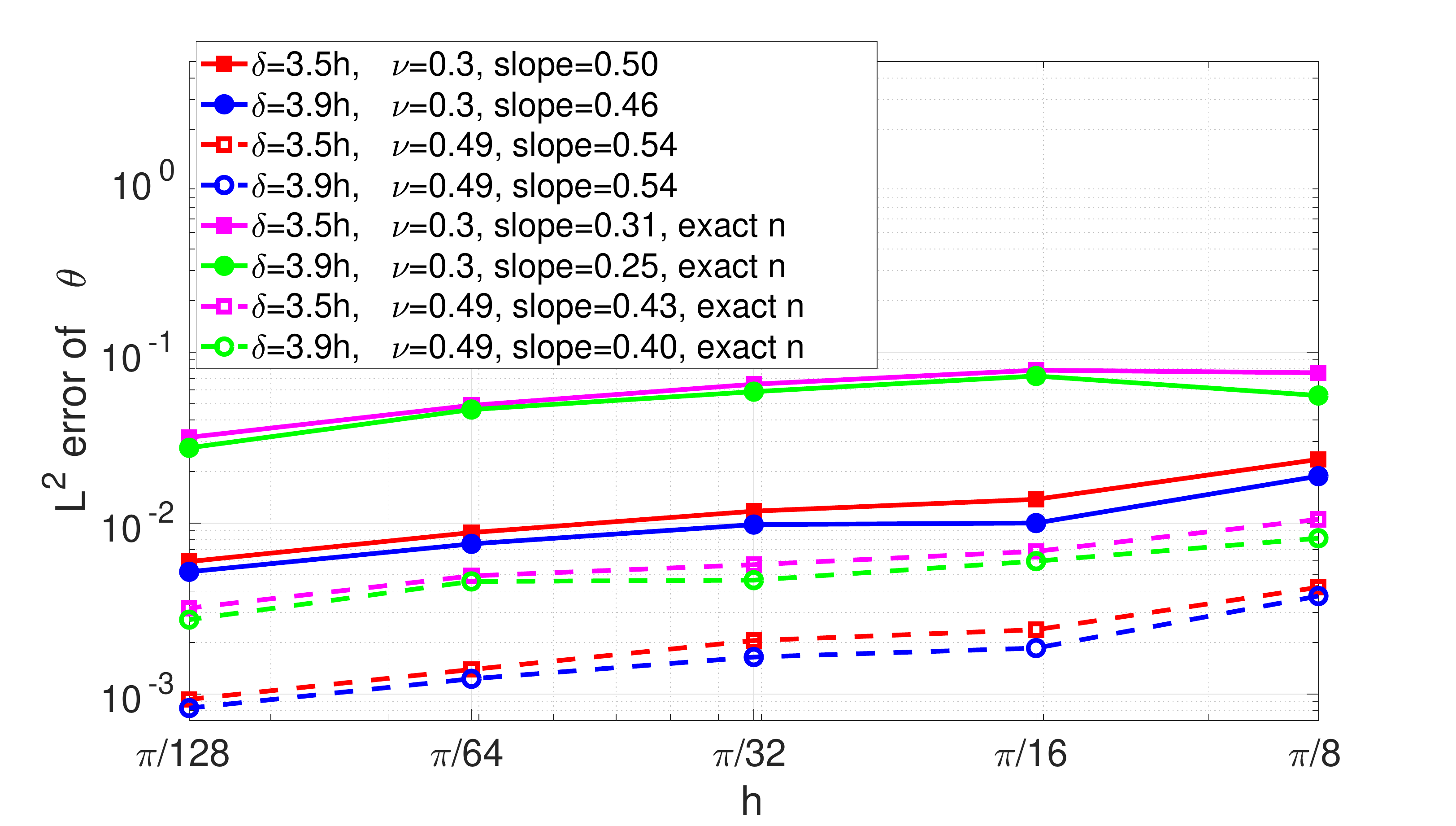}}
 \caption{Linear patch test with non-uniform discretizations: traction loads applied on a straight line. Left: $L^2(\omg)$ errors of displacement $\ub$. Right: $L^2(\omg)$ errors of dilatation $\theta$.
 }
 \label{fig:patch_linear_un}
\end{figure}

 \begin{figure}[!htb]\centering
 \subfigure{\includegraphics[width=0.49\textwidth]{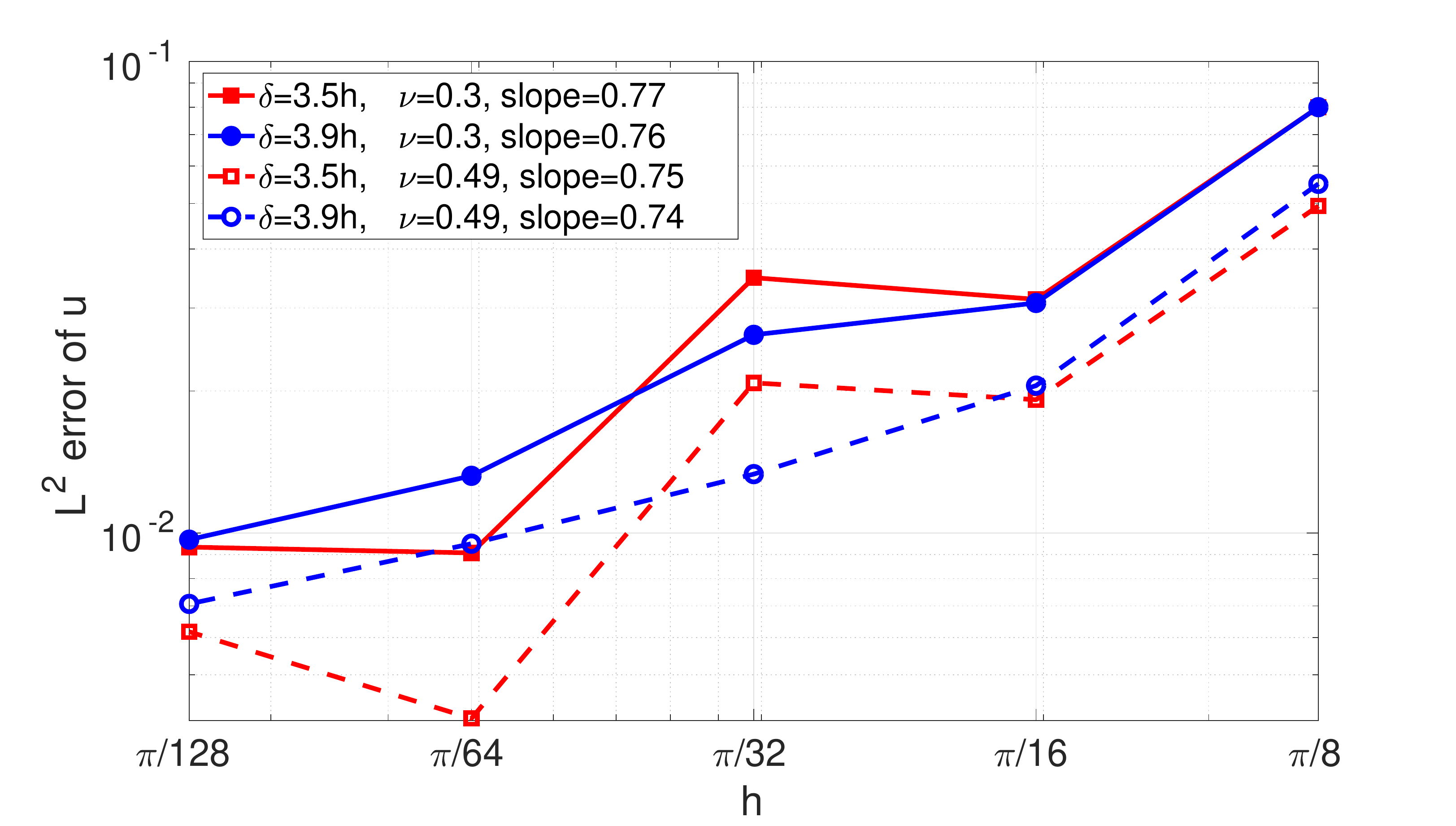}}
  \subfigure{\includegraphics[width=0.49\textwidth]{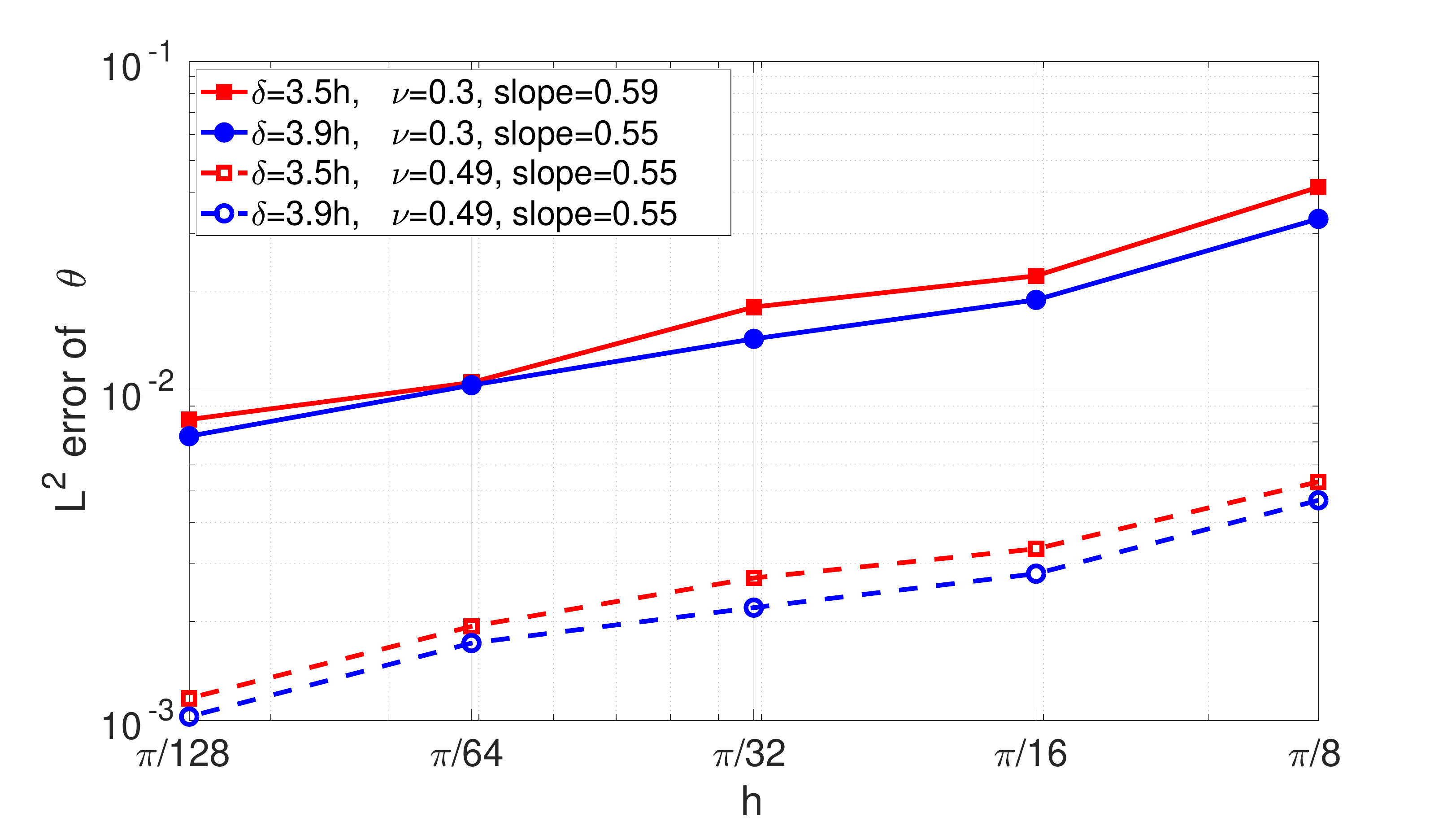}}
 \caption{Linear patch test with non-uniform discretizations: traction loads applied on boundary including a corner. Left: $L^2(\omg)$ errors of displacement $\ub$. Right: $L^2(\omg)$ errors of dilatation $\theta$.
 }
 \label{fig:patch_corner_un}
\end{figure}

\subsection{Manufactured solution test}

 \begin{figure}[!htb]\centering
 \subfigure{\includegraphics[width=0.49\textwidth]{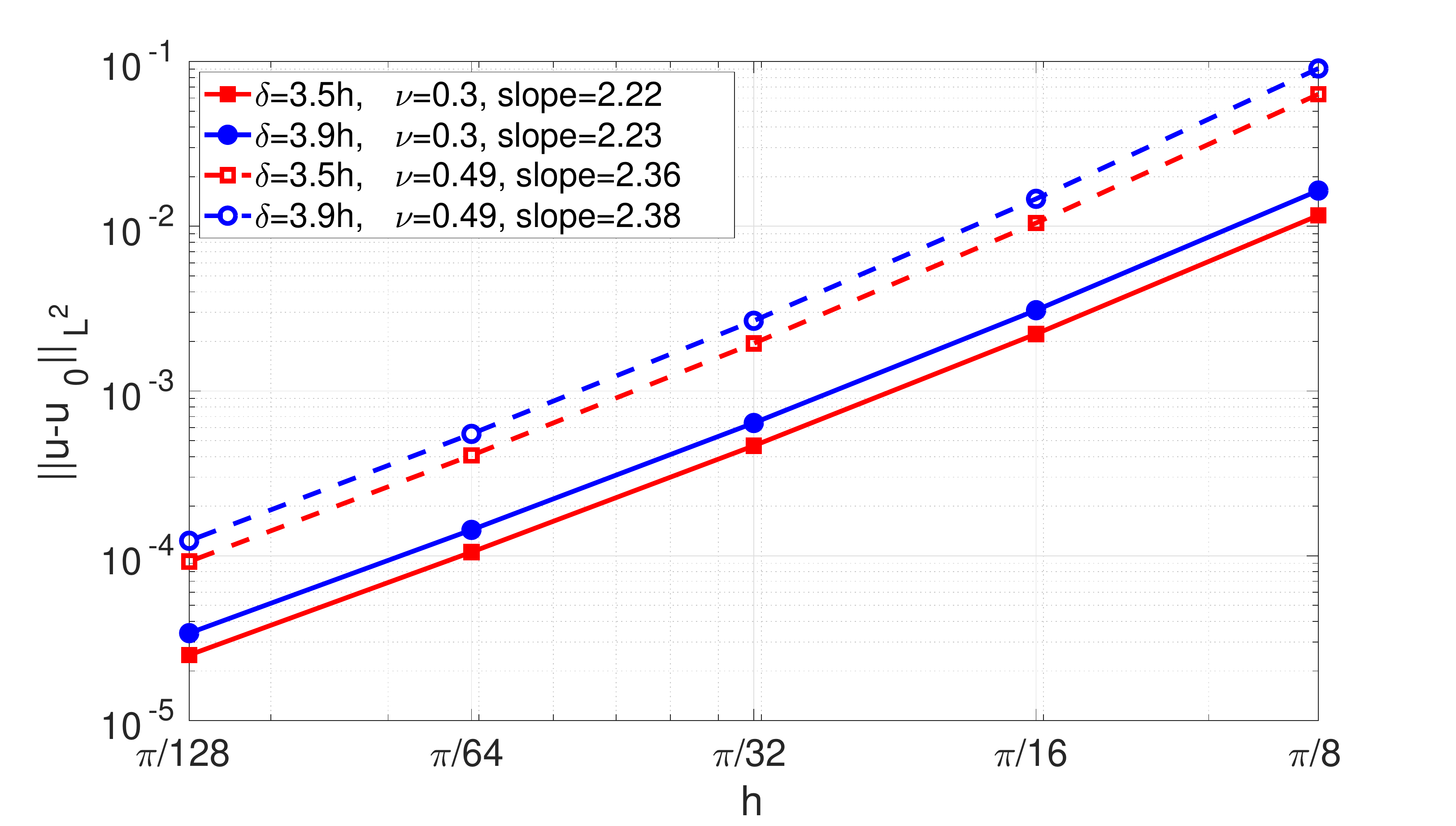}}
  \subfigure{\includegraphics[width=0.49\textwidth]{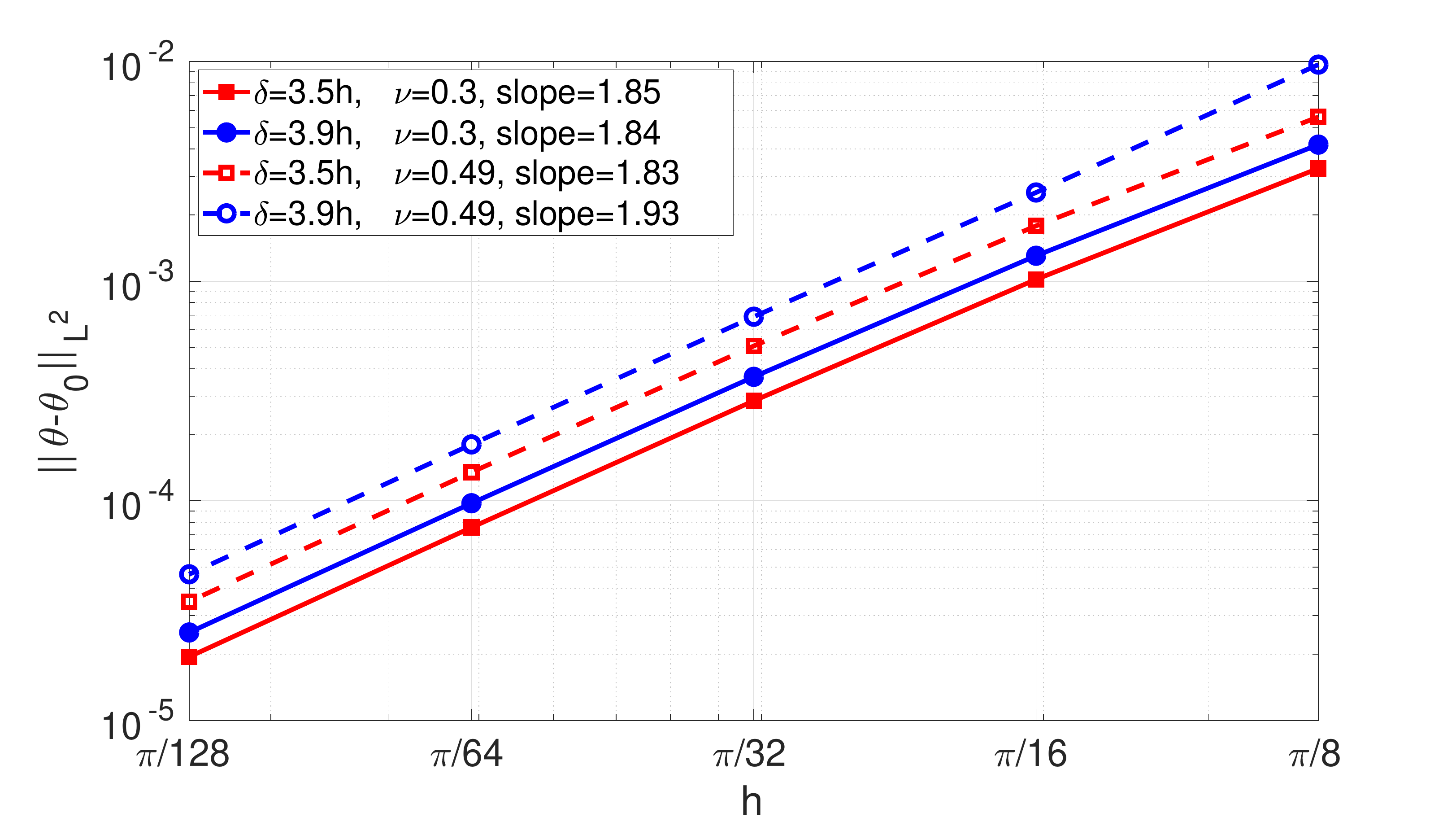}}
 \caption{M-convergence tests on a square domain with uniform discretizations and full Dirichlet-type boundary conditions. Left: the $L^2(\omg)$ difference between displacement $\ub$ and its local limit $\ub_0$. Right: the $L^2(\omg)$ difference between the nonlocal dilitation $\theta$ and its local limit $\theta_0=\nabla\cdot\ub_0$.
 }
 \label{fig:sin_Diri_s}
\end{figure}

 \begin{figure}[!htb]\centering
 \subfigure{\includegraphics[width=0.49\textwidth]{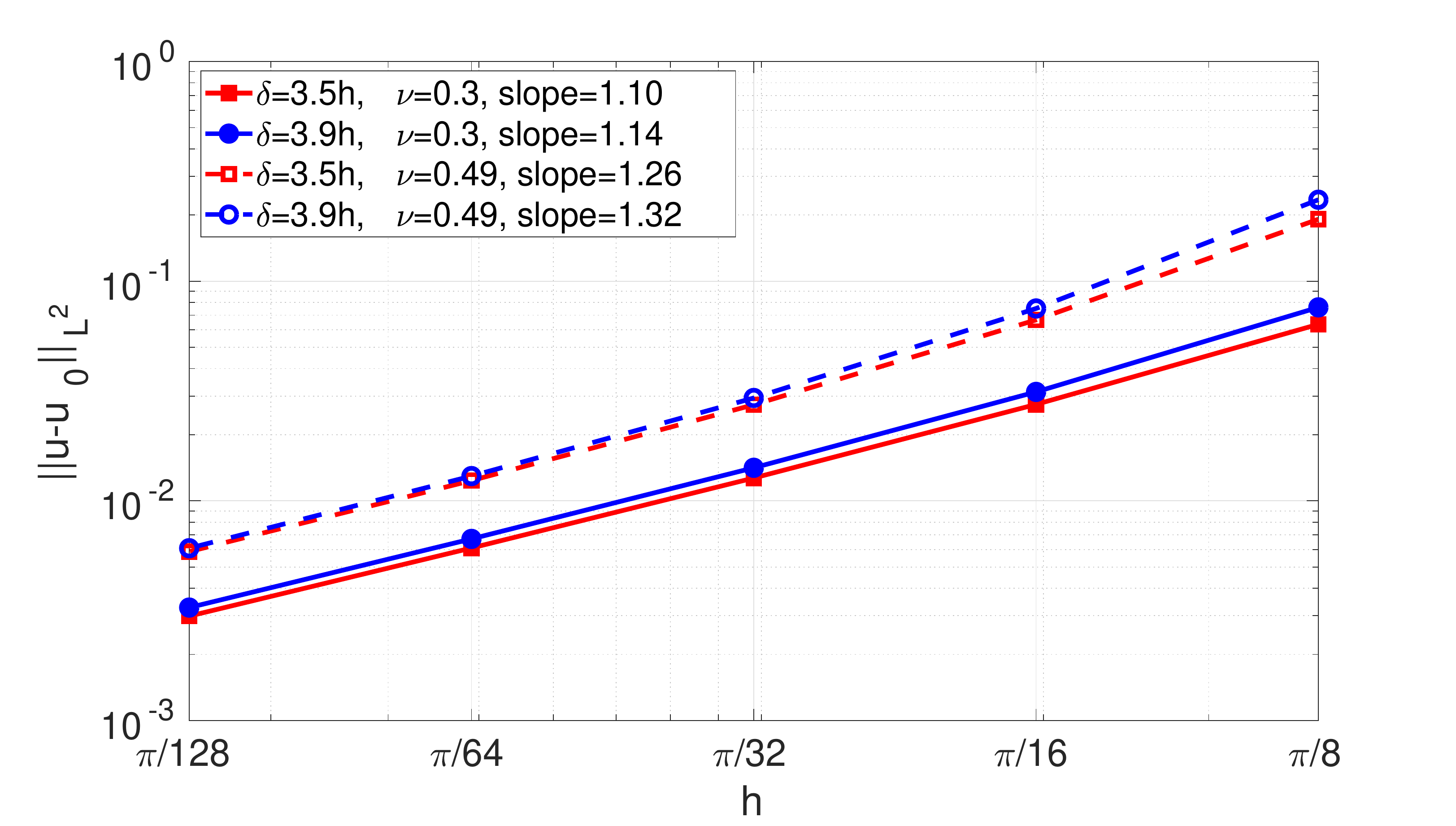}}
  \subfigure{\includegraphics[width=0.49\textwidth]{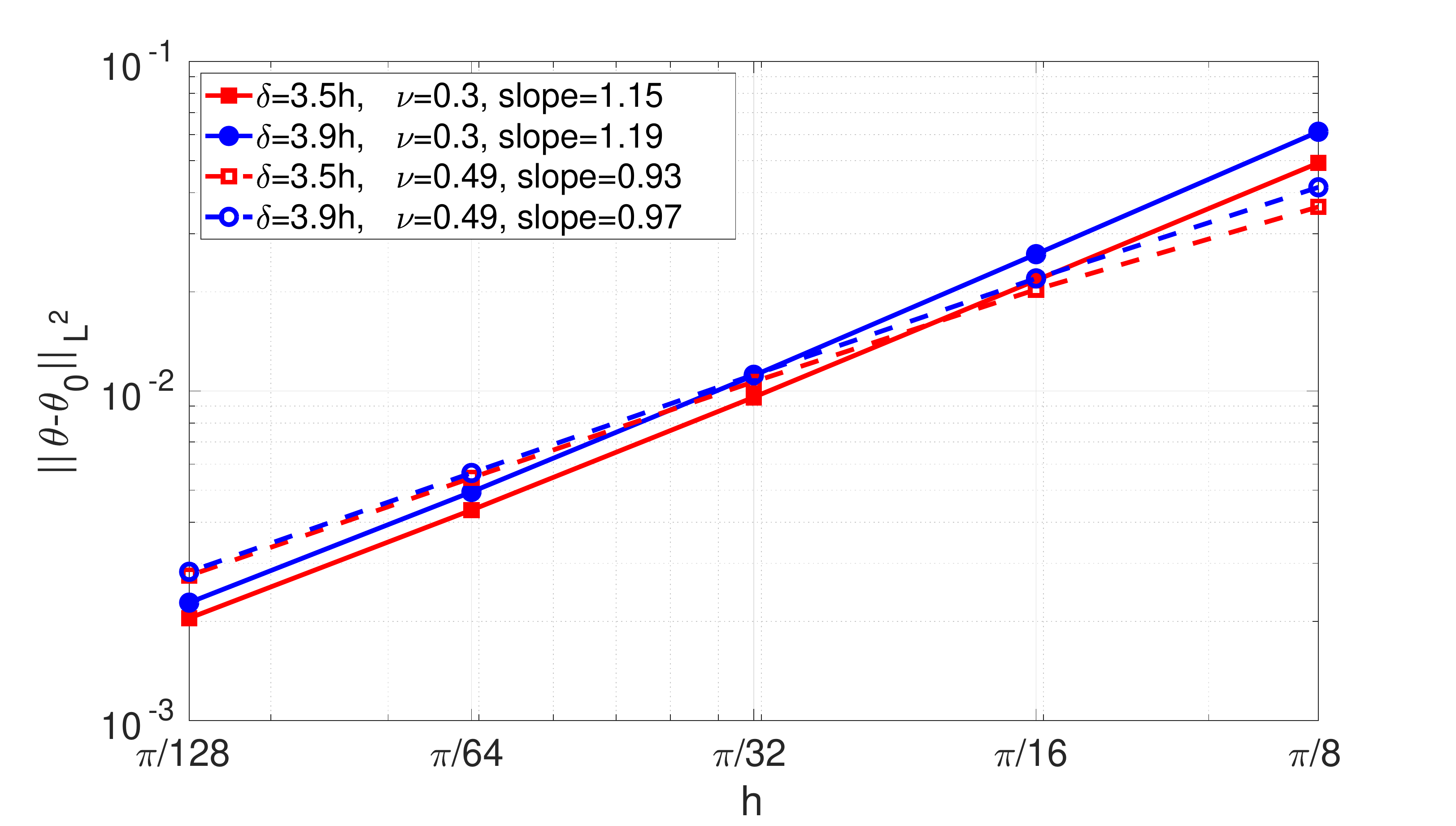}}
 \caption{M-convergence tests on a square domain with uniform discretizations and traction loads applied on a straight line. Left: the $L^2(\omg)$ difference between displacement $\ub$ and its local limit $\ub_0$. Right: the $L^2(\omg)$ difference between the nonlocal dilitation $\theta$ and its local limit $\theta_0=\nabla\cdot\ub_0$.
 }
 \label{fig:sin_linear_s}
\end{figure}

 \begin{figure}[!htb]\centering
 \subfigure{\includegraphics[width=0.49\textwidth]{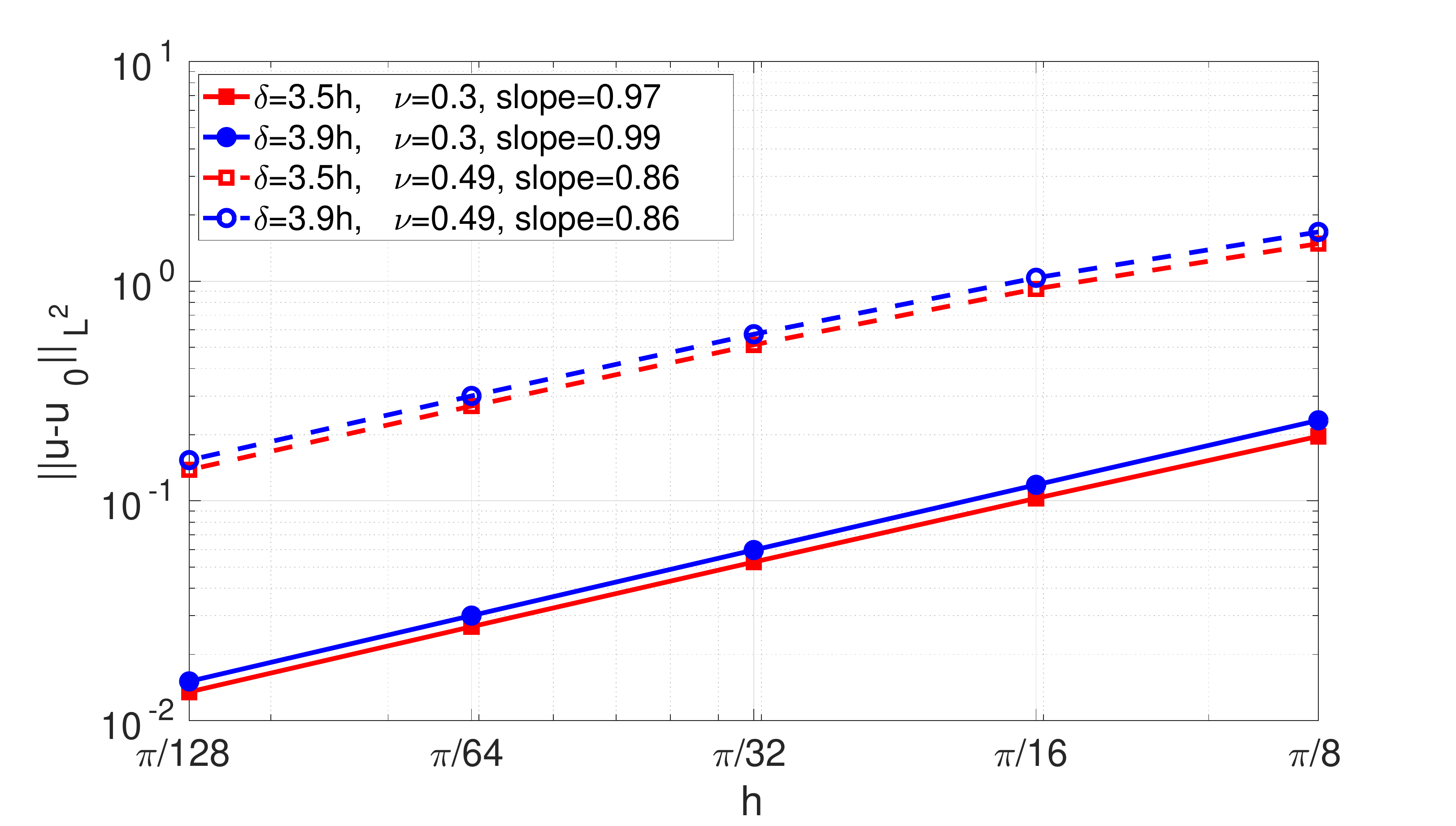}}
  \subfigure{\includegraphics[width=0.49\textwidth]{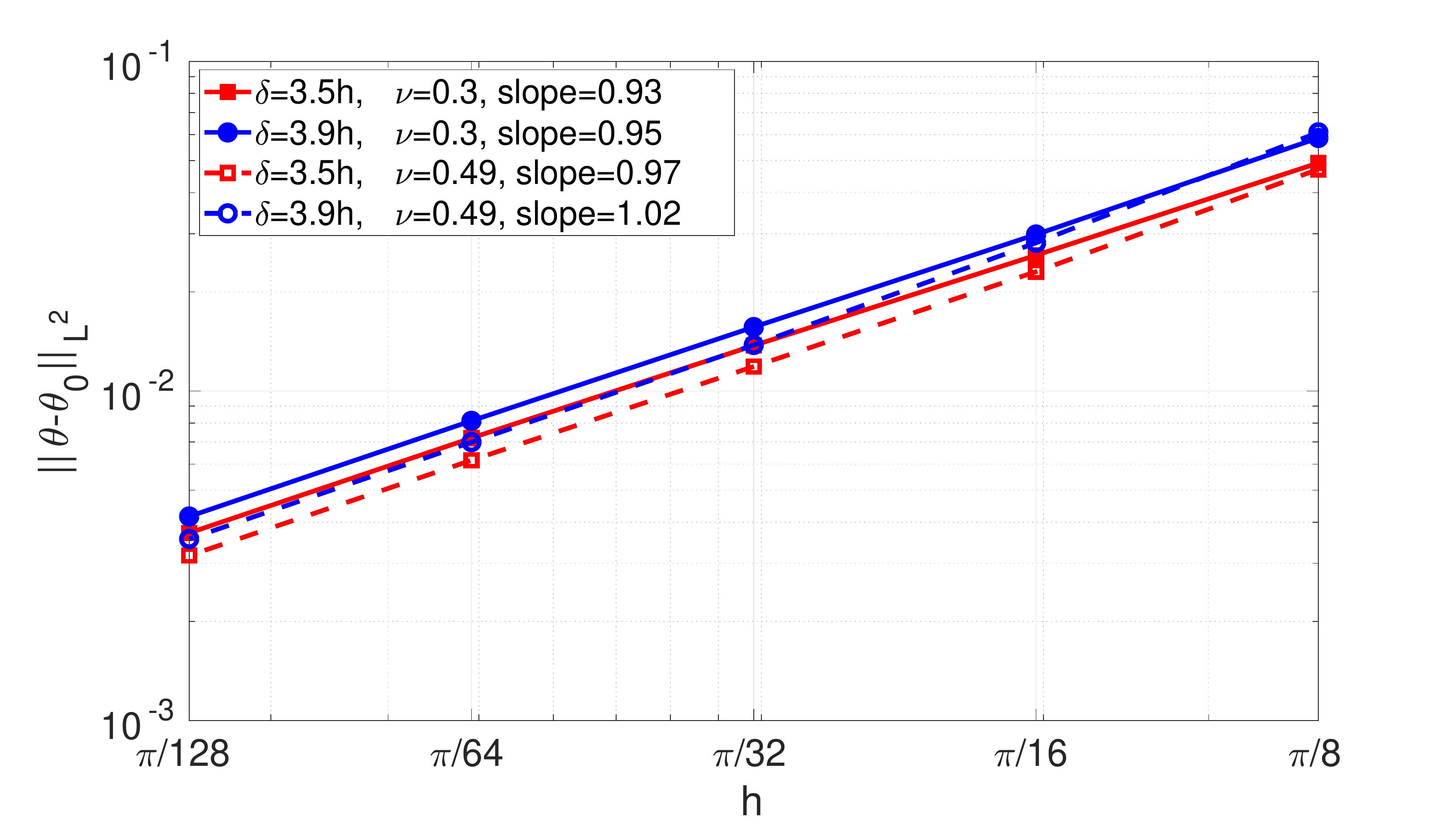}}
 \caption{M-convergence tests on a square domain with uniform discretizations and traction loads applied on boundary including a corner. Left: the $L^2(\omg)$ difference between displacement $\ub$ and its local limit $\ub_0$. Right: the $L^2(\omg)$ difference between the nonlocal dilitation $\theta$ and its local limit $\theta_0=\nabla\cdot\ub_0$.
 }
 \label{fig:sin_corner_s}
\end{figure}

 \begin{figure}[!htb]\centering
 \subfigure{\includegraphics[width=0.49\textwidth]{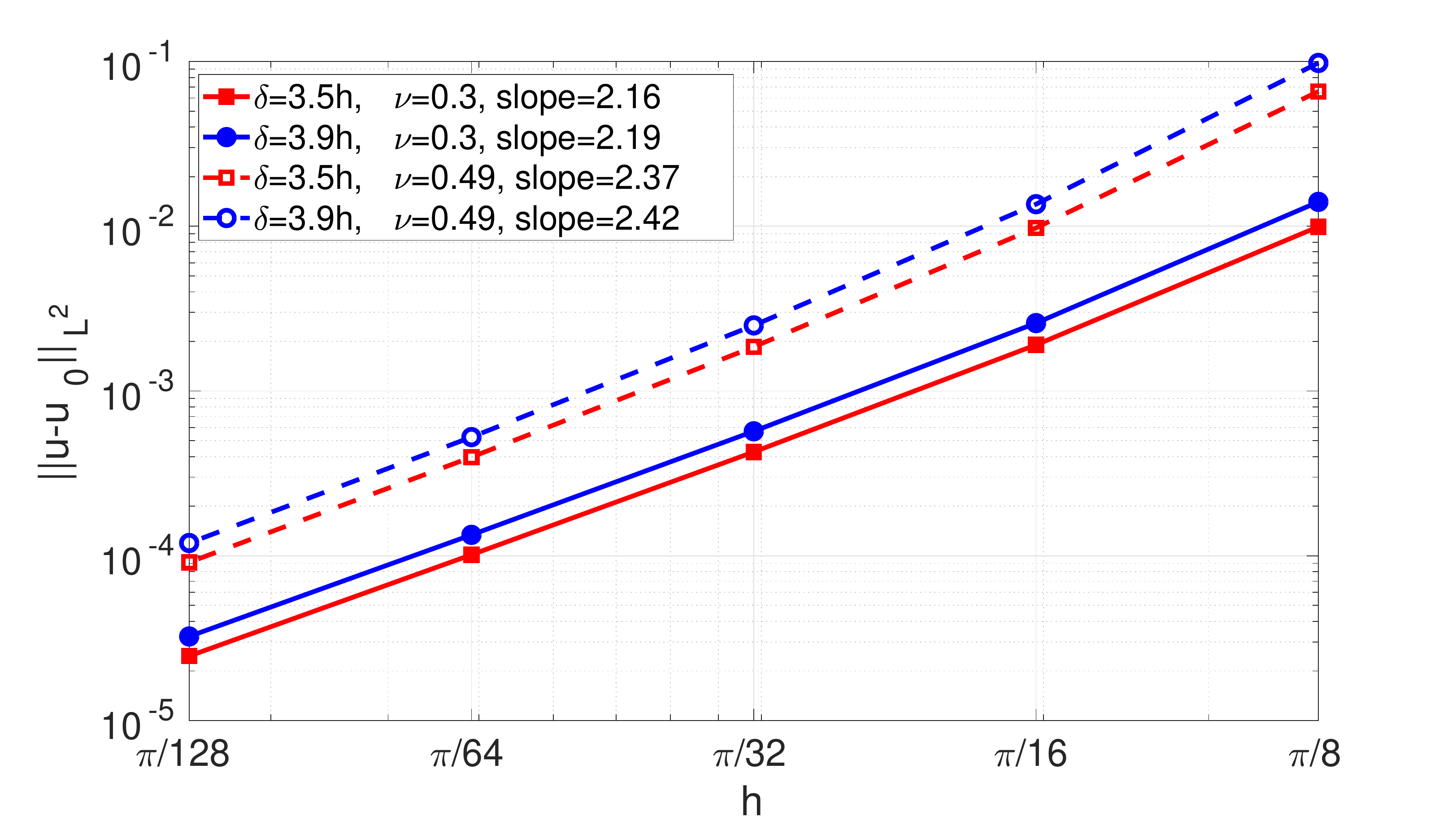}}
  \subfigure{\includegraphics[width=0.49\textwidth]{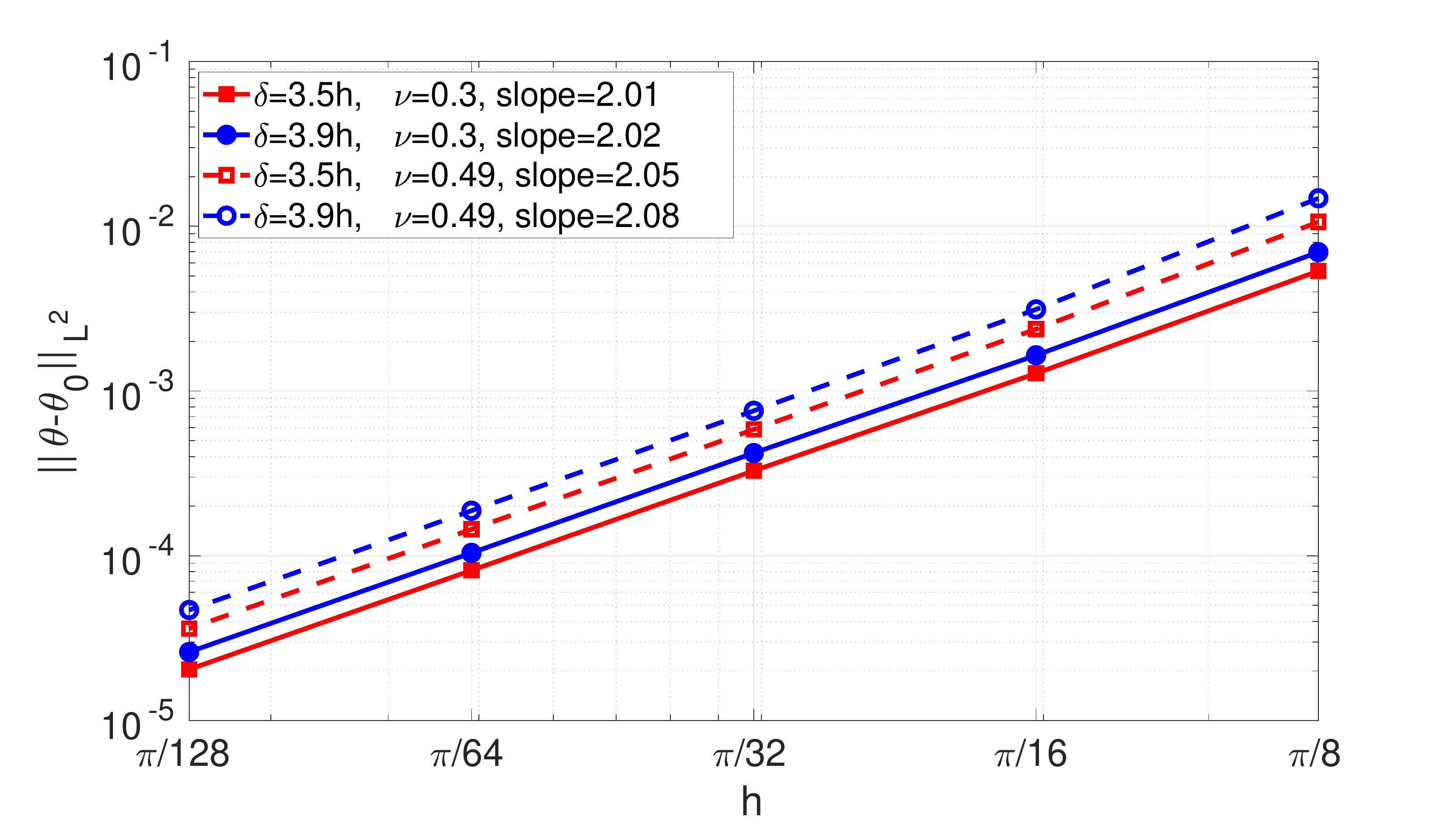}}
 \caption{M-convergence tests on a square domain with non-uniform discretizations and full Dirichlet-type boundary conditions. Left: the $L^2(\omg)$ difference between displacement $\ub$ and its local limit $\ub_0$. Right: the $L^2(\omg)$ difference between the nonlocal dilitation $\theta$ and its local limit $\theta_0=\nabla\cdot\ub_0$.
 }
 \label{fig:sin_Diri_un}
\end{figure}

 \begin{figure}[!htb]\centering
 \subfigure{\includegraphics[width=0.49\textwidth]{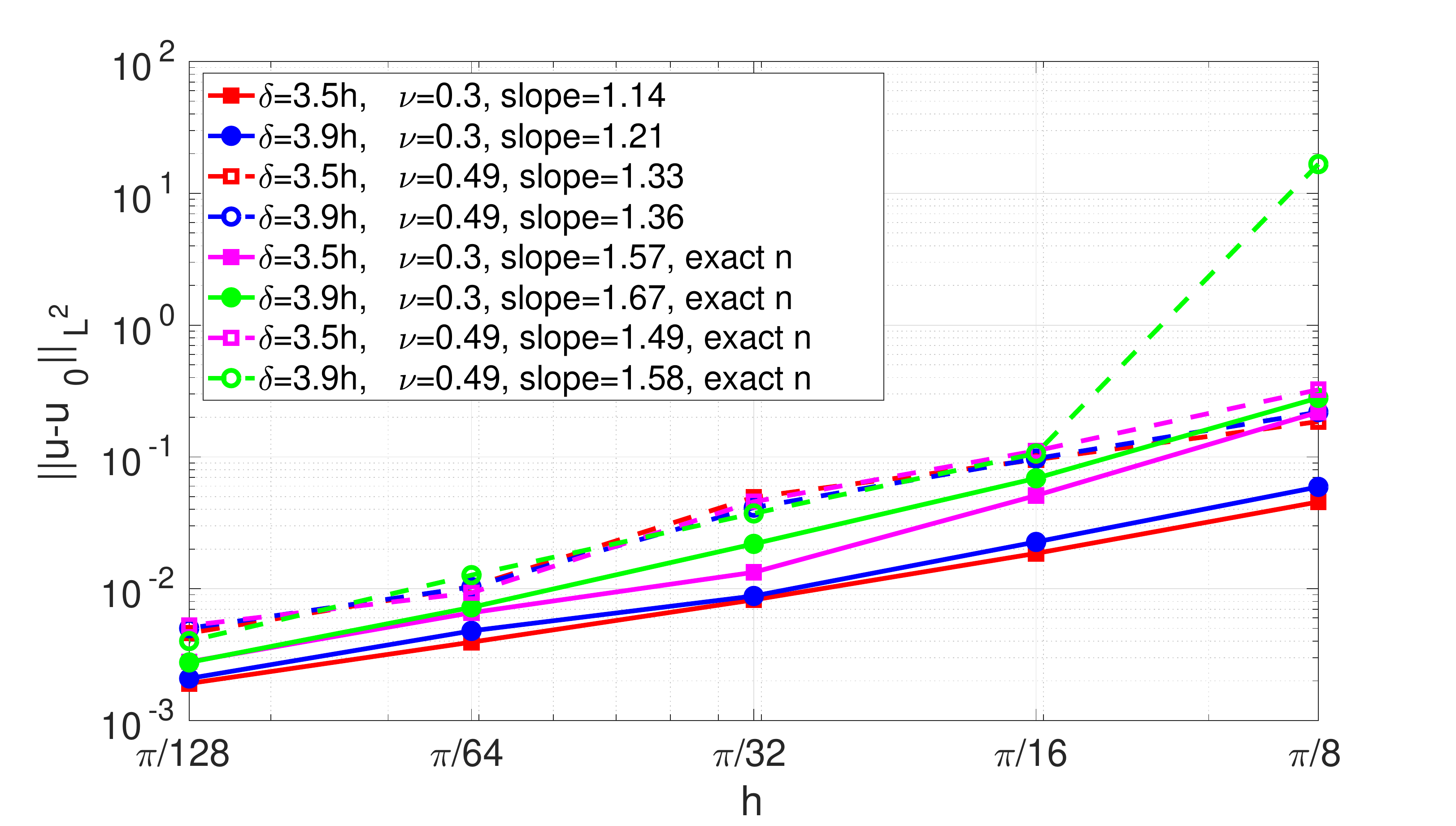}}
  \subfigure{\includegraphics[width=0.49\textwidth]{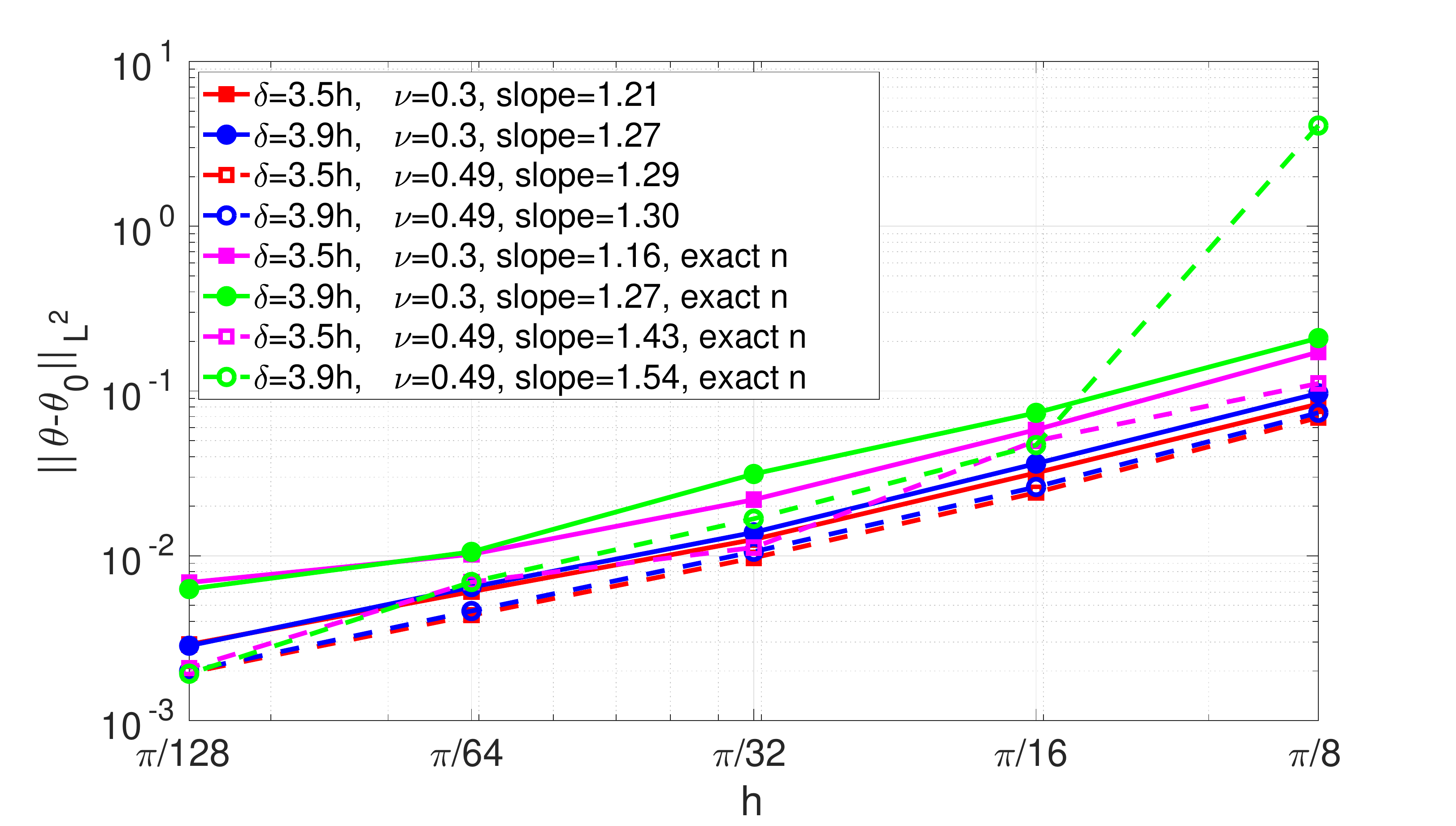}}
 \caption{M-convergence tests on a square domain with non-uniform discretizations and traction loads applied on a straight line. Left: the $L^2(\omg)$ difference between displacement $\ub$ and its local limit $\ub_0$. Right: the $L^2(\omg)$ difference between the nonlocal dilitation $\theta$ and its local limit $\theta_0=\nabla\cdot\ub_0$.
 }
 \label{fig:sin_linear_un}
\end{figure}

 \begin{figure}[!htb]\centering
 \subfigure{\includegraphics[width=0.49\textwidth]{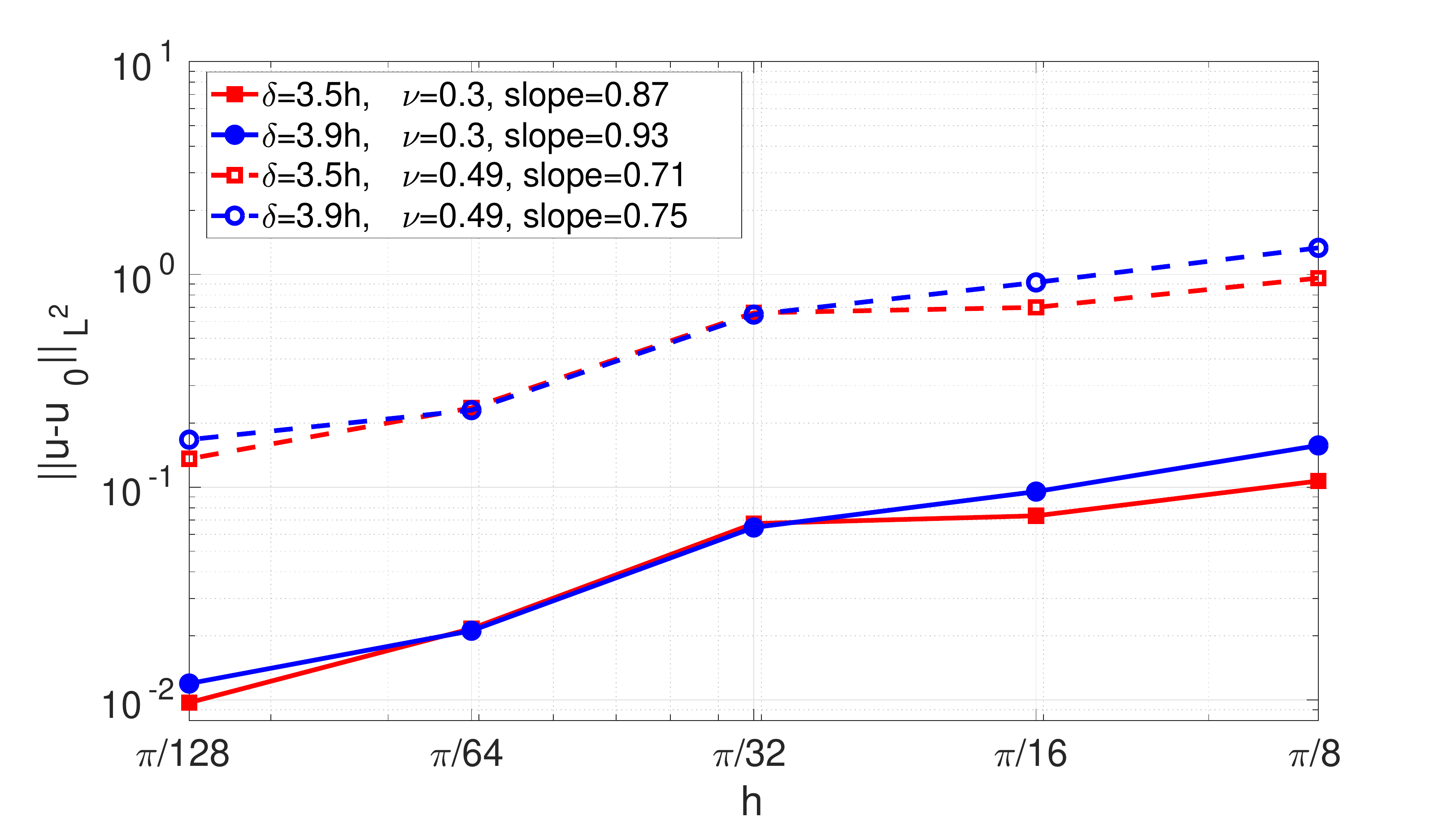}}
  \subfigure{\includegraphics[width=0.49\textwidth]{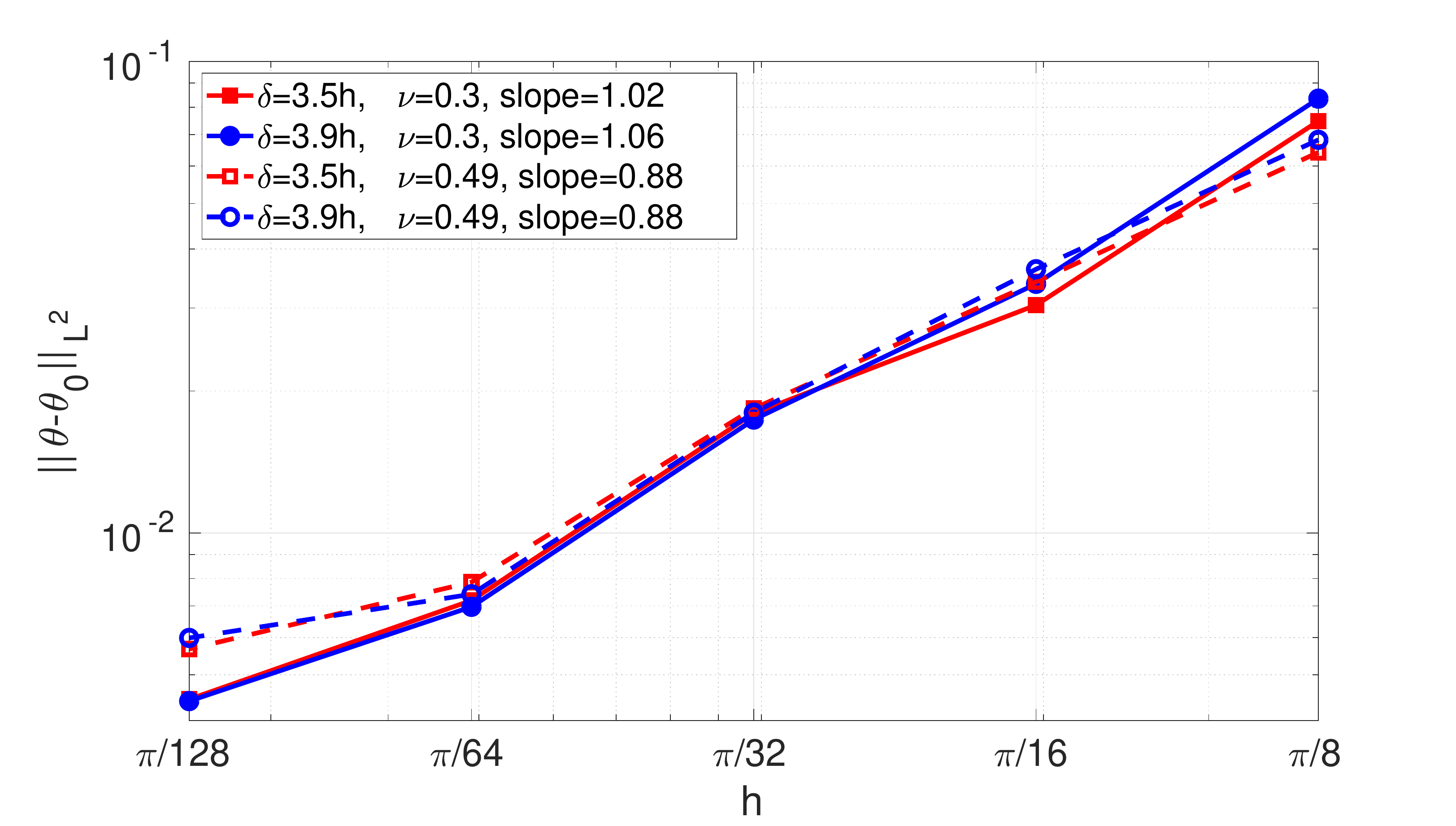}}
 \caption{M-convergence tests on a square domain with non-uniform discretizations and traction loads applied on boundary including a corner. Left: the $L^2(\omg)$ difference between displacement $\ub$ and its local limit $\ub_0$. Right: the $L^2(\omg)$ difference between the nonlocal dilitation $\theta$ and its local limit $\theta_0=\nabla\cdot\ub_0$.
 }
 \label{fig:sin_corner_un}
\end{figure}

\subsection{Traction loads on curvilinear free surfaces}

 \begin{figure}[!htb]\centering
 \subfigure{\includegraphics[width=0.49\textwidth]{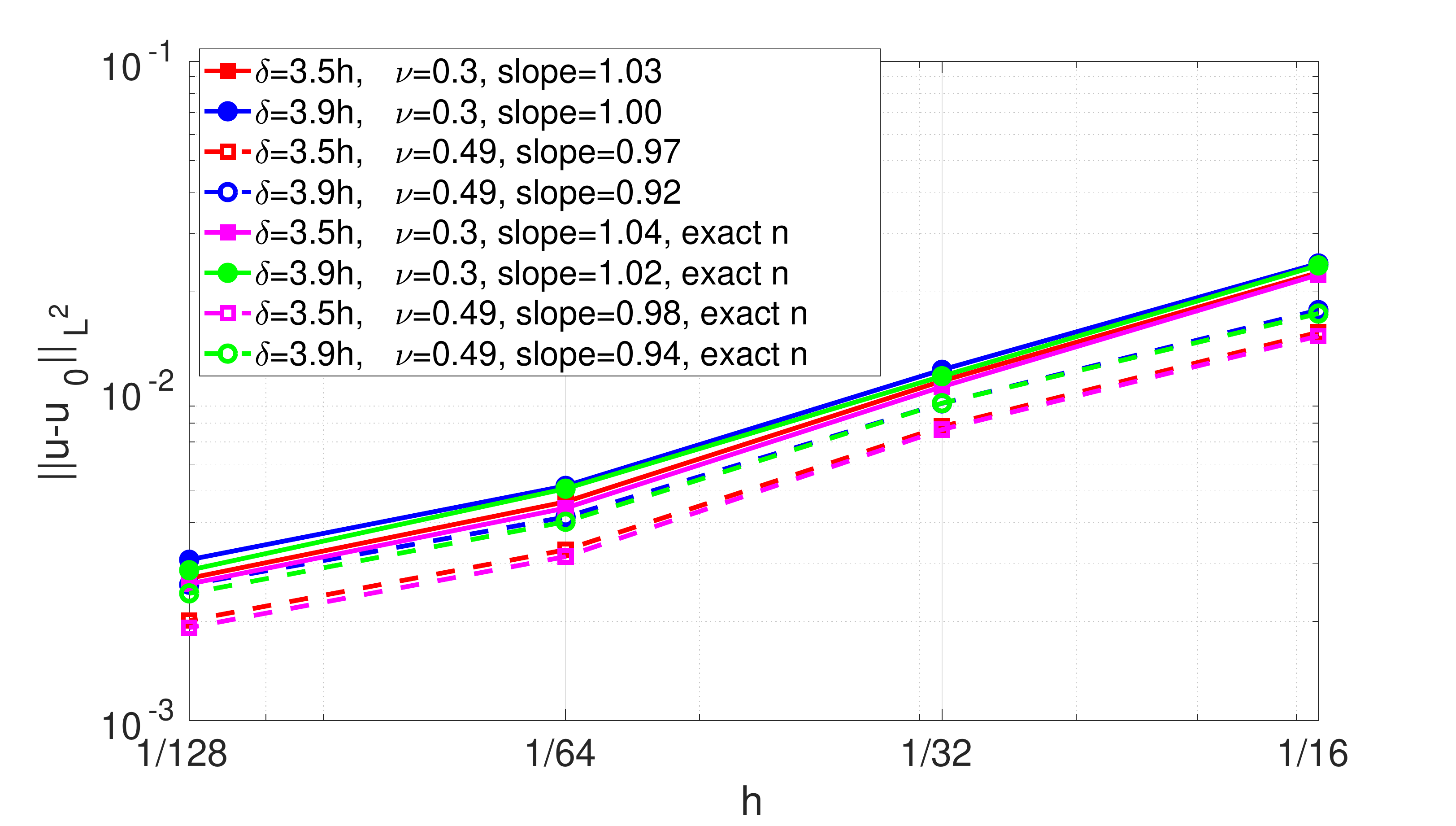}}
  \subfigure{\includegraphics[width=0.49\textwidth]{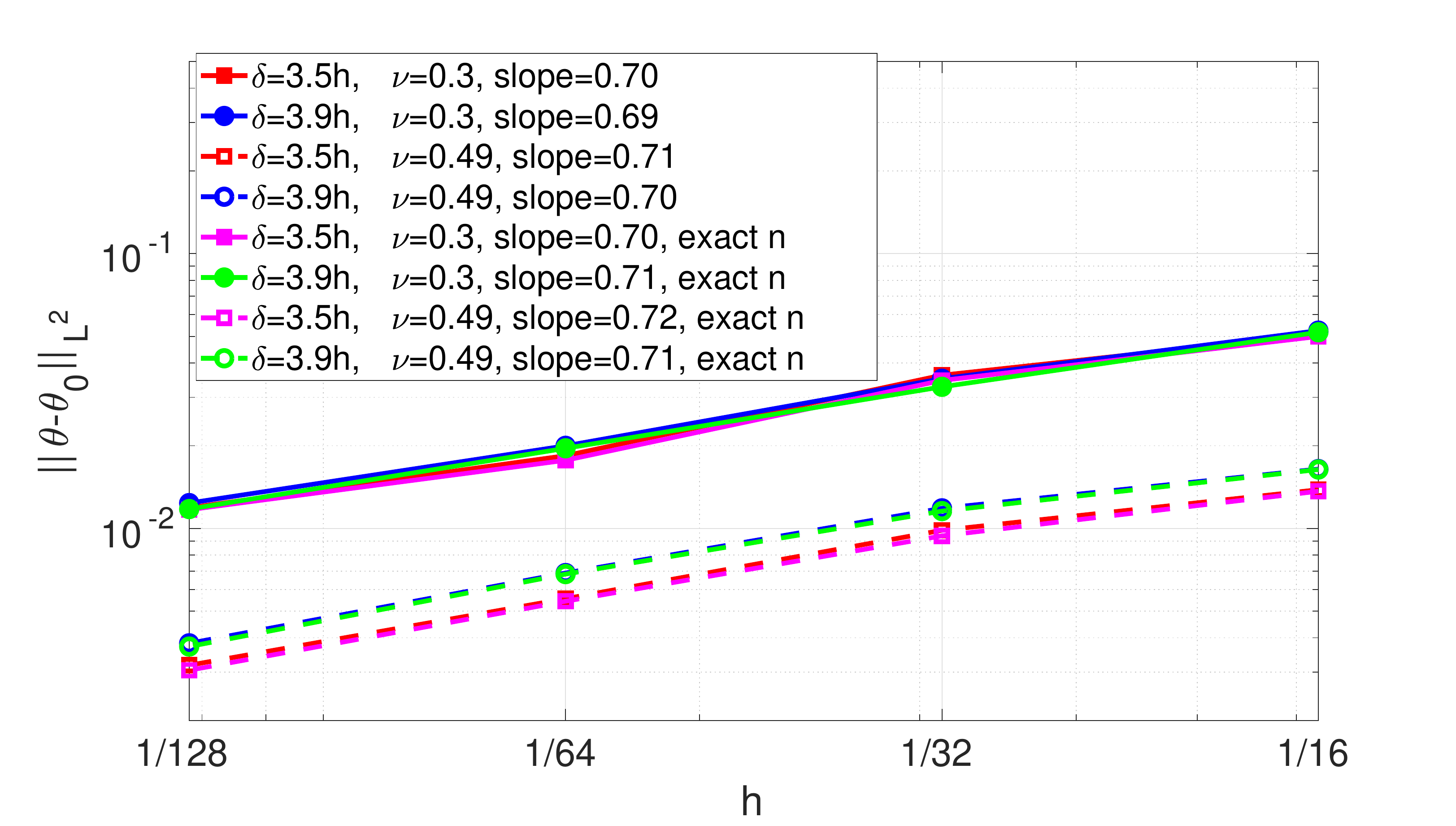}}
 \caption{M-convergence tests for a free-surface circular hole under remote loading with uniform discretizations. Left: the $L^2(\omg)$ difference between displacement $\ub$ and its local limit $\ub_0$. Right: the $L^2(\omg)$ difference between the nonlocal dilitation $\theta$ and its local limit $\theta_0=\nabla\cdot\ub_0$.
 }
 \label{fig:hole_s}
\end{figure}

 \begin{figure}[!htb]\centering
 \subfigure{\includegraphics[width=0.49\textwidth]{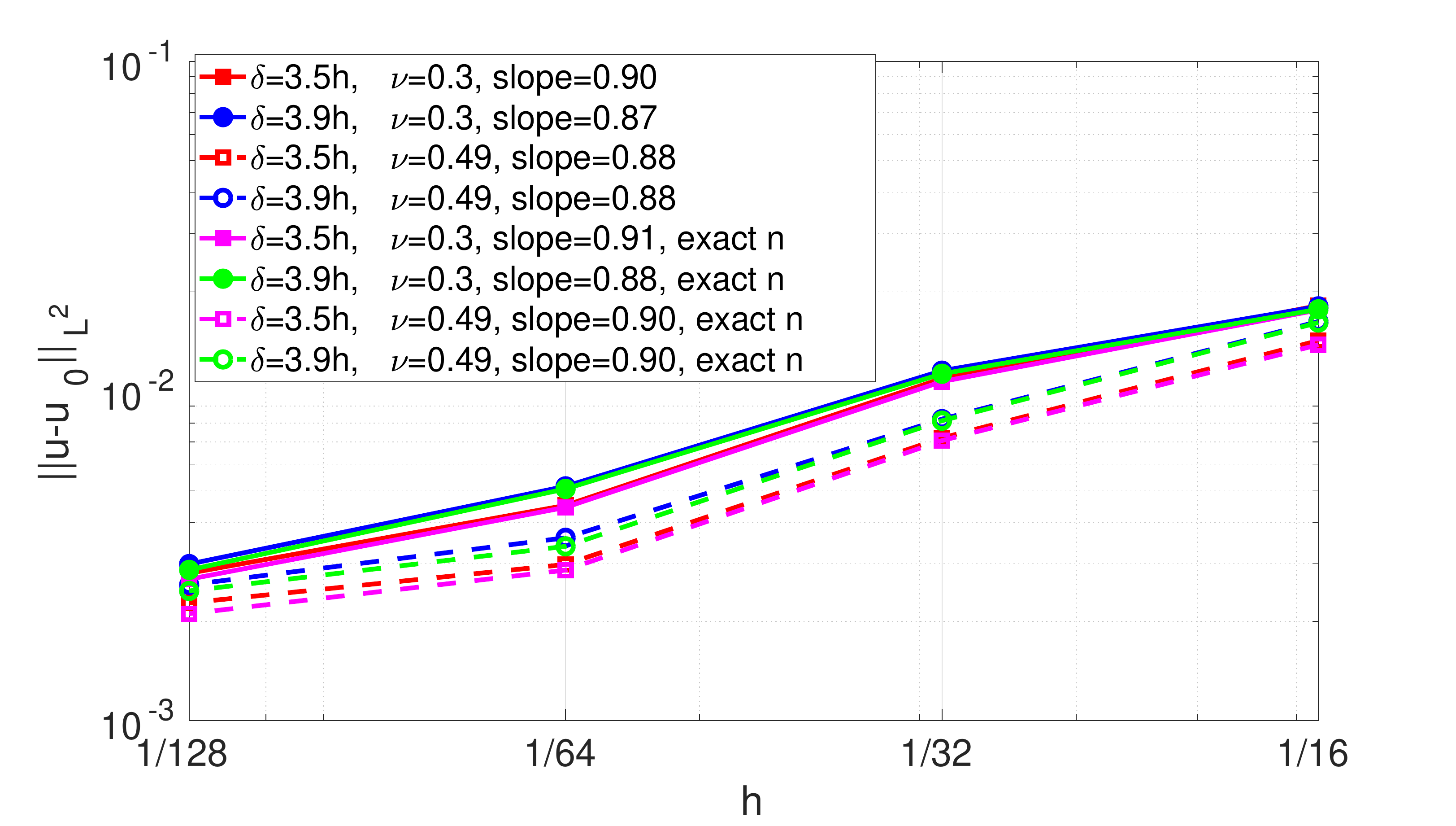}}
  \subfigure{\includegraphics[width=0.49\textwidth]{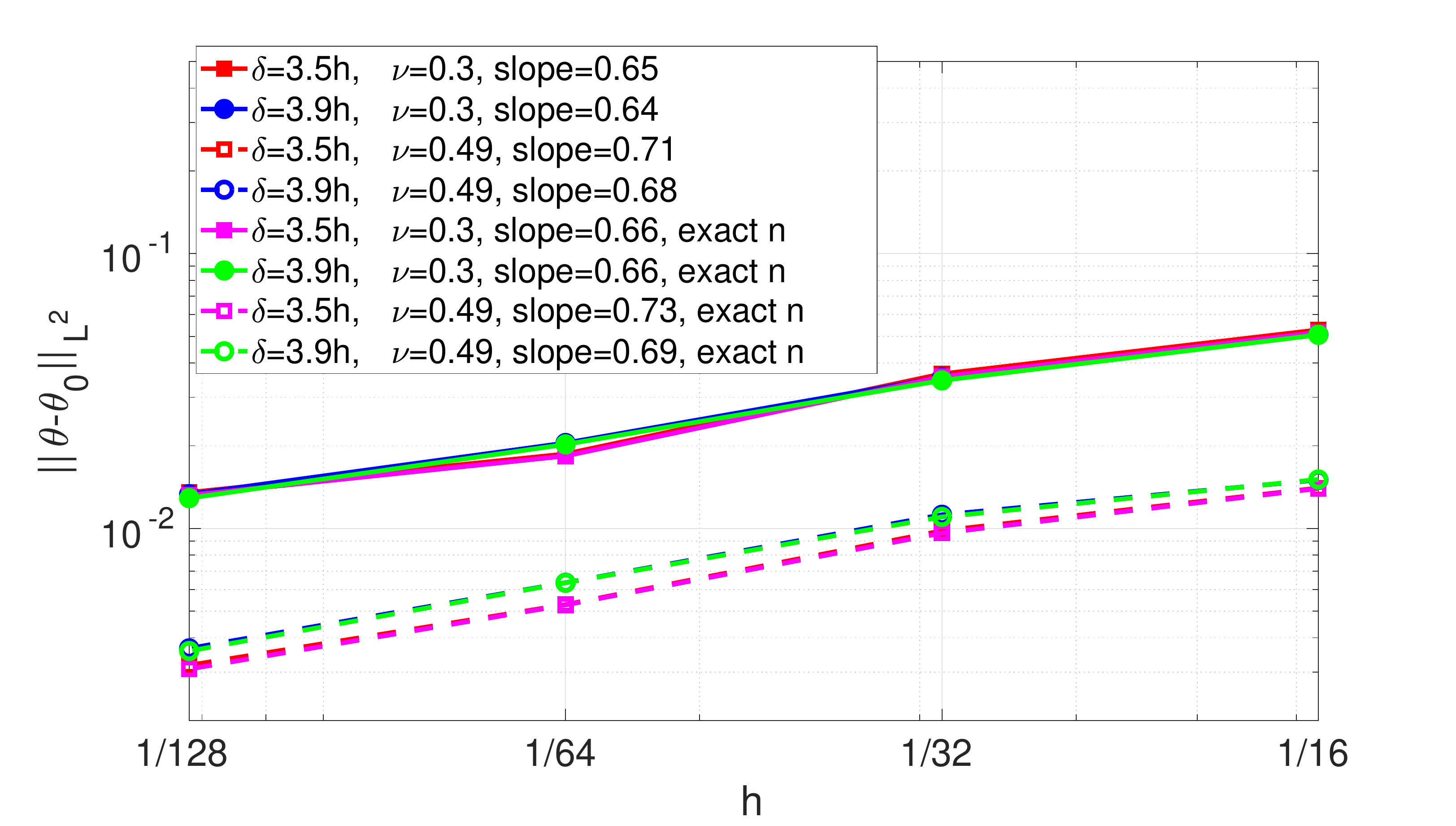}}
 \caption{M-convergence tests for a free-surface circular hole under remote loading with non-uniform discretizations. Left: the $L^2(\omg)$ difference between displacement $\ub$ and its local limit $\ub_0$. Right: the $L^2(\omg)$ difference between the nonlocal dilitation $\theta$ and its local limit $\theta_0=\nabla\cdot\ub_0$.
 }
 \label{fig:hole_un}
\end{figure}

 \begin{figure}[!htb]\centering
 \subfigure{\includegraphics[width=0.49\textwidth]{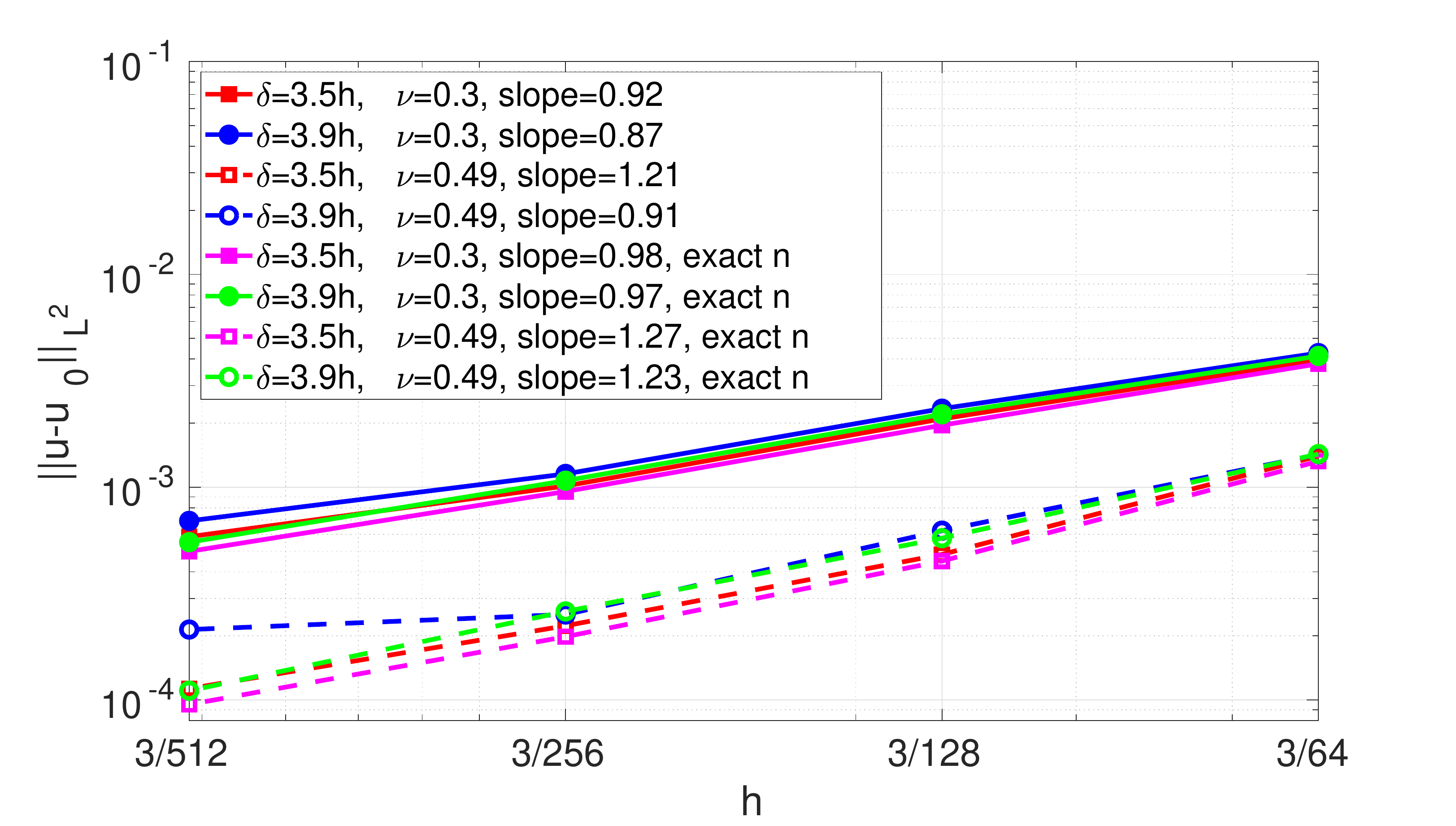}}
  \subfigure{\includegraphics[width=0.49\textwidth]{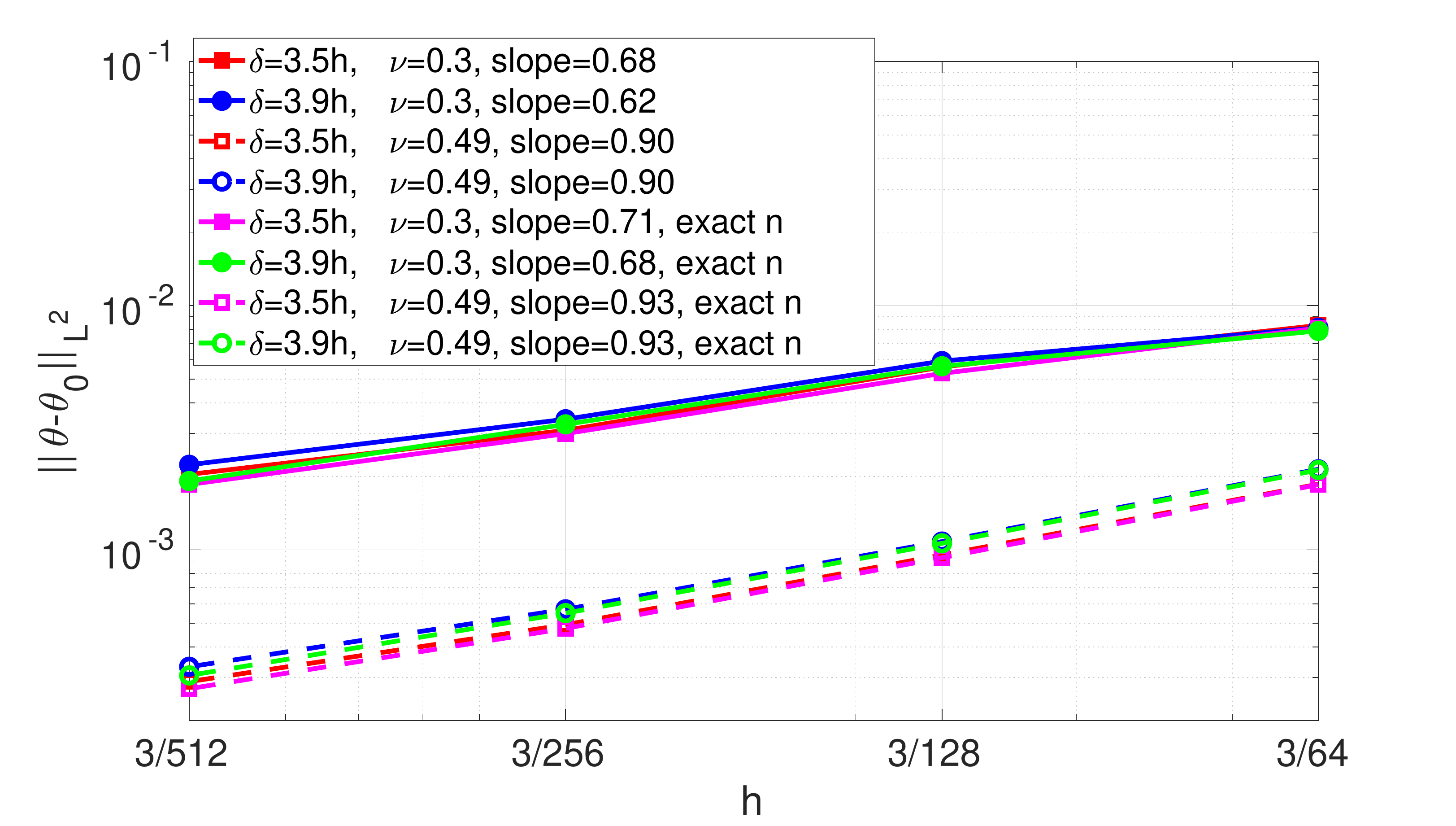}}
 \caption{M-convergence tests for a hollow disk under internal pressure with uniform discretizations. Left: the $L^2(\omg)$ difference between displacement $\ub$ and its local limit $\ub_0$. Right: the $L^2(\omg)$ difference between the nonlocal dilitation $\theta$ and its local limit $\theta_0=\nabla\cdot\ub_0$.
 }
 \label{fig:disk_s}
\end{figure}

 \begin{figure}[!htb]\centering
 \subfigure{\includegraphics[width=0.49\textwidth]{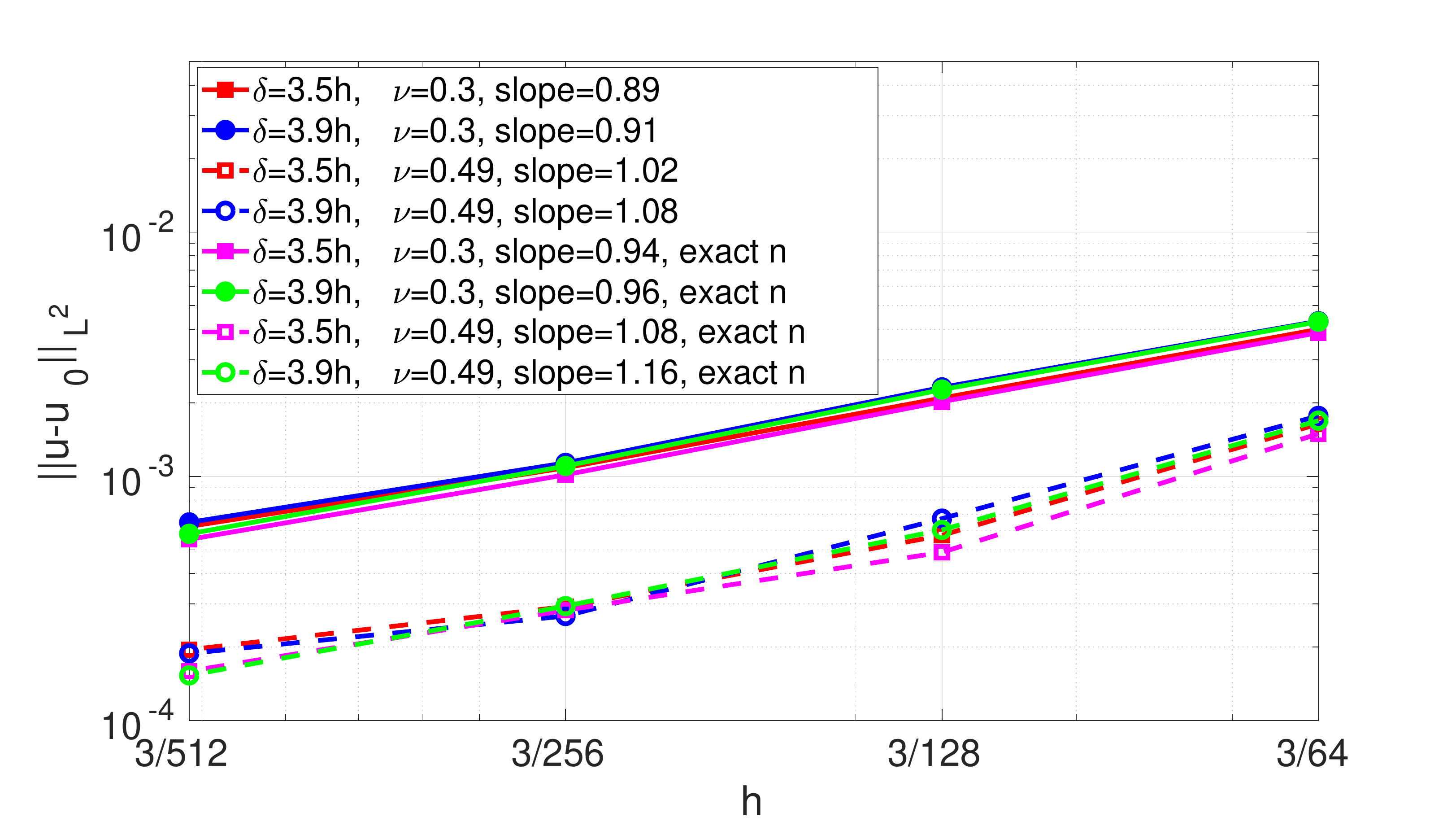}}
  \subfigure{\includegraphics[width=0.49\textwidth]{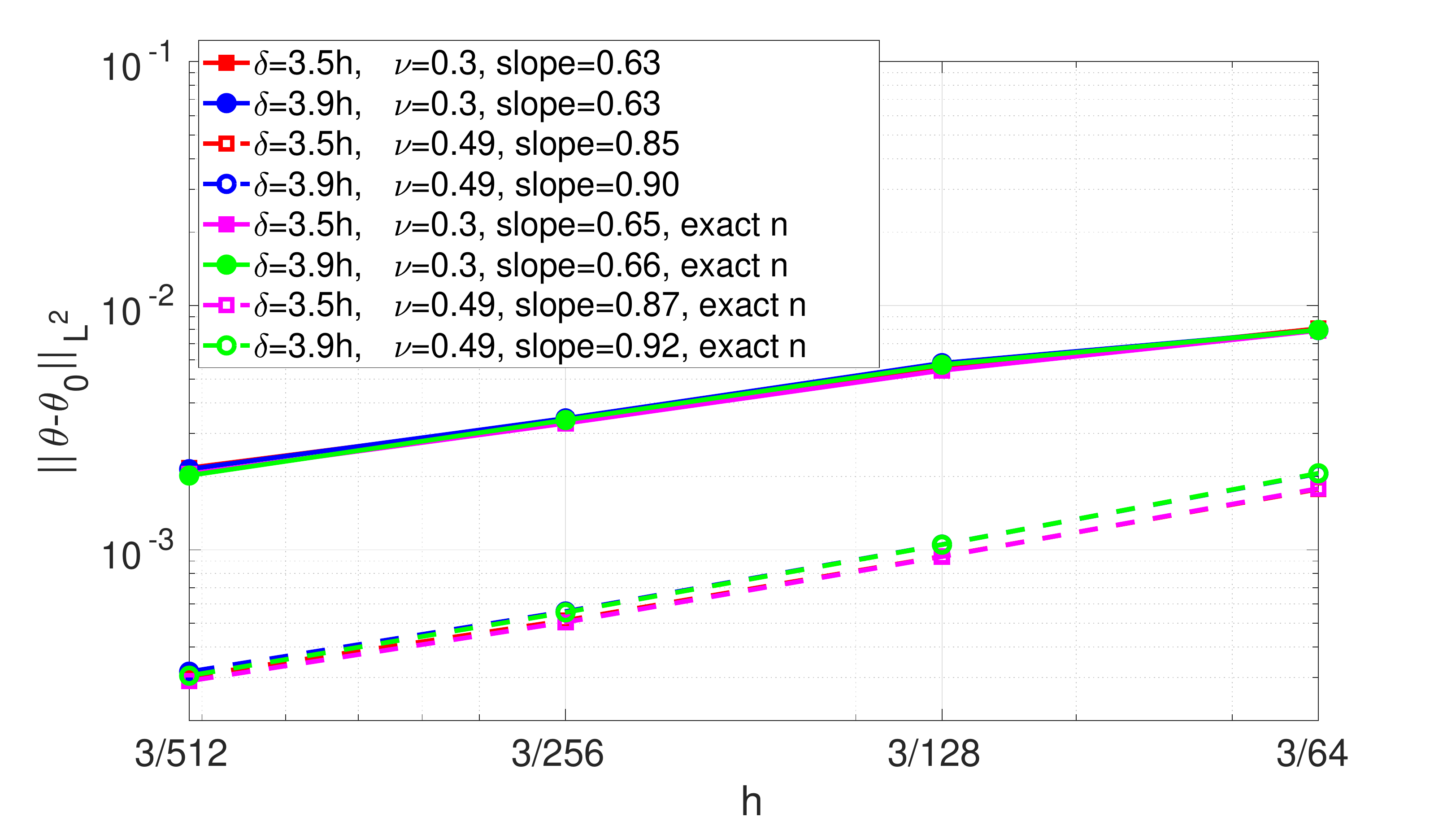}}
 \caption{M-convergence tests for a hollow disk under internal pressure with non-uniform discretizations. Left: the $L^2(\omg)$ difference between displacement $\ub$ and its local limit $\ub_0$. Right: the $L^2(\omg)$ difference between the nonlocal dilitation $\theta$ and its local limit $\theta_0=\nabla\cdot\ub_0$.
 }
 \label{fig:disk_un}
\end{figure}

\subsection{Composite materials with discontinuous material properties}

\begin{figure}[!htb]\centering
 \subfigure{\includegraphics[width=0.49\textwidth]{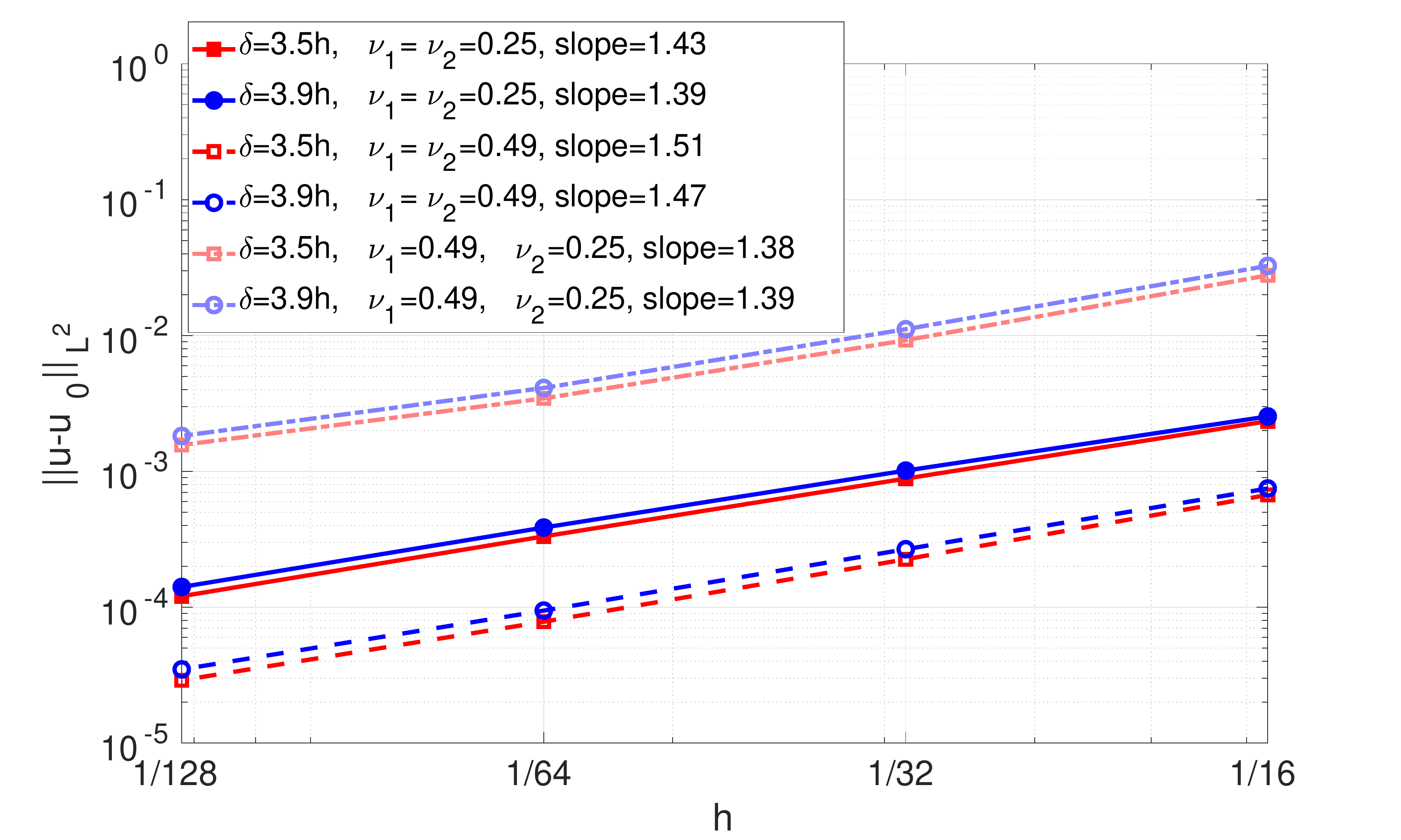}}
  \subfigure{\includegraphics[width=0.49\textwidth]{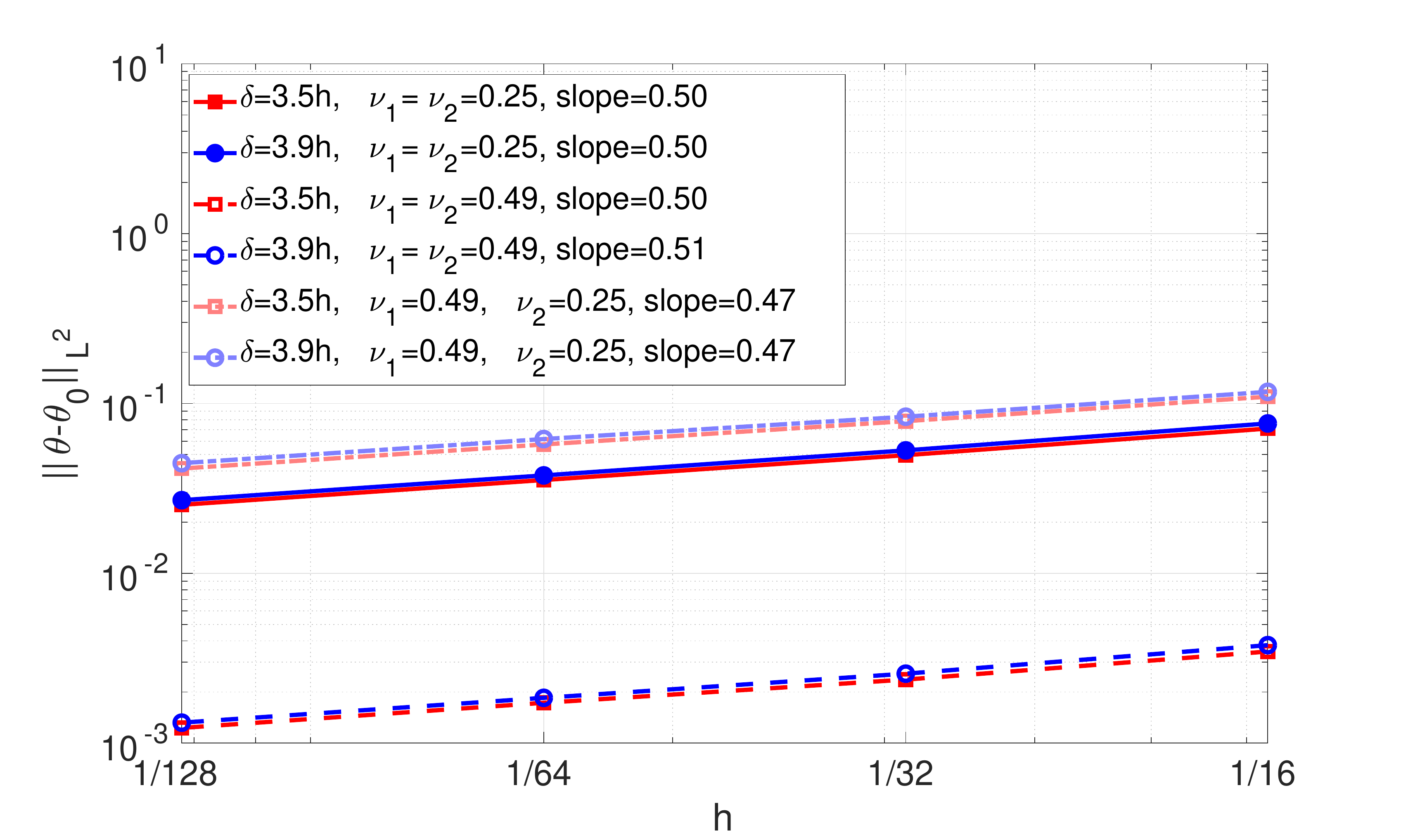}}
 \caption{M-convergence tests for composite materials with uniform discretizations. Left: the $L^2(\omg)$ difference between displacement $\ub$ and its local limit $\ub_0$. Right: the $L^2(\omg)$ difference between the nonlocal dilitation $\theta$ and its local limit $\theta_0=\nabla\cdot\ub_0$.
 }
 \label{fig:hetero_s}
\end{figure}

 \begin{figure}[!htb]\centering
 \subfigure{\includegraphics[width=0.49\textwidth]{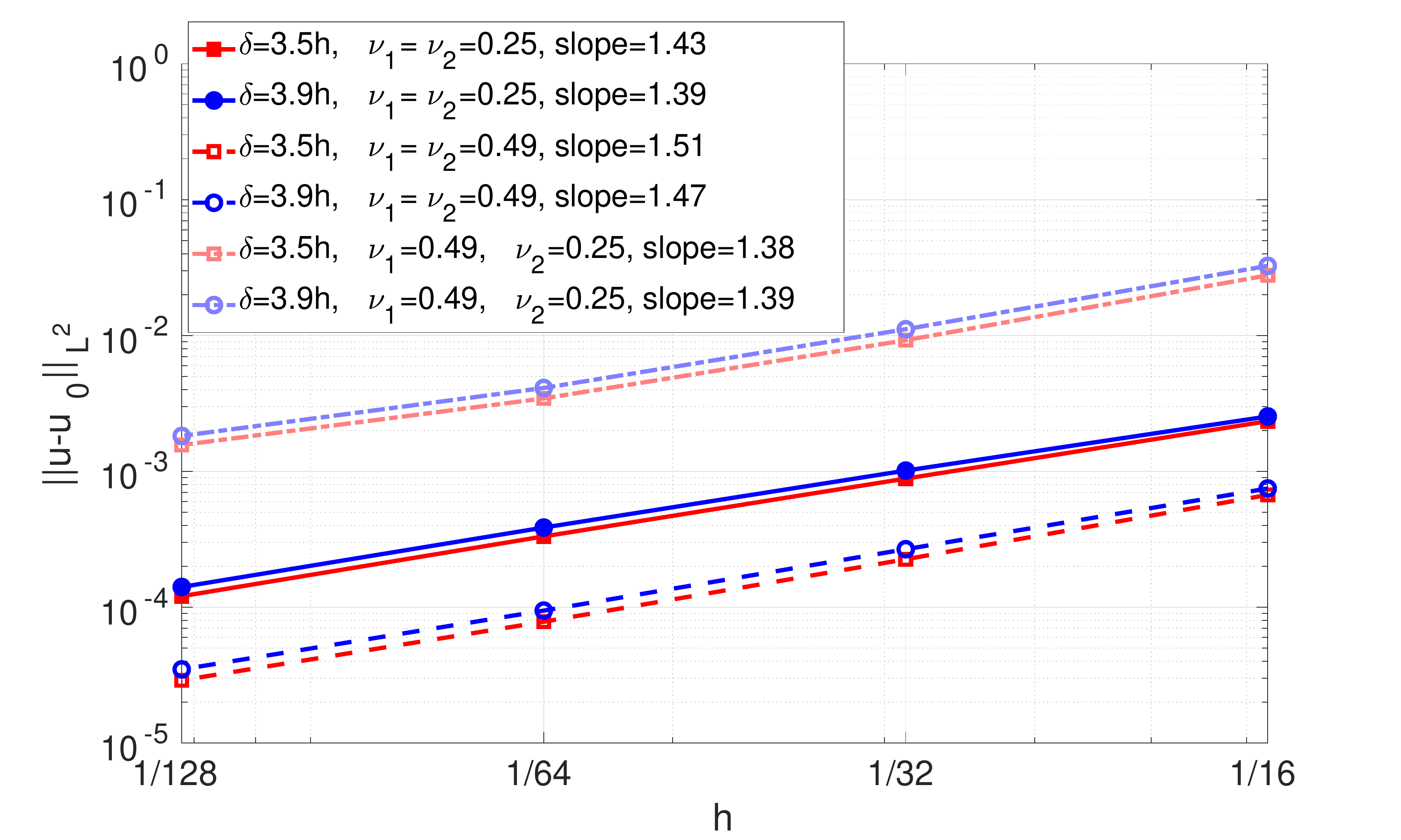}}
  \subfigure{\includegraphics[width=0.49\textwidth]{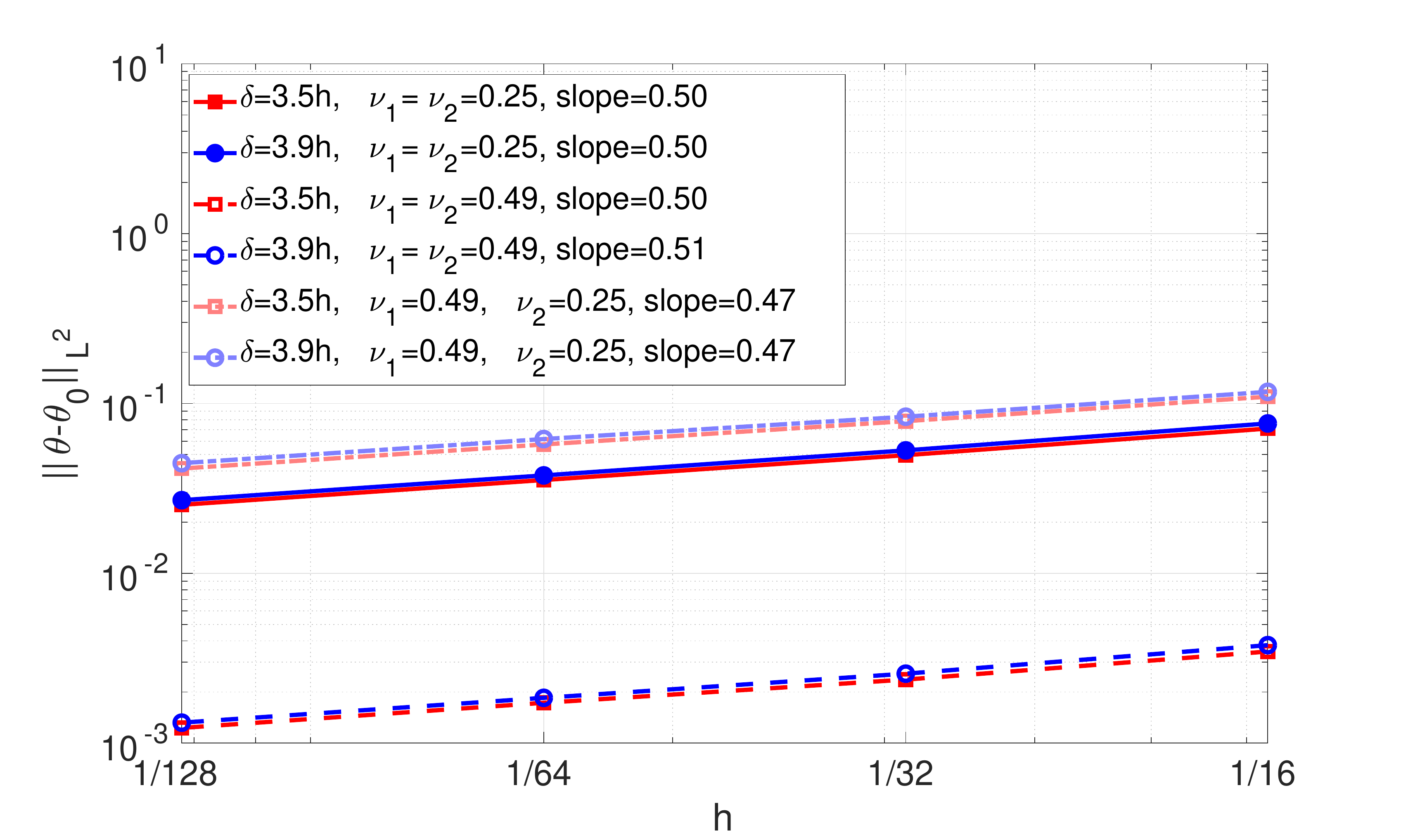}}
 \caption{M-convergence tests for composite materials with non-uniform discretizations. Left: the $L^2(\omg)$ difference between displacement $\ub$ and its local limit $\ub_0$. Right: the $L^2(\omg)$ difference between the nonlocal dilitation $\theta$ and its local limit $\theta_0=\nabla\cdot\ub_0$.
 }
 \label{fig:hetero_un}
\end{figure}

\pagebreak
\bibliography{arxiv}

\end{document}